\documentclass[11pt]{amsart}
\usepackage{geometry}                % See geometry.pdf to learn the layout options. There are lots.
\usepackage{graphicx}
\usepackage{amsmath,amssymb}
\usepackage{epstopdf}

\numberwithin{equation}{section}

\newtheorem{example}{Example}
\newtheorem{definition}{Definition}
\newtheorem{theorem}{Theorem}
\newtheorem{lemma}{Lemma}

\newtheorem{condition}{Condition}
\usepackage{hyperref}
\hypersetup{
    colorlinks,
    filecolor=black,
    linkcolor=black,
    urlcolor=black}
    \usepackage{subfigure}
\usepackage[lined,boxed]{algorithm2e}

\usepackage{float}
\usepackage{color}
\usepackage{lineno}
\modulolinenumbers[5]
\title{An artificial neural network approximation for Cauchy inverse problems}
 \author{Yixin Li, Xianliang Hu}
\address{School of Mathematical Science, Zhejiang University, $\textstyle{310027}$, Hangzhou, China}

\begin{document}
\maketitle
%\section{}
%\subsection{}

% \maketitle

\begin{abstract}
A novel artificial neural network method is proposed for solving Cauchy inverse problems. It allows multiple hidden layers with arbitrary width and depth, which theoretically yields better approximations to the inverse problems. In this research, the existence and convergence are shown to establish the well-posedness of neural network method for Cauchy inverse problems, and various numerical examples  are presented to illustrate its accuracy and stability. The numerical examples are from different points of view, including time-dependent and time-independent cases, high spatial dimension cases up to 8D, and cases with noisy boundary data and singular computational domain. Moreover, numerical results also show that neural networks with wider and deeper hidden layers could lead to better approximation for Cauchy inverse problems.

keywords:  Cauchy inverse problem; artificial neural network; well-posedness; high dimension; irregular domain; 
\end{abstract}

%\vskip .35cm   % for article
%{\bf Subject class:} 65M20, 65N22, 80A22

%%%%  1. First section
\section{Introduction}
\label{sec:introduction}
The approximation to Cauchy inverse problem is an important objective in the last few decades. Let us consider the following two classical cases: (I) for time independent problem, and (II) for time dependent problem:
Let $\Omega\subset\mathbb{R}^d$ be a  domain with continuous boundary $\partial \Omega$, where $d$ is the spatial dimension. It is worth to mention that $\Gamma$ is part but not all of  $\partial\Omega$, and the aim of Cauchy inverse problem is to recover solution $u$ on the rest of boundary $\partial\Omega/\Gamma$, with proper initial and boundary conditions $h,f,g$. 

\begin{equation}
\label{eq:elliptic}
(I) 
\begin{cases}
\mathcal{L} u(\mathbf{x}) = 0 & \mathbf{x} \ \ $in$  \ \ \Omega\\
u(\mathbf{x}) = f & \mathbf{x} \ \ $on$ \ \ \Gamma\\
\frac{\partial u(\mathbf{x})}{\partial \mathbf{n}} = g& \mathbf{x} \ \ $on$  \ \ \Gamma
\end{cases}
\end{equation}
and
\begin{equation}
\label{eq:parabolic}
(II) 
\begin{cases}
\frac{\partial u(\mathbf{x},t)}{\partial t} + \mathcal{L} u(\mathbf{x},t) = 0 & [\mathbf{x},t] \ \ $in$  \ \ \Omega \times\mathcal{T}\\
u(\mathbf{x},t) = f & [\mathbf{x},t] \ \ $on$ \ \ \Gamma\times\mathcal{T}\\
\frac{\partial u(\mathbf{x},t)}{\partial \mathbf{n}} = g& [\mathbf{x},t] \ \ $on$  \ \ \Gamma\times\mathcal{T} \\
u(\mathbf{x},0) = h& \mathbf{x}\ \ $in$ \ \ \Omega
\end{cases}
\end{equation}
where $\mathbf{n}$ is the outer unit normal with respect to $\partial\Omega$, $\mathcal{L}$ is a linear operator and $\mathcal{T} = [0,T]$ represents time. 

There are many works concerning the implementation and analysis for their numerical methods. For time-independent case \eqref{eq:elliptic}, Carleman-type estimates for the discrete scheme are used for Laplace's case in \cite{Santosa1991}. After that, many authors propose various algorithms for Cauchy inverse problem for Laplace's equation, such as conjugate gradient method\cite{Lesnic2000}, Backus-Gilbert algorithm\cite{Hon_2001}, regularization methods\cite{Reinhardt1999}, and some methods from linear algebra\cite{Nachaoui2002,Nachaoui2004}. Meanwhile some convergence and stability analysis are constructed. Chakib et al.\cite{Chakib2001} proved the existence of a solution to Cauchy problem, and it is first proved that the desired solution is the unique fixed point of some appropriate operator in \cite{Cimeti_re_2001}. For time-dependent case \eqref{eq:parabolic}, Besela et al. \cite{Besala1966} proved the uniqueness of solutions with direct case in 1966. After that there are many authors considering various types of stable numerical algorithms in different fields such as heat equation\cite{Cannon1967,Elden1987}, Helmholtz equation\cite{Berntsson2014} and time-space fractional diffusion equation\cite{Sakamoto2011, SLi2018}. 

The key point of numerical methods for Cauchy inverse problem is the ways to treat the ill-posedness. It is well known that there exists at most one solution to the above two Cauchy problem. However, they are typically ill-posed, which means that a small change of the initial data may induce large changes of the solutions. It is referred to \cite{Isakov1998} and the reference there in for more details on this issue. Regularization method is an effective and general technique to deal with the ill-posedness of inverse problem \cite{Ito2014}. During the last few decades, different regularization methods have been proposed to solve various PDE inverse problem, such as Tikhonov regularization method \cite{Jin_2008, Tomoya2008}, boundary element method \cite{Cheng_2014,Marin2003}, variational method\cite{Jin2010} and dynamical regularization algorithm \cite{Zhang_2018}, etc. As for Cauchy problem, some numerical analysis and experiments on the regularization for different equations are proposed, such as Laplace's equation\cite{Bourgeois_2005, Wei2013}, elliptic equation\cite{Feng_2013},  Helmholtz equation\cite{Qin2010, Berntsson2017} and so on.

Artificial neural networks(ANN) methods for approximating physical models described by PDE systems, as well as other kinds of non-linear problems, have attracted significant interest. Lagaris et al. \cite{Lagaris1998, Lagaris2000} used this idea early in 1998 for low-dimensional solutions. Then the idea was extended in several follow-up works on various direct problems, including high order differential equations\cite{Malek2006} and partial differential equations\cite{Aarts2001}. The approaches of these methods essentially rely on the universal approximation property of ANN, which was proved in the pioneering work \cite{Cybenko1989} for one single hidden layer, and then extended and refined in \cite{Kurt1989,Kurt1991}. 
 
It is well known that deeper networks approximate things better in the field of deep learning. In this sense, deep neural network(DNN) is also popular to solve PDE, especially for high-dimensional PDE problems like Hamilton-Jacobi-Bellman equation \cite{Justin2018,Giuseppe2017}. Recently, physics informed neural network models \cite{Raissi2019} are developed to solve PDE by demonstrated approach on the nature and arrangement of the available data, which are effective for various problems, including fractional ADEs \cite{Pang2019}, stochastic problems \cite{Zhang2019} and so on. A method to solve unknown governing equations with DNN is proposed in \cite{Xiu2019}. Long et.al \cite{Dong2019} propose a new deep neural network, named PDE-Net 2.0, to discover (time-dependent) PDEs. 

There are some other interesting topics combining ANN to enhance the performance of traditional methods. Mishra\cite{Mishra2018} combined existing finite differential method with ANN and White et al.\cite{White2019} used neural network surrogate in topology optimization. Li et al. \cite{Li2018} recast the training in deep learning as a control problem which is allowed to formulate necessary optimality conditions in continuous time using the Pontryagin’s maximum principle (PMP). Moreover, Yan and Zhou\cite{Yan2019} propose an adaptive procedure to construct a multi-fidelity polynomial chaos surrogate model in inverse problems.

In this paper, we propose a novel numerical method for solving Cauchy inverse problems using artificial neural network. As the spatial dimension grows, the computational cost of ANN method grows not so quickly as the traditional numerical methods. Within the proposed approach, we use a neural network instead of a linear combination of Lagrangian basis functions to represent the solution of PDEs, and impose the PDE constraint and boundary conditions via a collocation type method. %The parameters of networks are trained to satisfy the different operators, initial and boundary conditions at randomly sampled input datas with stochastic iterative methods, e.g., ADAM algorithm.  Numerical examples show that the method is good to deal with the ill-posedness for problem \eqref{eq:elliptic} and \eqref{eq:parabolic}, especially in the form of data completion. 

The rest of this paper is organized as follows. In Section \ref{sec:model}, we describe the neural network model for solving PDEs with linear operators and some initial and boundary conditions. In Section \ref{sec:analysis}, the convergence theorems are discussed in details. We prove the denseness and m-denseness of a network with which ensure the approximation capabilities of multi-hidden layer networks. Then we prove a theorem about convergence of ANN to approach the Cauchy inverse problems. Numerical examples are presented in Section \ref{sec:numericalexample}. We use the physical model's information(operator or initial and boundary datas with noise) rather than any other exact or experiment solutions to train the neural networks. At last some conclusions are given in Section \ref{sec:conclusion}.

\section{Artificial neural network(ANN) method for Cauchy inverse problem}
\label{sec:model}
\begin{figure}[H]
	\centering
	\includegraphics[height=5cm]{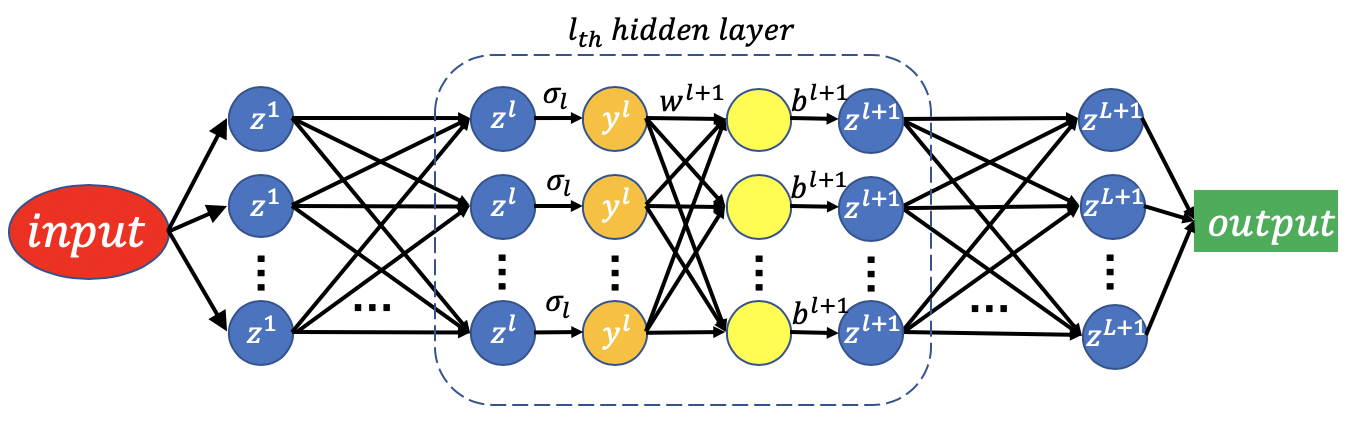}\qquad
	\caption{Structures of ANN with $L$ hidden layers}
	\label{fig:Neural_Networks}
	\end {figure}

Let us consider deep, fully connected feedforward ANNs to solve the Cauchy inverse problem. Given a network consisting of $L$ hidden layers. For convenience, the input and output layer are denoted as layer $0$ and layer $L+1$, respectively. There are some nonlinear functions being used in the hidden layers, says activation functions $\sigma$. The network defined above can mathematically be regarded as a mapping $\mathbb{R}^N \to \mathbb{R}$. Fig.~\ref{fig:Neural_Networks} shows the structure of such a network.
As it can be seen, in layer $l, l = 0, 1, \dots, L$, let $\mathbf{w}^l,\mathbf{b}^l$ denote the weights and bias and $\sigma_l$ be the activation functions. With above definitions, $\mathbf{z}^l$, the inputs of layer $l$ can be represented as
\begin{eqnarray*}
\label{eq:layerin}
\mathbf{z}^{l+1} &= &\mathbf{w}^{l+1}\mathbf{y}^{l} + \mathbf{b}^{l+1},\\
\mathbf{y}^{l+1} &=& \sigma_{l+1}(\mathbf{z}^{l+1})
\end{eqnarray*}
We use the notation $\mathbf{w} = \{\mathbf{w}^1, \mathbf{w}^2,\dots,\mathbf{w}^{L+1}\}, \mathbf{b} = \{\mathbf{b}^1,\mathbf{b}^2,\dots,\mathbf{b}^{L+1}\},$ inputs $\mathbf{x}$. For the simplicity, the notation of outputs is defined as
$$ \mathbf{y}^{L+1}:= NET(\mathbf{x};\mathbf{w},\mathbf{b}),$$ which is used to indicate that network takes $\mathbf{x} \in \mathbb{R}^N$ as input and parametrized by the weights and biases $\mathbf{w}$, $\mathbf{b}$.

\subsection{Network model for time-independent problem \eqref{eq:elliptic}}
The main idea of this method is to find a solution $\bar{u}$ for Problem \eqref{eq:elliptic} in the form of network output $NET(\mathbf{x};\mathbf{w},\mathbf{b})$. Defining the cost function
\begin{equation}
\label{eq:J_elliptic}
J(\bar{u}) = \Vert \mathcal{L}\bar{u}\Vert^2_{L_2(\Omega)} + \Vert \bar{u}-f\Vert^2_{L_2(\Gamma)} + \Vert \frac{\partial \bar{u}}{\partial \mathbf{n}} - g\Vert^2_{L_2(\Gamma)},
\end{equation}
then ANN approach for problem \eqref{eq:elliptic} can be written as
\begin{equation}
\label{eq:ANNpro_independent}
\begin{split}
&\displaystyle\min_{\mathbf{w},\mathbf{b}} J(\bar{u})\\
& s.t. \ \bar{u} = NET(\mathbf{x};\mathbf{w},\mathbf{b}).
\end{split}
\end{equation}
The equivalence of problem \eqref{eq:elliptic} and \eqref{eq:ANNpro_independent} will soon be shown in the next section. Here, let us first introduce the back propagation algorithm(gradient based method) to solve problem \eqref{eq:ANNpro_independent}. 
Denote $\mathbf{x}_{in} = [x_1,x_2,\dots, x_N]$ as random sampling in space $\Omega$, among which there are $N_o, N_d, N_n$ sampling points belonging to $\Omega, \Gamma$(Dirichlet boundary), $\Gamma$(Neumann boundary), respectively and it is required that $N_o + N_d + N_n = N$. For the purpose of verifying the stability of the approximation, certain statistical noise is added manually to the label data $f,g$, such that
\begin{equation*}
\Vert f^\delta - f\Vert_\Gamma \leq \delta, \ \ \Vert g^\delta - g\Vert_\Gamma\leq\delta,
\end{equation*}
where $\delta$ represents the level of statistical noise. For the ease of representation, the cost function \eqref{eq:J_elliptic} is written in discrete form as
\begin{equation}
\label{eq:J_discrete}
\begin{split}
J(\bar{u})  &= J_o(\bar{u}) + J_d(\bar{u}) + J_n(\bar{u})\\
&= \displaystyle\sum_{i=1}^{N_o}\left( \mathcal{L}\bar{u}_i\right)^2 + \sum_{i=1}^{N_d}\left(\bar{u}_i - f^\delta_i\right)^2 + \sum_{i = 1}^{N_n}\left(\frac{\partial \bar{u}_i}{\partial \mathbf{n}} - g_i^\delta\right)^2,
\end{split}
\end{equation}
where $\bar{u}_i = NET(\mathbf{x}_i;\mathbf{w},\mathbf{b}), f^\delta_i = f^\delta(\mathbf{x}_i)$ and $g^\delta_i = g^\delta(\mathbf{x}_i)$. To this point, the back propagation can be formulated as 
\begin{equation}
\label{eq:backpropa_ell}
\begin{split}
\frac{\partial J(\bar{u})}{\partial \mathbf{w}} &= \frac{\partial J_o(\bar{u})}{\partial \mathbf{w}} + \frac{\partial J_d(\bar{u})}{\partial \mathbf{w}} + \frac{\partial J_n(\bar{u})}{\partial \mathbf{w}} \\
&= 2\left(\sum_{i=1}^{N_o} \mathcal{L}\bar{u}_i \frac{\partial \mathcal{L} \bar{u}_i }{\partial\mathbf{w}} + \sum_{i = 1}^{N_d} \left(\bar{u}_i - f_i^\delta\right) \frac{\partial \bar{u}_i }{\partial\mathbf{w}} + \sum_{i=1}^{N_n} \left( \frac{\partial \bar{u}_i }{\partial \mathbf{n}}- g_i^\delta\right) \frac{\partial^2 \bar{u}_i }{\partial \mathbf{n}\partial\mathbf{w}}\right)
\end{split}
%\frac{\partial J(\bar{u})}{\partial \mathbf{b}} &= &\frac{\partial J_o(\bar{u})}{\partial \mathbf{b}} + \frac{\partial J_d(\bar{u})}{\partial \mathbf{b}} + \frac{\partial J_n(\bar{u})}{\partial \mathbf{b}}.
\end{equation}
Similarly, 
\begin{equation}
\begin{split}
\frac{\partial J(\bar{u})}{\partial \mathbf{b}} &= \frac{\partial J_o(\bar{u})}{\partial \mathbf{b}} + \frac{\partial J_d(\bar{u})}{\partial \mathbf{b}} + \frac{\partial J_n(\bar{u})}{\partial \mathbf{b}} \\
&= 2\left(\sum_{i=1}^{N_o} \mathcal{L}\bar{u}_i \frac{\partial \mathcal{L} \bar{u}_i }{\partial\mathbf{b}} + \sum_{i = 1}^{N_d} \left(\bar{u}_i - f_i^\delta\right) \frac{\partial \bar{u}_i }{\partial\mathbf{b}} + \sum_{i=1}^{N_n} \left( \frac{\partial \bar{u}_i }{\partial \mathbf{n}}- g_i^\delta\right) \frac{\partial^2 \bar{u}_i }{\partial \mathbf{n}\partial\mathbf{b}}\right)
\end{split}
%\frac{\partial J(\bar{u})}{\partial \mathbf{b}} &= &\frac{\partial J_o(\bar{u})}{\partial \mathbf{b}} + \frac{\partial J_d(\bar{u})}{\partial \mathbf{b}} + \frac{\partial J_n(\bar{u})}{\partial \mathbf{b}}.
\end{equation}
We supply the details to compute the $\mathcal{L} \bar{u},  \frac{\partial \bar{u} }{\partial \mathbf{n}}, \bar{u}$ and their corresponding back propagation in \ref{App:backpro}. 

To summarize, the structure of ANN method to solve time independent problem is shown in the following figure \ref{fig:ANN_schematic}
\begin{figure}[H]
\centering
\includegraphics[height=5cm]{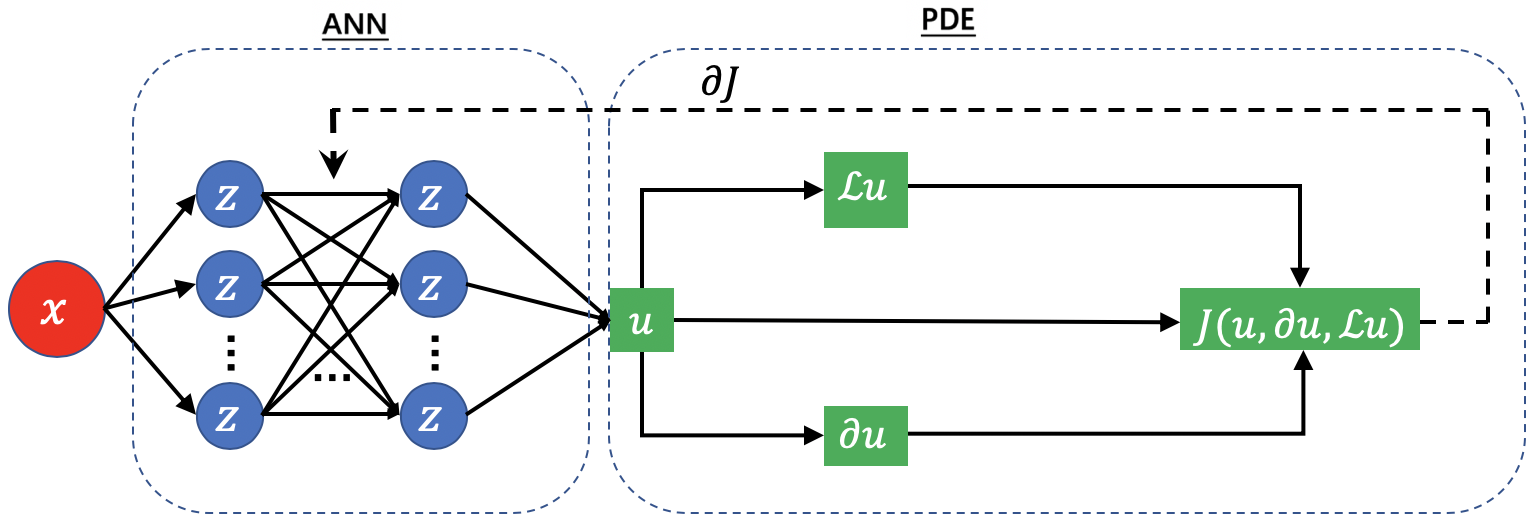}\qquad
\caption{Schematic of the ANN for solving Cauchy problem on time independent case}
\label{fig:ANN_schematic}
\end {figure}

\subsection{Network model for time-dependent problem}
The main idea of this method is to find a solution $\bar{u}$ for Problem \eqref{eq:parabolic} in the form of network output $NET(\mathbf{x},t;\mathbf{w},\mathbf{b})$. Defining the cost function
\begin{equation}
\label{eq:J_parabolic}
J(\bar{u}) = \Vert \mathcal{L}\frac{\partial \bar{u}}{\partial t} - \bar{u}\Vert^2_{L_2(\Omega\times\mathcal{T})} + \Vert \bar{u}-f\Vert^2_{L_2(\Gamma\times\mathcal{T})} + \Vert \frac{\partial \bar{u}}{\partial \mathbf{n}} - g\Vert^2_{L_2(\Gamma\times\mathcal{T}) }+ \Vert \bar{u} - h\Vert^2_{L_2(\Omega)},
\end{equation}
then ANN approach for problem \eqref{eq:parabolic} can be written as
\begin{equation}
\label{eq:ANNpro_dependent}
\begin{split}
& \displaystyle\min_{ \mathbf{w},\mathbf{b}} J(\bar{u}) \\
& \bar{u} = NET(\mathbf{x},t;\mathbf{w},\mathbf{b}).
\end{split}
\end{equation}
The equivalence of problem \eqref{eq:parabolic} and \eqref{eq:ANNpro_dependent} will soon be shown in the next section. Here, let us first introduce back propagation algorithm to solve problem \eqref{eq:ANNpro_dependent}. 
Denote $\mathbf{x}_{in} = [x_1,x_2,\dots, x_N]$ as random sampling in space $\Omega\times\mathcal{T}$, in which there are $N_o, N_d, N_n, N_t$ sampling points belonging to $\Omega\times\mathcal{T}, \Gamma\times\mathcal{T}$(Dirichlet boundary), $\Gamma\times\mathcal{T}$(Neumann boundary), $\Omega$, respectively, and it is required that $N_o + N_d + N_n + N_t= N$. For the purpose of verifying the stability of the approximation, certain statistical noise is added manually to the label data $f,g,h$, such that
\begin{equation*}
\Vert f^\delta - f\Vert_\Gamma \leq \delta, \ \ \Vert g^\delta - g\Vert_\Gamma\leq\delta, \ \ \Vert h^\delta - h\Vert_{\Omega} \leq\delta.
\end{equation*}
where $\delta$ represents the level of statistical noise. For the ease of representation, the cost function \eqref{eq:J_parabolic} is written in discrete form as
\begin{equation}
\label{eq:J_discrete_para}
\begin{split}
J(u)  &= J_o(u) + J_d(u) + J_n(u) + J_t(u)\\
&= \displaystyle\sum_{i=1}^{N_o} \left(\frac{\partial \bar{u}_i}{\partial t} + \mathcal{L}\bar{u}_i \right)^2 + \displaystyle\sum_{i=1}^{N
_d} \left(\bar{u}_i - f^\delta_i\right)^2 + \sum_{i = 1}^{N_n} \left( \frac{\partial \bar{u}_i}{\partial\mathbf{n}} - g^\delta_i\right)^2+  \displaystyle\sum_{i=1}^{N_t}\left(\bar{u}_i - h^\delta_i\right)^2 ,
\end{split}
\end{equation}
where $\bar{u}_i = NET(\mathbf{x}_i,t_i;\mathbf{w},\mathbf{b}), f^\delta_i = f^\delta(\mathbf{x}_i), g^\delta_i = g^\delta(\mathbf{x}_i)$ and $h^\delta_i = h^\delta(\mathbf{x}_i)$. To this point, the back propagation can be formulated as 
\begin{equation}
\label{eq:backpropa_para}
\begin{split}
\frac{\partial J(\bar{u})}{\partial \mathbf{w}} &= \frac{\partial J_o(\bar{u})}{\partial \mathbf{w}} + \frac{\partial J_d(\bar{u})}{\partial \mathbf{w}} + \frac{\partial J_n(\bar{u})}{\partial \mathbf{w}} +  \frac{\partial J_t(\bar{u})}{\partial \mathbf{w}}\\
&=  2\big(\sum_{i=1}^{N_o} \left(\frac{\partial \bar{u}_i}{\partial t}+\mathcal{L}\bar{u}_i \right) \left(\frac{\partial^2 \bar{u}_i}{\partial t\partial\mathbf{w}}+\frac{\partial \mathcal{L}\bar{u}_i }{\partial\mathbf{w}}\right) \\
&+ \sum_{i=1}^{N_n} \left( \frac{\partial \bar{u}_i }{\partial \mathbf{n}}- g_i^\delta\right)   \frac{\partial^2 \bar{u}_i }{\partial \mathbf{n}\partial\mathbf{w}} + \sum_{i = 1}^{N_d} \left( \bar{u}_i - f_i^\delta\right) \frac{\partial \bar{u}_i }{\partial\mathbf{w}} + \sum_{i = 1}^{N_t} \left( \bar{u}_i - h_i^\delta\right) \frac{\partial \bar{u}_i }{\partial\mathbf{w}} \big)
\end{split}
\end{equation}
Similarly, 
\begin{equation}
\begin{split}
\frac{\partial J(\bar{u})}{\partial \mathbf{b}} &= \frac{\partial J_o(\bar{u})}{\partial \mathbf{b}} + \frac{\partial J_d(\bar{u})}{\partial \mathbf{b}} + \frac{\partial J_n(\bar{u})}{\partial \mathbf{b}}  + \frac{\partial J_t(\bar{u})}{\partial \mathbf{b}},\\
&= 2\big(\sum_{i=1}^{N_o} \left(\frac{\partial \bar{u}_i}{\partial t}+\mathcal{L}\bar{u}_i \right) \left(\frac{\partial^2 \bar{u}_i}{\partial t\partial\mathbf{b}}+\frac{\partial \mathcal{L}\bar{u}_i }{\partial\mathbf{b}}\right) + \sum_{i=1}^{N_n} \left( \frac{\partial \bar{u}_i }{\partial \mathbf{n}}- g_i^\delta\right) \frac{\partial^2 \bar{u}_i }{\partial \mathbf{n}\partial\mathbf{b}}\\
& + \sum_{i = 1}^{N_d} \left( \bar{u}_i - f_i^\delta\right) \frac{\partial \bar{u}_i }{\partial\mathbf{b}} + \sum_{i = 1}^{N_t} \left( \bar{u}_i - h_i^\delta\right) \frac{\partial \bar{u}_i }{\partial\mathbf{b}} \big)
\end{split}
\end{equation}

We supply the details to compute the $\frac{\partial \bar{u}}{\partial t}, \mathcal{L} \bar{u},  \frac{\partial \bar{u}}{\partial \mathbf{n}}, \bar{u}$ and their corresponding back propagation in \ref{App:backpro}. To summarize, the structure of ANN method to solve time-independent problem is shown in the following figure \ref{fig:ANN_para}
\begin{figure}[H]
\centering
\includegraphics[height=5cm]{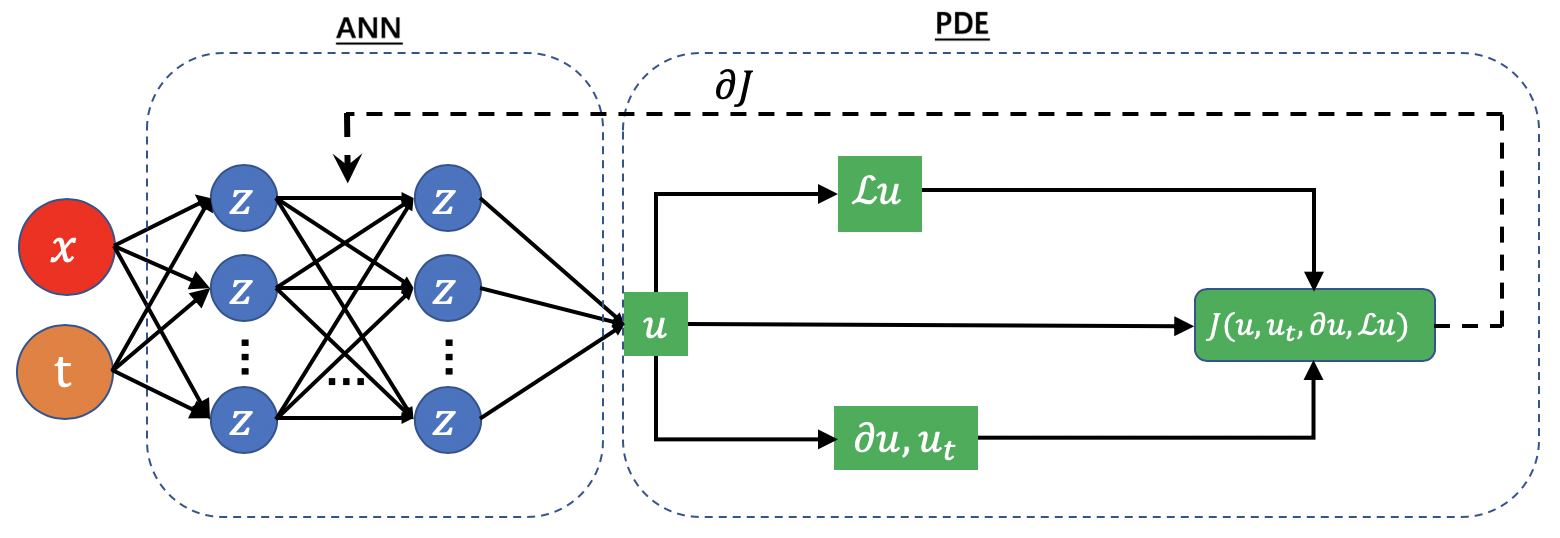}\qquad
\caption{Schematic of the ANN for solving Cauchy problem of time dependent equation}
\label{fig:ANN_para}
\end {figure}

\subsection{Training algorithm for networks}
\label{sec:algo}
The original algorithm for neural networks is gradient decent(GD) method. In this approximation it can be formulated as:

\begin{eqnarray*}
\mathbf{w}^{n+1} &= &\mathbf{w}^{n} - \Delta t \frac{\partial J(\bar{u})}{\partial \mathbf{w}},\\
\mathbf{b}^{n+1} &= &\mathbf{b}^{n} - \Delta t \frac{\partial J(\bar{u})}{\partial \mathbf{b}},
\end{eqnarray*}
where $n$ is the iterations and $\Delta t$ is time step. It is well know that ADAM algorithm is a stable and fast stochastic algorithm in the field of optimization. In this sense, ADAM algorithm is used in this research, and we would like to remark here that GD method can not reach a satisfying result in our numerical experiments. The main formulas for weights $\mathbf{w}$ are shown as following, and formula for bias $\mathbf{b}$ is similar to it.

\begin{equation}
\label{eq:ADAM}
\begin{cases}
\mathbf{w}^{n+1} = \mathbf{w}^{n} -  \Delta t\frac{v_w^{n+1}}{\sqrt{s_w^{n+1}} + \epsilon} \frac{\partial J(\bar{u})}{\partial \mathbf{w}}\\
v^{n+1}_w  = \left(\beta_1v^n_w + (1-\beta_1) \frac{\partial J(\bar{u})}{\partial \mathbf{w}^{n}}\right)/\left(1-\beta_1^{n}\right)\\
s^{n+1}_w = \left(\beta_2s^n_w + (1-\beta_2)\left(\frac{\partial J(\bar{u})}{\partial \mathbf{w}^{n}}\right)^2\right)/\left(1-\beta_2^{n} \right)
\end{cases}
\end{equation}
where $v_w,s_w$ is the matrix of parameters. $\beta_1$ and $\beta_2$ are constant closed to $1$ and $\epsilon$ is a small constant. To summarize, the ANN algorithm to solve the Cauchy problem is constructed as follows:

\begin{algorithm}[H]
\label{algorithm}
 \SetAlgoLined
 \textit{Input}: input data $\mathbf{x}^i \in \Omega$ or $\Omega\times\mathcal{T}$\; 
  \textit{Input}: target data $f^\delta, g^\delta, h^\delta$\;
  \textit{Input}: the number of hidden layers $L$\;
 Initialize the structure of neural networks\;
 Initialize weights $\mathbf{w}^l$, bias $\mathbf{b}^l$ and other parameters\;
  \While{$Iter\le np$}
 {Calculate outputs $\mathbf{y}^{L+1}$ with networks\;
 Calculate the cost function $J$ given by \eqref{eq:J_elliptic} or \eqref{eq:J_parabolic}\;
 Calculate gradient by \eqref{eq:backpropa_ell} or \eqref{eq:backpropa_para}\;
  Update $\mathbf{w}^l$ and $\mathbf{b}^l$ by ADAM algorithm \eqref{eq:ADAM}\;
 }
  \textit{Output}: the history of cost function $J$\;
  \textit{Output}: solution $\bar{u}$ of PDEs\;
\caption{ANN algorithm to solve the Cauchy problem}
\end{algorithm}

\section{Convergence of the neural network approximation}
\label{sec:analysis}
In this section, we discuss some conclusions on equivalence between PDE problem \eqref{eq:parabolic} and optimization problem \eqref{eq:ANNpro_dependent}. To fulfill this, the definitions of dense and m-dense networks  following \cite{Kurt1991} are necessary.
\begin{definition}[denseness]
A network $Net(\mathbf{x};\mathbf{w},\mathbf{b})$ is dense, if it satisfies
\begin{equation}
\Vert Net(\mathbf{x};\mathbf{w},\mathbf{b}) - f(\mathbf{x})\Vert\leq \epsilon, \ \ \forall f \in C(\bar{\Omega}),
\end{equation}
\end{definition}
\begin{definition}[m-denseness]
A network $Net(\mathbf{x};\mathbf{w},\mathbf{b})$ is m-dense, if it satisfies
\begin{equation}
\displaystyle\max_{\vert\alpha\vert\leq m}\Vert \nabla^\alpha Net(\mathbf{x};\mathbf{w},\mathbf{b}) - \nabla^\alpha f(\mathbf{x})\Vert \leq\epsilon.
\end{equation}
\end{definition}

The proof is carried out in two steps. In the first step(Section \ref{subsec:dense} and Section \ref{subsec:mdense}), we show that networks with $L$ hidden layers are dense and m-dense in design domain $\Omega$. In the second step(Section \ref{subsec:equivalence}), equivalence between PDE problem \eqref{eq:parabolic} and optimization problem \eqref{eq:ANNpro_dependent} are given. It is worth to mention that all the proofs in Section \ref{subsec:dense} and Section \ref{subsec:mdense} only depend on properties of networks. Suppose networks with $L$ hidden layers can be regarded as a mapping like
\begin{equation*}
\mathcal{A}_d^n(\sigma) = \{\xi(\mathbf{x};t):\mathbb{R}^{d+1} \to \mathbb{R}|\xi(\mathbf{x};t) = \displaystyle\sum_{i = 1}^{n}\xi_i \sigma(\mathbf{w}_i^{L}z^{L}(\mathbf{x};t) + b_i^{L})\},
\end{equation*}
where $\sigma$ is the sigmoid function $d$ is the spatial dimensions and $n$ is the number of units in the $L_{th}$ layer. It is straightforward to define $ \mathcal{A}_d(\sigma) = \displaystyle\cup_{n = 1}^{\infty}\mathcal{A}_d^n(\sigma)$ for briefness. 
\subsection{The denseness of $\mathcal{A}_d(\sigma)$}
\label{subsec:dense}
Consider a bounded set $\Omega$ in $\mathbb{R}^d$ with boundary $\partial\Omega$.
Kurt Hornik has proved the denseness and m-denseness of single hidden layer networks in \cite{Kurt1991}, and it will be extend ed to multi hidden layers' type in theorem \ref{thm:0_dense}. Let us consider an important lemma at first.
\begin{lemma}
\label{lm:l_approach}
Define $l_{th}(l = 1,2,\dots,L+1)$ hidden layer neural network function as $$\mathbf{z}^l(\mathbf{x};t) = \mathbf{w}^l\sigma(\mathbf{z}^{l-1}(\mathbf{x};t)) + \mathbf{b}^l, \ \ \mathbf{z}^l \in \mathbb{R}^{n_l}$$ for all $\epsilon > 0,$ there exists $\mathbf{A}^l =  [a^l_1, a^l_2,\dots,a^l_d]^T \in \mathbb{R}^{d\times n_l},$ such that 
\begin{equation}
\label{eq:lm1}
 \Vert \mathbf{A}^l\sigma(\mathbf{z}^l(\mathbf{x};t)) - \mathbf{x}\Vert_{\bar{\Omega}} := \displaystyle\sup_{i = 1}^d \vert a_i^l\sigma(\mathbf{z}^l(\mathbf{x};t)) - x_i\vert < \epsilon 
 \end{equation}
\end{lemma}
$\mathbf{Proof. }$
Let us use the method of induction to verify this lemma. 

I. Verify that equation \eqref{eq:lm1} holds when $l = 1$.

 Following theorem 2 in \cite{Kurt1991}, it is clear that for any $\epsilon > 0$ and $x_i \in \Omega$, there exists $a_i^1 \in \mathbb{R}^{n_1}$ such that 
\begin{equation*}
\label{eq:l1}
\Vert \mathbf{A}^1\sigma(\mathbf{z}^1(\mathbf{x};t)) - \mathbf{x}\Vert_{\bar{\Omega}} = \vert a_i^1\sigma(\mathbf{z}^1(\mathbf{x};t)) - x_i\vert < \epsilon,
\end{equation*}
which verify equation \eqref{eq:lm1}.

II. Assume equation \eqref{eq:lm1} is true for $l=k$ to verify that it also holds when $l = k+1$.

Fix $\mathbf{A}^{k+1}$, since sigmoid function $\sigma$ satisfies the Lipschitz continuity, it yields that
\begin{equation*}
\begin{split}
\displaystyle\sup_i\vert a_i^{k+1}\sigma(\mathbf{z}^{k+1}(\mathbf{x};t)) - x_i\vert &= \sup_i\vert a_i^{k+1}\sigma(\mathbf{w}^{k+1}\sigma(\mathbf{z}^k(\mathbf{x};t) )+ \mathbf{b}^k) - x_i\vert \\
&= \sup_i\vert \sum_ja_{ij}^{k+1}\sigma(\mathbf{w}_j^{k+1}\sigma(\mathbf{z}^k(\mathbf{x};t) )+ \mathbf{b}_j^k) - x_i\vert \\
&\leq \sup_i\vert \sum_ja_{ij}^{k+1}(\sigma(\mathbf{w}_j^{k+1}\sigma(\mathbf{z}^k(\mathbf{x};t) )+ \mathbf{b}_j^k) -x_i)\vert \\
&+ \sup_i\vert \sum_j a^{k+1}_{ij} x_i- x_i\vert \\
&\leq \sup_i \vert \sum_j a^{k+1}_{ij}\epsilon\vert +  \sup_i\vert \sum_j a^{k+1}_{ij} x_i- x_i\vert\\
&\leq \epsilon, \ \ by \ choosing \sum_j a^{k+1}_{ij} = 1,
\end{split}
\end{equation*}
which completes the proof of lemma \ref{lm:l_approach}. \qed

With the above lemma, we can extend theorem in \cite{Kurt1991} into multi-hidden layers neural networks, sees in the following theorem:
\begin{theorem}
	\label{thm:0_dense}
	For sigmoid function $\sigma$, network $\mathcal{A}_d(\sigma)$ is dense in $C(\bar{\Omega}\times\mathcal{T})$.
\end{theorem}
$\mathbf{Proof. }$
According to theorem 1 in \cite{Kurt1991}, it follows that
\begin{equation}
\Vert A_1\sigma(\mathbf{w}^1\mathbf{x} + b^1) - f(\mathbf{x})\Vert \leq \epsilon, \ \ \forall f \in C(\bar{\Omega}\times\mathcal{T}). 
\end{equation}
It is obviously that
\begin{equation}
\begin{split}
\Vert\mathcal{A}_d(\sigma) - f(\mathbf{x})\Vert &=  \Vert A_{L+1}\sigma(\mathbf{z}^{L+1}(\mathbf{x};t)) - f(\mathbf{x})\Vert \\
& \leq \Vert A_1\sigma(\mathbf{w}^1\mathbf{x} + \mathbf{b}^1) - f(\mathbf{x})\Vert \\
&+ \Vert  A_{L+1}\sigma(\mathbf{z}^{L+1}(\mathbf{x};t)) - A_1\sigma(\mathbf{w}^1\mathbf{x} + \mathbf{b}^1)\Vert\\
&\leq 2\epsilon, \ \ \forall f \in C(\bar{\Omega}\times\mathcal{T}).
\end{split}
\end{equation}
Hence the statements in theorem \ref{thm:0_dense} are proved. \qed

\subsection{The m-denseness of $ \mathcal{A}_d(\sigma)$}
\label{subsec:mdense}
Let us consider an important lemma at first.
\begin{lemma}
\label{lm:l_approach_m}
Define $l_{th}$ hidden layer neural network function as $$\mathbf{z}^l(\mathbf{x};t)  = \mathbf{w}^l\sigma(\mathbf{z}^{l-1}(\mathbf{x};t)) + \mathbf{b}^l, \ \ \mathbf{z}^l \in \mathbb{R}^{n_l}$$ then for all $ \epsilon > 0,$ there exists $\mathbf{A}^l =  [a^l_1, a^l_2,\dots,a^l_d]^T \in \mathbb{R}^{d\times n^l},$ such that 
\begin{equation}
\label{eq:lm2}
 \displaystyle\max_{\vert\alpha\vert\leq m}\sup_{\mathbf{x}\in\bar{\Omega}}\vert \nabla^\alpha\mathbf{A}^l\sigma(\mathbf{z}^l(\mathbf{x};t)) - \nabla^\alpha\mathbf{A}^1\sigma(\mathbf{w}^1\mathbf{x}+\mathbf{b}^1)\vert < \epsilon. 
 \end{equation}
\end{lemma}
$\mathbf{Proof. }$

Let us use the method of induction to verify this lemma. 

I. Verify that equation \eqref{eq:lm2} holds when $l = 1$.

 it is clear that for any $\epsilon > 0$ and $x_i \in \Omega$, there establish 
\begin{equation}
\label{eq:l1}
 \displaystyle\max_{\vert\alpha\vert\leq m}\sup_{\mathbf{x}\in\bar{\Omega}}\vert \nabla^\alpha\mathbf{A}^l\sigma(\mathbf{z}^l(\mathbf{x};t)) - \nabla^\alpha\mathbf{A}^1\sigma(\mathbf{w}^1\mathbf{x}+\mathbf{b}^1)\vert  = 0 < \epsilon,
 \end{equation}
which verify equation \eqref{eq:lm2}.

II. Assume equation \eqref{eq:lm2} is true for $l=k$ to verify that it also holds when $l = k+1$.

Fix $\mathbf{A}^{k+1}$, since sigmoid function $\sigma$ satisfies the Lipschitz continuity and is bounded, it yields that
\begin{equation*}
\begin{split}
 &\ \ \ \displaystyle\max_{\vert\alpha\vert\leq m}\sup_{\mathbf{x}\in\bar{\Omega}}\vert \nabla^\alpha\mathbf{A}^{k+1}\sigma(\mathbf{z}^{k+1}(\mathbf{x};t)) - \nabla^\alpha\mathbf{A}^1\sigma(\mathbf{w}^1\mathbf{x}+\mathbf{b}^1)\vert\\
 &=  \displaystyle\max_{\vert\alpha\vert\leq m}\sup_{\mathbf{x}\in\bar{\Omega}}\sup_i\vert \nabla^\alpha a_i^{k+1}\sigma(\mathbf{z}^l(\mathbf{x};t)) - \nabla^\alpha a_i^1\sigma(\mathbf{w}^1\mathbf{x}+\mathbf{b}^1)\vert\\
&= \displaystyle\max_{\vert\alpha\vert\leq m}\sup_{\mathbf{x}\in\bar{\Omega}}\sup_i\vert \nabla^\alpha(\sum_ja_{ij}^{k+1}\sigma(\mathbf{w}_j^{k+1}\sigma(\mathbf{z}^k(\mathbf{x};t) )+ \mathbf{b}_j^k) - \sum_j a^1_{ij}\sigma(\mathbf{w}^1_j\mathbf{x} + \mathbf{b}^1_j))\vert \\
&\leq \displaystyle\max_{\vert\alpha\vert\leq m}\sup_{\mathbf{x}\in\bar{\Omega}}\sup_i\nabla^\alpha\big(\vert \sum_ja_{ij}^{k+1}(\sigma(\mathbf{w}_j^{k+1}\sigma(\mathbf{z}^k(\mathbf{x};t) )+ \mathbf{b}_j^k) -\sigma(\mathbf{w}^1_j\mathbf{x}+ \mathbf{b}^1_j))\vert\\
& + \vert \sum_j (a^{k+1}_{ij}-a^1_{ij}) \sigma(\mathbf{w}^1_j\mathbf{x} + \mathbf{b}^1_j))\vert\big) \\
&\leq L\epsilon \ \ by \ choosing \sum_j a^{k+1}_{ij} = \sum_ja^1_{ij} =  1,
\end{split}
\end{equation*}
which completes the proof of lemma \ref{lm:l_approach_m}. \qed

With the above lemma, we can extend the theorem 3 in \cite{Kurt1991} into multi-hidden layers neural networks, which yields following theorem:

\begin{theorem}
	\label{thm:m_dense}
	For sigmoid function $\sigma \in C^m(\bar{\Omega}\times\mathcal{T})$, we have $\mathcal{A}_d(\sigma)$ is uniformly m-dense on $C^m(\bar{\Omega}\times\mathcal{T})$.
\end{theorem}
$\mathbf{Proof. }$
According to theorem 3 in \cite{Kurt1991}, it follows that
\begin{equation}
\displaystyle\max_{\vert\alpha\vert\leq m}\Vert \nabla^\alpha A_1\sigma(\mathbf{w}^1\mathbf{x} + b^1) - \nabla^\alpha f(\mathbf{x})\Vert \leq \epsilon, \forall f \in C^m(\bar{\Omega}). 
\end{equation}
It is obviously that
\begin{equation}
\begin{split}
&\displaystyle\max_{\vert\alpha\vert\leq m}\Vert \nabla^\alpha A_{L+1}\sigma(\mathbf{z}^{L+1}(\mathbf{x};t)) - \nabla^\alpha f(\mathbf{x})\Vert \\
& \leq \displaystyle\max_{\vert\alpha\vert\leq m}(\Vert \nabla^\alpha A_{L+1}\sigma(\mathbf{z}^{L+1}(\mathbf{x};t)) - \nabla^\alpha A_{1}\sigma(\mathbf{w}^{1}\mathbf{x} + b^{1}))\Vert \\
& + \Vert \nabla^\alpha A_1\sigma(\mathbf{w}^1\mathbf{x} + b^1) - \nabla^\alpha f(\mathbf{x})\Vert)\\
& \leq (L+1)\epsilon, \forall f \in C^m(\bar{\Omega}\times\mathcal{T}). 
\end{split}
\end{equation}
Hence the statements in theorem \ref{thm:m_dense} are proved. \qed

%\begin{proposition}
%\label{pr:unique}
%There exist a unique solution $u^p$ in to Problem \eqref{eq:Cauchy}, moreover $u^p \in H^2(\Omega)$.
%\end{proposition}
\subsection{Equivalence between PDE problem \eqref{eq:parabolic} and optimization problem \eqref{eq:ANNpro_dependent}}
\label{subsec:equivalence}
Let us assume initially that problem \eqref{eq:parabolic} owns the following conditions
\begin{condition}
\label{lm:unique}
\begin{itemize}
\item There exist a unique solution $u^p$ in Problem \eqref{eq:parabolic}, moreover $u^p \in H^2(\Omega\times\mathcal{T})$.
\item $\mathcal{L}$ is Lipschitz continuous on $\bar{\Omega}\times\mathcal{T}$.
\item $h \in C^{2}(\bar{\Omega}\times\mathcal{T}) $ and its first derivative bounded in $\bar{\Omega}\times\mathcal{T}$.
\item $\partial\Omega \in C^2$.
\end{itemize}
\end{condition}

It follows two important theorems according to these conditions :

\begin{theorem}
\label{thm:conv_E}
For all $ \epsilon > 0$, there exists a series of neural network approach $\psi^n(\mathbf{x};t;\mathbf{w},\mathbf{b})$ such that $J(\psi^n) < K\epsilon$, where $K = \max(\vert \Omega\times\mathcal{T}\vert, \vert \Gamma\times\mathcal{T}\vert,\vert\Omega\vert)$ and $J(\psi^n)$ is in the case of equation \eqref{eq:J_parabolic}.
\end{theorem}
$\mathbf{Proof. }$
Let us define $\psi\in\mathcal{A}_d(\sigma)$ as a neural network approach. Assume that operator $\mathcal{L} \in C^m(\mathbb{R}^d)$ and sigmoid function $\sigma_l \in C^{m}(\mathbb{R}^{d}\times\mathcal{T})$ is non-constant and bounded in $[0,1]$. It is clear that $\bar{\Omega}$ is a compact subset in $\mathbb{R}^d$. According to theorem \ref{thm:m_dense}, there follows that
\begin{equation}
\displaystyle\max_{|\alpha\leq m|}\sup_{\mathbf{x}\in\bar{\Omega}} \vert \nabla^\alpha u(\mathbf{x};t) - \nabla^\alpha \psi\vert < \frac{\epsilon}{4}, \ \ for \ all \ u\in C^m(\mathbb{R}^d\times\mathcal{T}),
\end{equation}
which yields
\begin{equation}
\label{eq:L1_control}
\begin{split}
\displaystyle\sup_{\mathbf{x}\in\Gamma,t\in\mathcal{T}}\left\vert u - \psi\right\vert + \sup_{\mathbf{x}\in\Gamma,t\in\mathcal{T}}\left\vert \frac{\partial u}{\partial\mathbf{n}} - \frac{\partial\psi}{\partial\mathbf{n}}\right\vert + \sup_{\mathbf{x}\in\Omega,t\in\mathcal{T}}\left\vert \frac{\partial u}{\partial t} + \mathcal{L}u  - (\frac{\partial\psi}{\partial t} + \mathcal{L}\psi)\right\vert + \sup_{\mathbf{x}\in\Omega,t=0}\left\vert u-\psi\right\vert \\
 < 4\sum_{\alpha = 0}^m\sup_{\mathbf{x}\in\bar{\Omega}} \vert \nabla^\alpha u(\mathbf{x};t) - \nabla^\alpha \psi\vert < \epsilon
 \end{split}
\end{equation}
Let $u^p(\mathbf{x})$ be a solution of Problem \eqref{eq:parabolic} and using the conclusion of equation \eqref{eq:L1_control} and H$\ddot{o}$lder inequality, there establishes 
\begin{equation}
\begin{split}
J(\psi) & = \left\Vert \frac{\partial \psi}{\partial t} + \mathcal{L}\psi - (\frac{\partial u^p}{\partial t}+\mathcal{L}u^p)\right\Vert^2_{L^2(\Omega\times\mathcal{T})} + \left\Vert \frac{\partial\psi}{\partial\mathbf{n}} - \frac{\partial u^p}{\partial\mathbf{n}}\right\Vert^2_{L^2(\Gamma\times\mathcal{T})}\\
& + \Vert \psi - u^p\Vert^2_{L^2(\Gamma\times\mathcal{T})} +  \Vert \psi - u^p\Vert^2_{L^2(\Omega\times(t=0))}\\
&\leq \vert \Gamma\times\mathcal{T}\vert \displaystyle\sup_{\mathbf{x}\in\Gamma,t\in\mathcal{T}}\left\vert u^p - \psi\right\vert^2 +  \vert \Omega\vert \displaystyle\sup_{\mathbf{x}\in\Omega,t=0}\left\vert u^p - \psi\right\vert^2\\
&+ \vert\Gamma\times\mathcal{T}\vert\sup_{\mathbf{x}\in\Gamma,t\in\mathcal{T}}\left\vert \frac{\partial u^p}{\partial\mathbf{n}} - \frac{\partial\psi}{\partial\mathbf{n}}\right\vert^2 + \vert\Omega\times\mathcal{T}\vert\sup_{\mathbf{x}\in\Omega,t\in\mathcal{T}}\left\vert \frac{\partial \psi}{\partial t} + \mathcal{L}\psi - (\frac{\partial u^p}{\partial t}+\mathcal{L}u^p)\right\vert^2 \\
&< K\epsilon^2 < K\epsilon \ (when \ \epsilon\to 0),
\end{split}
\end{equation}
which completes the proof of theorem \ref{thm:conv_E}. \qed

\begin{theorem}
\label{thm:conv_u}
The series of neural network approach $\{\psi^n\}$ converges to the solution $u^*$ of \eqref{eq:parabolic}, as $n\to\infty$.
\end{theorem}
$\mathbf{Proof. }$
Since $\psi^n$ is a solution of \eqref{eq:ANNpro_dependent}, it is clear that $J(\psi^n) \leq J(\psi)$ for all $\psi \in \mathcal{A}_d(\sigma)$. In particular, there establishes $0\leq J(\psi^n) \leq J(\psi) \leq K\epsilon$. When $\epsilon\to 0$, it involves that
\begin{eqnarray*}
(\frac{\partial}{\partial t} + \mathcal{L})\psi^n = g^n   & in  & \Omega\times\mathcal{T}, \\
\psi^n - f = g^n_f &on& \Gamma\times\mathcal{T},\\
\frac{\partial\psi^n}{\partial\mathbf{n}} - g = g_g^n &on& \Gamma\times\mathcal{T}\\
\psi^n(x;0) - h = g^n_h & in & \Omega,
\end{eqnarray*}
for some $g^n, g^n_f, g^n_g, g_h^n$ such that
\begin{equation}
\Vert g^n\Vert^2_{L^2(\Omega\times\mathcal{T})} + \Vert g^n_f\Vert^2_{L^2(\Gamma\times\mathcal{T})} + \Vert g^n_g\Vert^2_{L^2(\Gamma\times\mathcal{T})} + \Vert g^n_h\Vert^2_{L^2(\Omega\times(t=0))} \to 0,\ \ as \ \ n\to\infty
\end{equation}
Assume condition \ref{lm:unique} establishes and $\psi^n \in L^2(\Omega\times\mathcal{T})$. It is clear that $\{\psi^n\}$ is uniformly bounded with respect to $n$ in $L^\infty(\mathcal{T},L^2(\Omega))$, which imply that there exist a subsequence, denoted by $\{\hat{\psi}^n\}$, converging to some $u$ in the weak-* sense in $L^\infty(\mathcal{T},L^2(\Omega))$. 

Next following condition \ref{lm:unique} and theorem 7.3 in \cite{Justin2018}, it is obviously that there exists a constant $C<\infty$ such that $$ \int_{\bar{\Omega}\times\mathcal{T}} \left(\frac{\partial}{\partial t}\mathcal{L}\right)\hat{\psi}^nd\mathbf{x}dt < C,$$ which lead that $\{\bar{\psi}^n\}$ converges almost everywhere to $u$ in $\bar{\Omega}\times\mathcal{T}.$ Then it can be proved that $\psi^n$ is the solution of problem \eqref{eq:parabolic} when $n\to\infty$, which completes the proof of theorem \ref{thm:conv_u}. \qed

To this end, we have proved that problem \eqref{eq:parabolic} is equivalent to problem \eqref{eq:ANNpro_dependent}(as well as problem \eqref{eq:elliptic} and \eqref{eq:ANNpro_independent}) if condition \ref{lm:unique} holds. In addition, the neural network solution is convergence to the exact solution.

\section{Numerical examples}
\label{sec:numericalexample}
In this section, we present extensive numerical results to demonstrate the ANN method for Cauchy problem. Firstly examples of low and high dimensional problems are displayed to verify the accuracy of this method both on time- dependent and independent cases.
\subsection{Numerical validation of time-dependent case}
\begin{figure}[H]
\centering
\includegraphics[height=5.5cm]{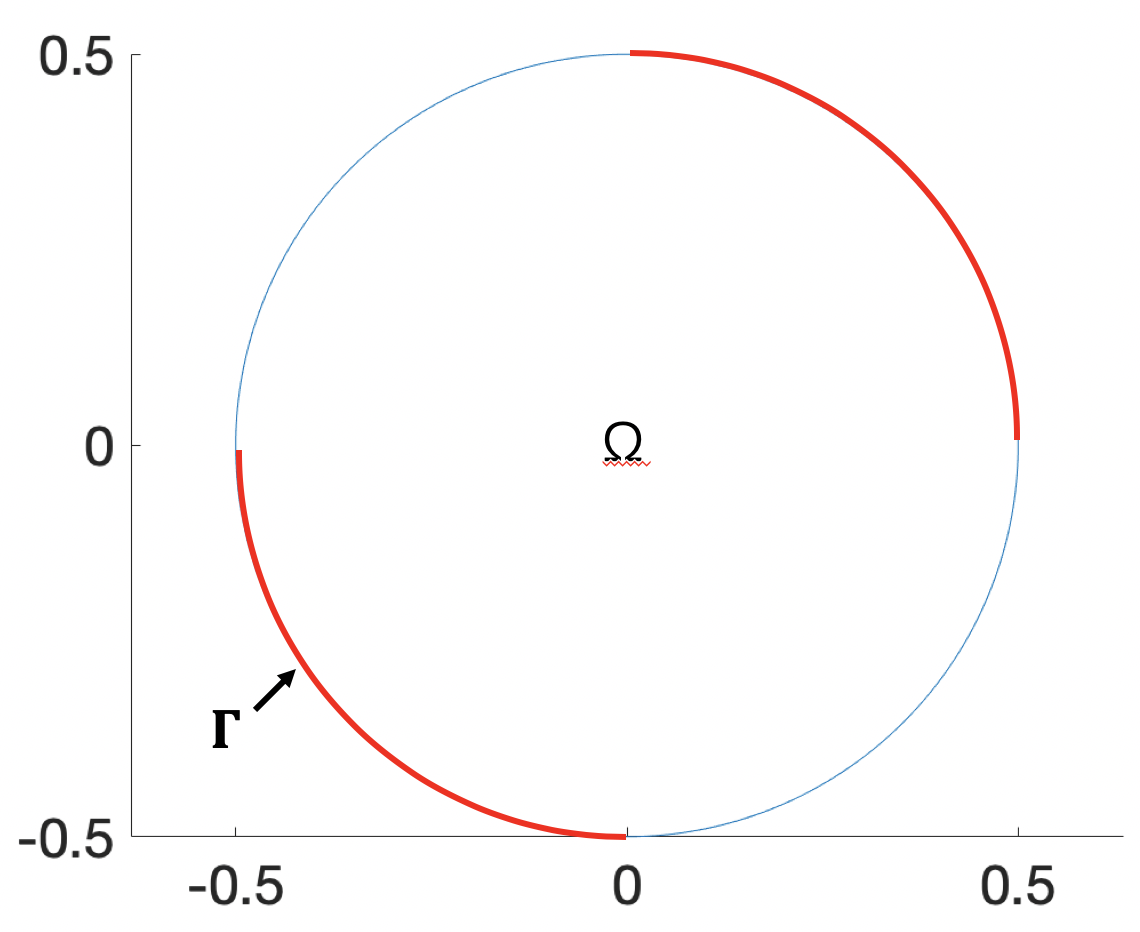}\qquad
\includegraphics[height=5.5cm]{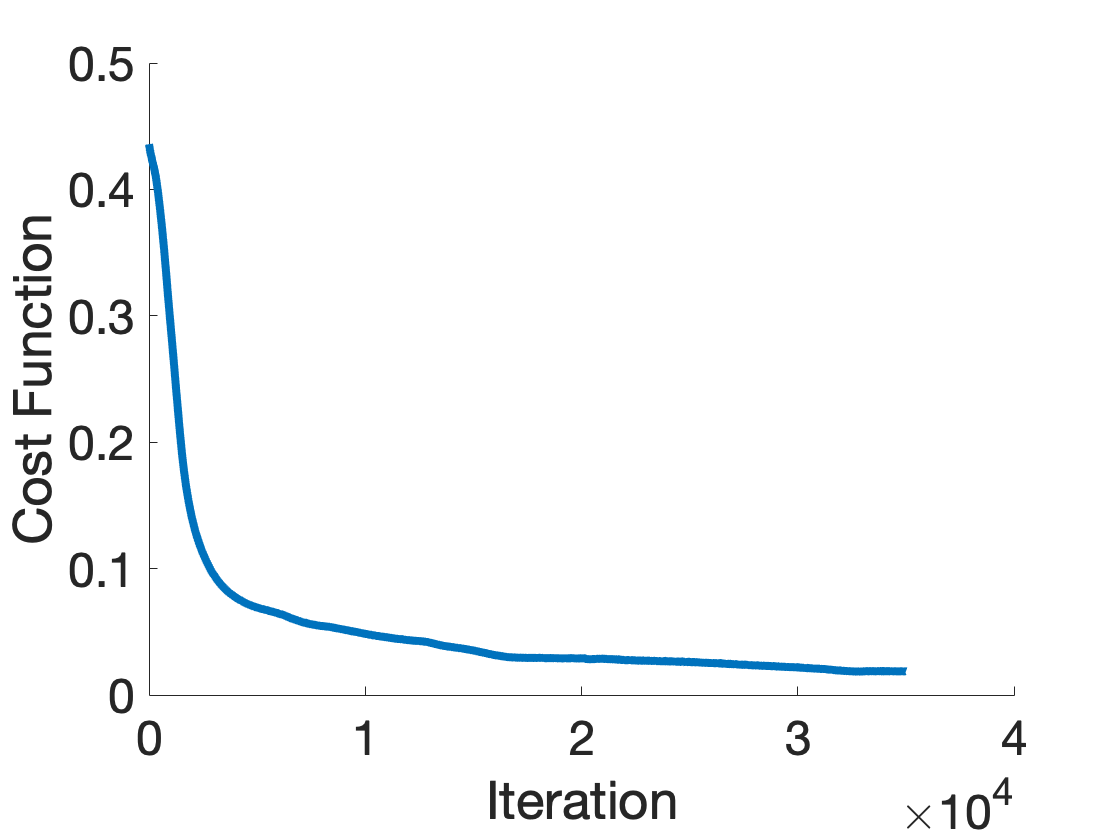}\qquad
\caption{(Left)The design area and boundary conditions of example \ref{em:parabolic}; (Right)The convergence history of cost function during iteration where steps = $10^{-4}$.}
\label{fig:condition_para}
\end {figure}

Let operator $\mathcal{L}$ in problem \eqref{eq:elliptic} be $\Delta$(Laplace operator). The structure of ANN is chosen with layers $[d,120,20,14,12,10,1]$, where $d$ is the dimension of input data. By setting initial $\mathbf{w}$ and $\mathbf{b}$ randomly in $(-1,1)$ and $\beta_1 = 0.9,\beta_2 = 0.999$ for parameters in ADAM algorithm, example of 2D case time-dependent problem is illustrated to
\begin{example}[parabolic case]
\label{em:parabolic}
The equation of time dependent problem is given as
\begin{equation}
\begin{cases}
\frac{\partial u(\mathbf{x};t)}{\partial t}+\Delta u(\mathbf{x};t) = 0 & \mathbf{x},t \ \ $in$  \ \ \Omega, \mathcal{T} \\
u(\mathbf{x};t) = e^{x_1}sin(x_2)cos(t) & \mathbf{x},t \ \ $on$ \ \ \Gamma, \mathcal{T}\\
\frac{\partial u(\mathbf{x})}{\partial \mathbf{n}} = [e^{x_1}sin(x_2), e^{x_1}cos(x_2)]*\mathbf{n}& \mathbf{x},t \ \ $on$  \ \ \Gamma,\mathcal{T}\\
u(\mathbf{x};0) = e^{x_1}sin(x_2) & \mathbf{x} \ \ $in$ \ \ \Omega
\end{cases}
\end{equation}
 domain $\Omega$ and boundary $\Gamma$ is shown in the Fig. \ref{fig:condition_para}(Left) and $ \mathcal{T} = [0,\frac{\pi}{2}]$
\end{example}

There are 10000 and 2500 points randomly sampling in $\Omega \times \mathcal{T}$ and  $\Gamma \times \mathcal{T}$, respectively. The level of noise $\delta$ is set to be $\% 1$. Fig.~\ref{fig:condition_para}(Right) shows the convergence history of the cost function, illustrating the accuracy of ANN method for time-independent case Cauchy inverse problem. More details of errors between exact and computed solution are presented in Tab.~\ref{fig:approach_para}. Four different time $t$ is chosen to  show that solution of ANN approximation is similar to the exact one. 

\begin{table}[h!]
\begin{center}
\begin{tabular}{|c|c|c|c|c|}
\hline
\ 
&$t = \frac{\pi}{5}  $    
&$t = \frac{3\pi}{10}  $ 
&$t = \frac{2\pi}{5}   $
&$t = \frac{\pi}{2}$ \\
\hline
exact solution
&\raisebox{-.5\height}{\includegraphics[width=1.3in]{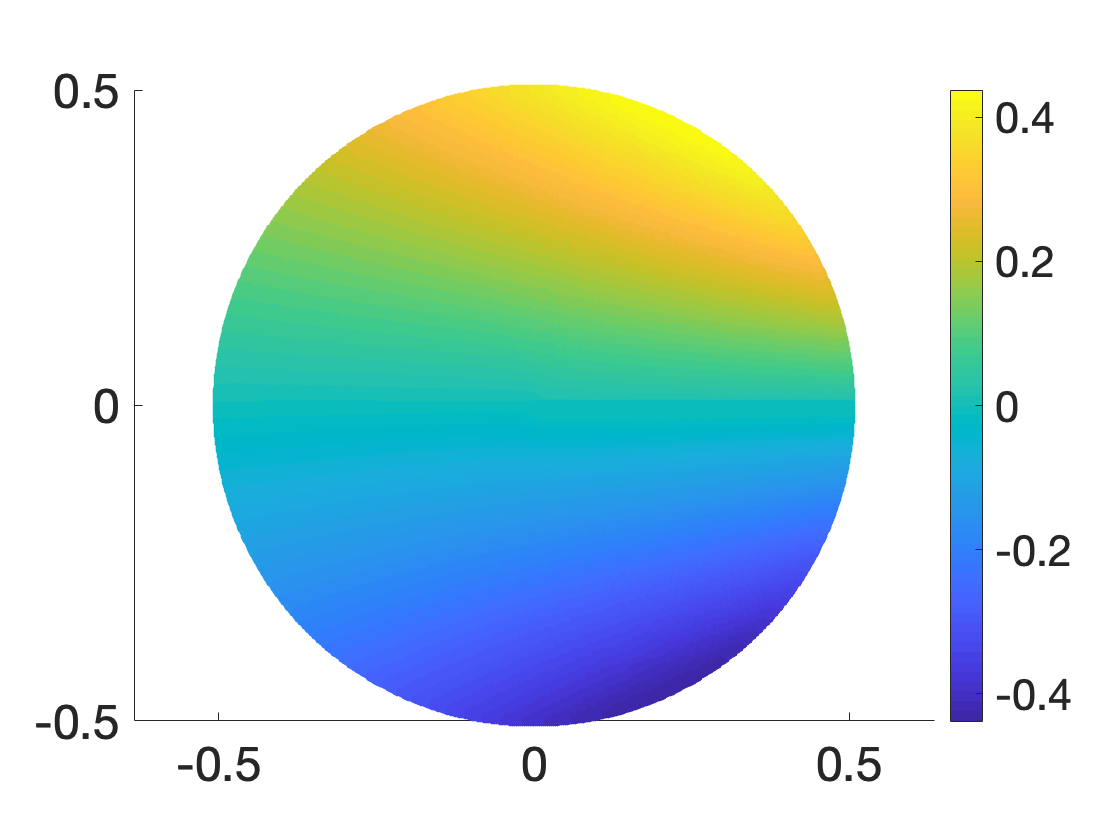}}                
&\raisebox{-.5\height}{\includegraphics[width=1.3in]{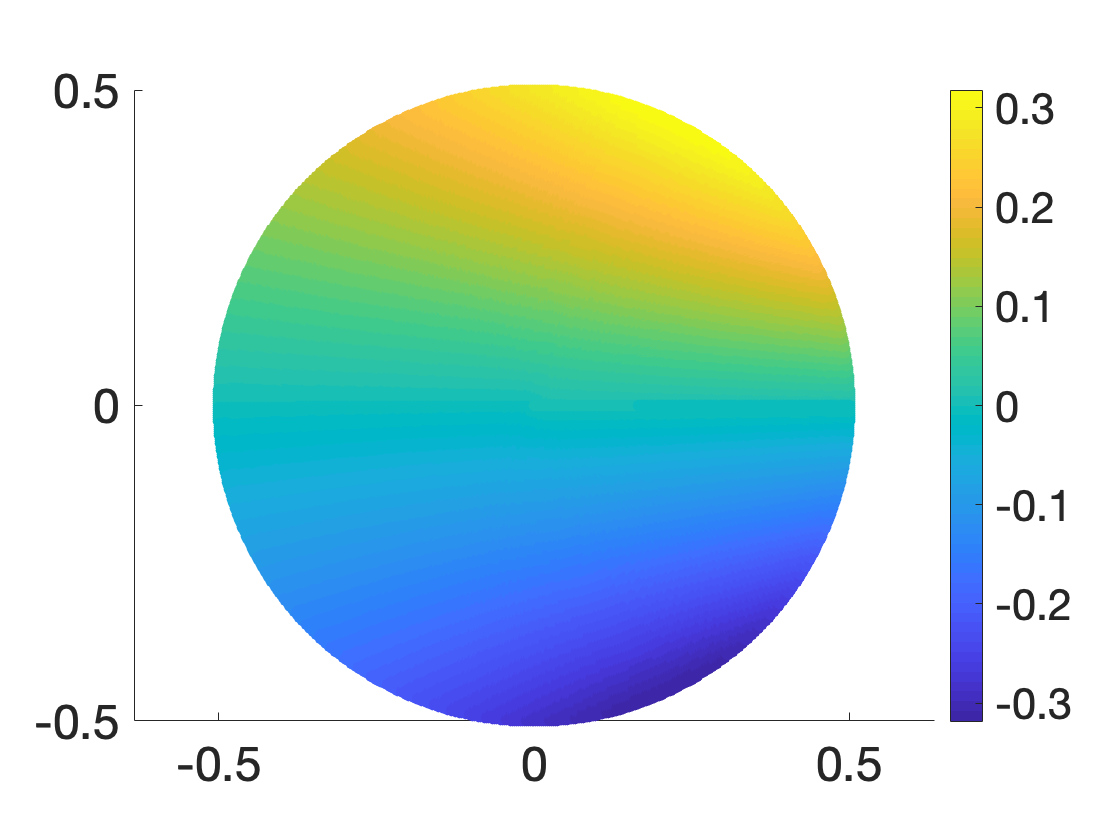}}  
&\raisebox{-.5\height}{\includegraphics[width=1.3in]{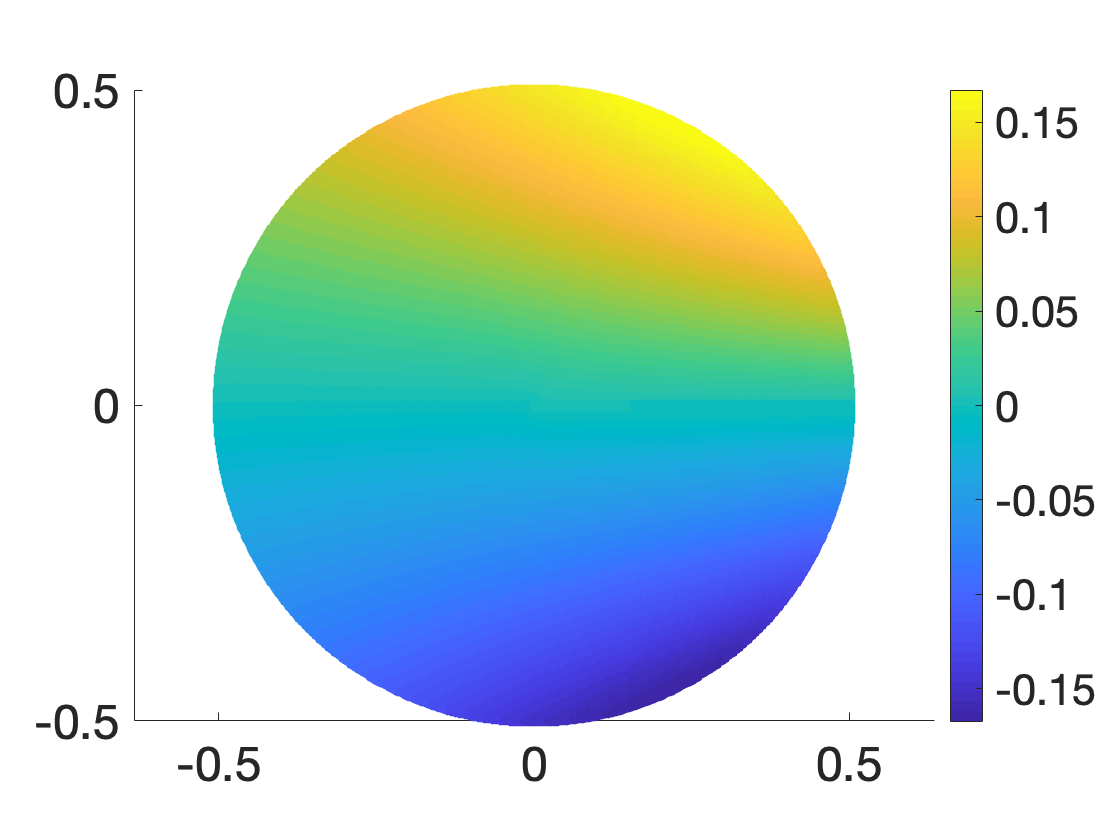}} 
&\raisebox{-.5\height}{\includegraphics[width=1.3in]{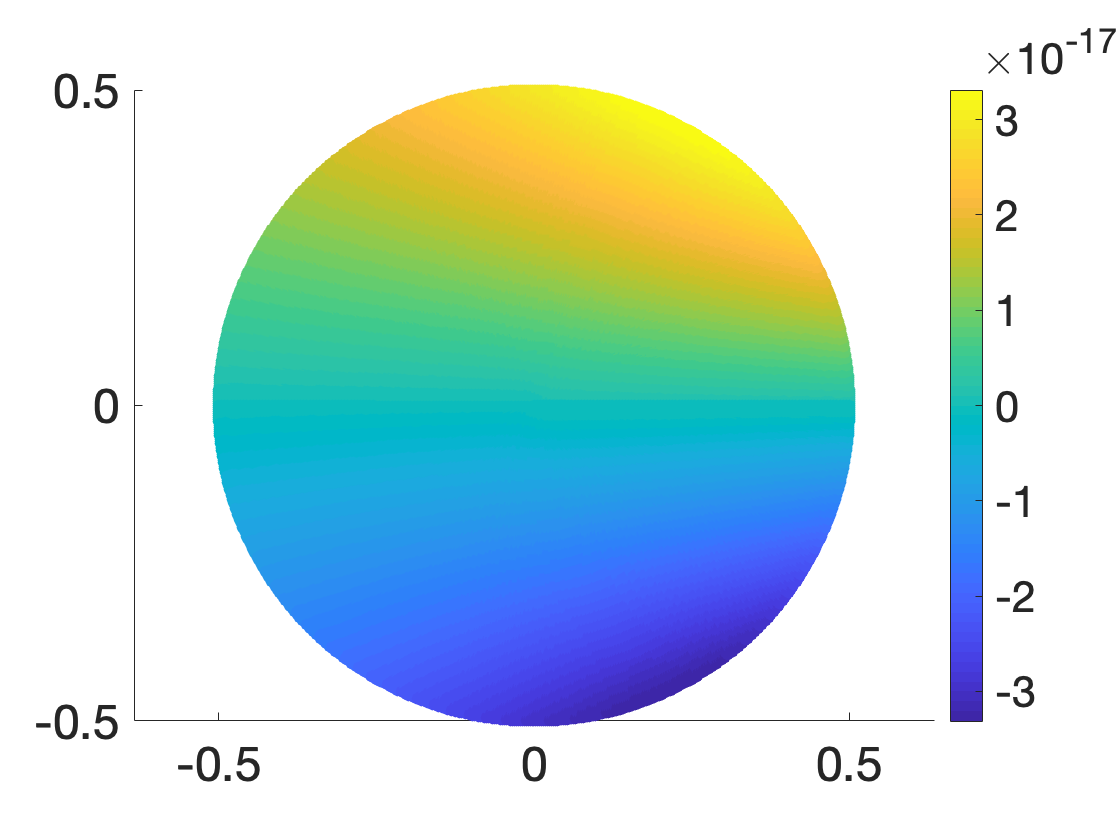}} \\
\hline
errors
&\raisebox{-.5\height}{\includegraphics[width=1.3in]{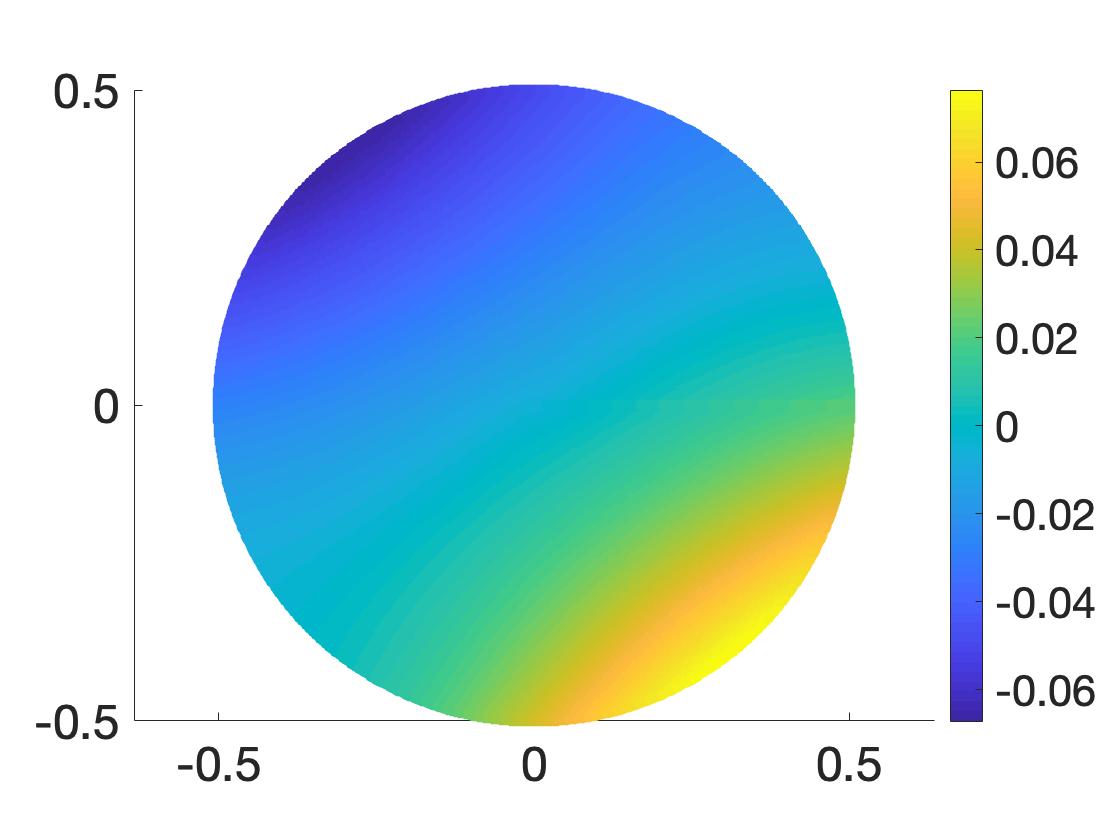}}                
&\raisebox{-.5\height}{\includegraphics[width=1.3in]{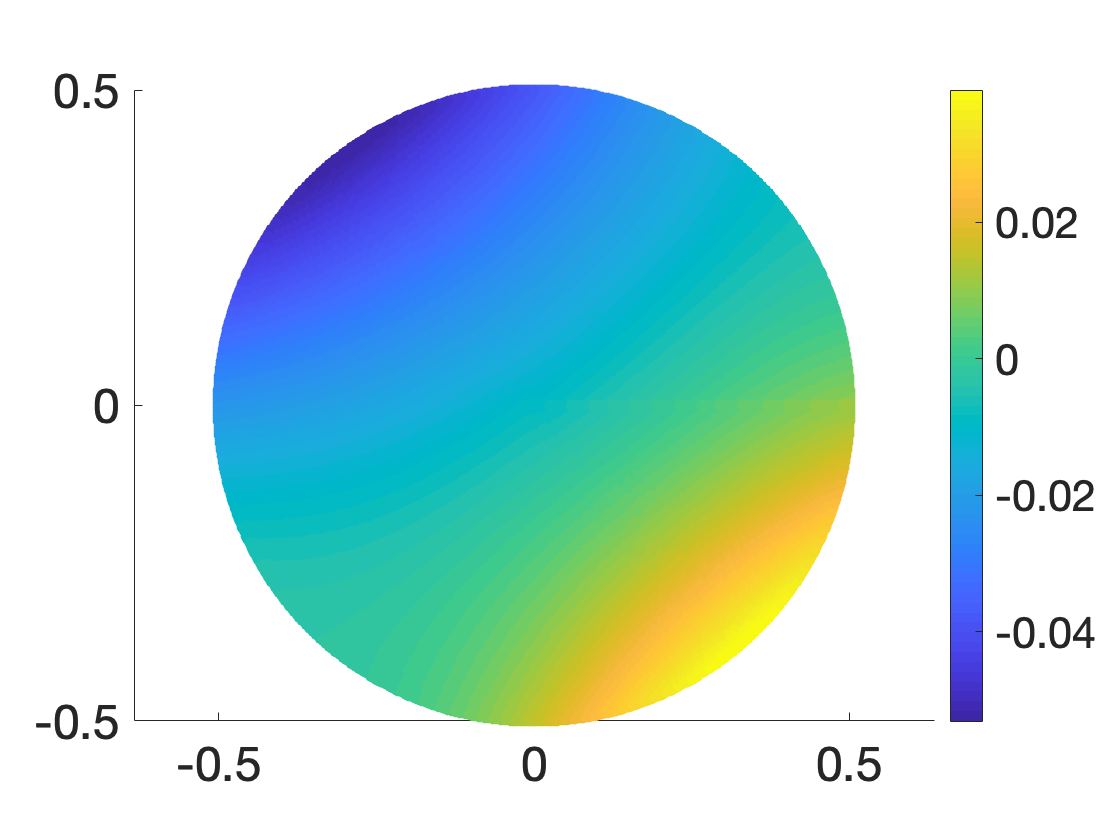}}  
&\raisebox{-.5\height}{\includegraphics[width=1.3in]{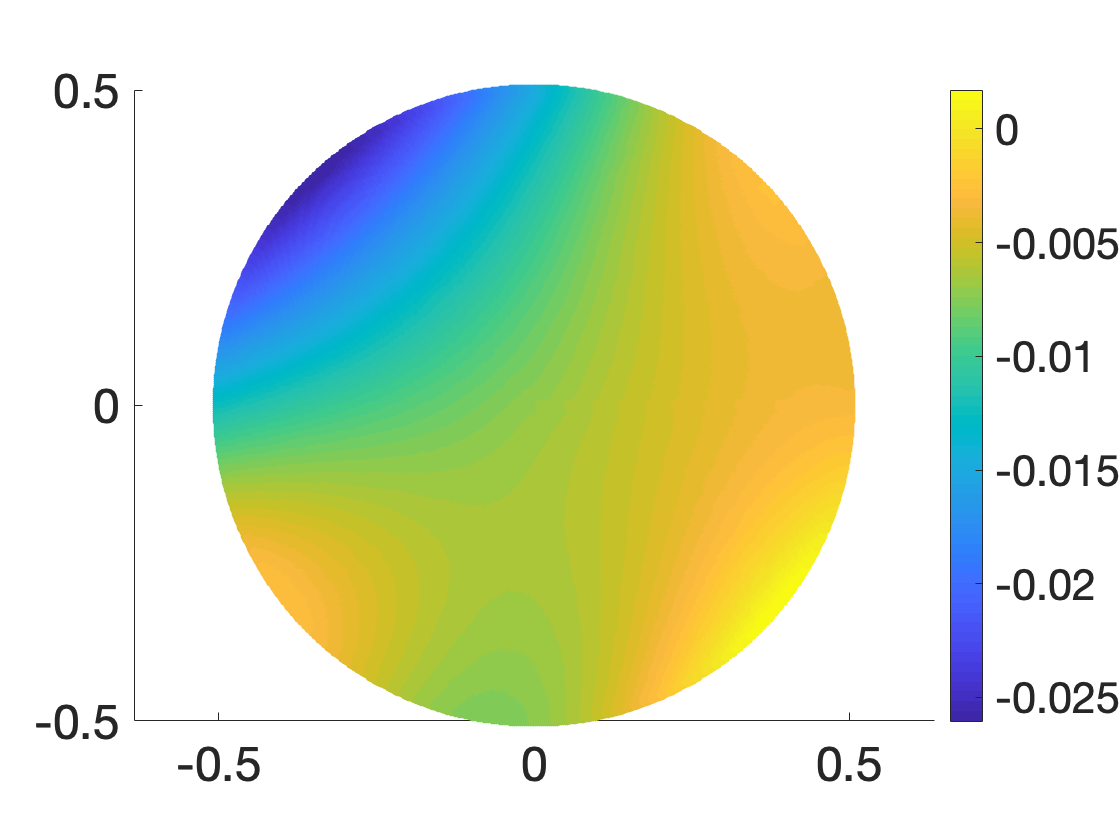}}  
&\raisebox{-.5\height}{\includegraphics[width=1.3in]{para_t5_err.png}} \\
\hline

\end{tabular}
\end{center}
\caption{The exact solution in one section and errors  between exact and computed solution}
\label{fig:approach_para}
\end{table}

%\begin{figure}
%\centering
%\subfigure[t = $\frac \pi 5$]{    \label{fig:approach_para3:a}
%		\includegraphics[width=2.4in]{para_t3_appr.png}
%	}
%	\subfigure[t = $\frac \pi 5$]{    \label{fig:approach_para3:b}
%		\includegraphics[width=2.4in]{para_t3_err.png}
%	}
%	\subfigure[t =$ \frac {3\pi}{10}$]{  \label{fig:approach_para4:a}
%		\includegraphics[width=2.4in]{para_t4_appr.png}
%	}
%		\subfigure[t = $\frac {3\pi}{10}$]{  \label{fig:approach_para4:b}
%		\includegraphics[width=2.4in]{para_t4_err.png}
%	}
%			\subfigure[t = $\frac {2\pi}{5}$]{  \label{fig:approach_para5:a}
%		\includegraphics[width=2.4in]{para_t5_appr.png}
%	}
%			\subfigure[t = $\frac {2\pi}{5}$]{  \label{fig:approach_para5:b}
%		\includegraphics[width=2.4in]{para_t5_err.png}
%	}
%			\subfigure[t = $\frac {\pi}{2}$]{  \label{fig:approach_para6:a}
%		\includegraphics[width=2.4in]{para_t6_appr.png}
%	}
%			\subfigure[t = $\frac {\pi}{2}$]{  \label{fig:approach_para6:b}
%		\includegraphics[width=2.4in]{para_t6_err.png}
%	}
%	\caption{(Left)The exact solution in one section;(Right)errors  between exact and computed solution}
%	\label{fig:approach_para} %% label for entire figure
%\end{figure}

\subsection{Numerical validations of time-independent case} 
\begin{figure}[H]
\centering
\includegraphics[height=5.5cm]{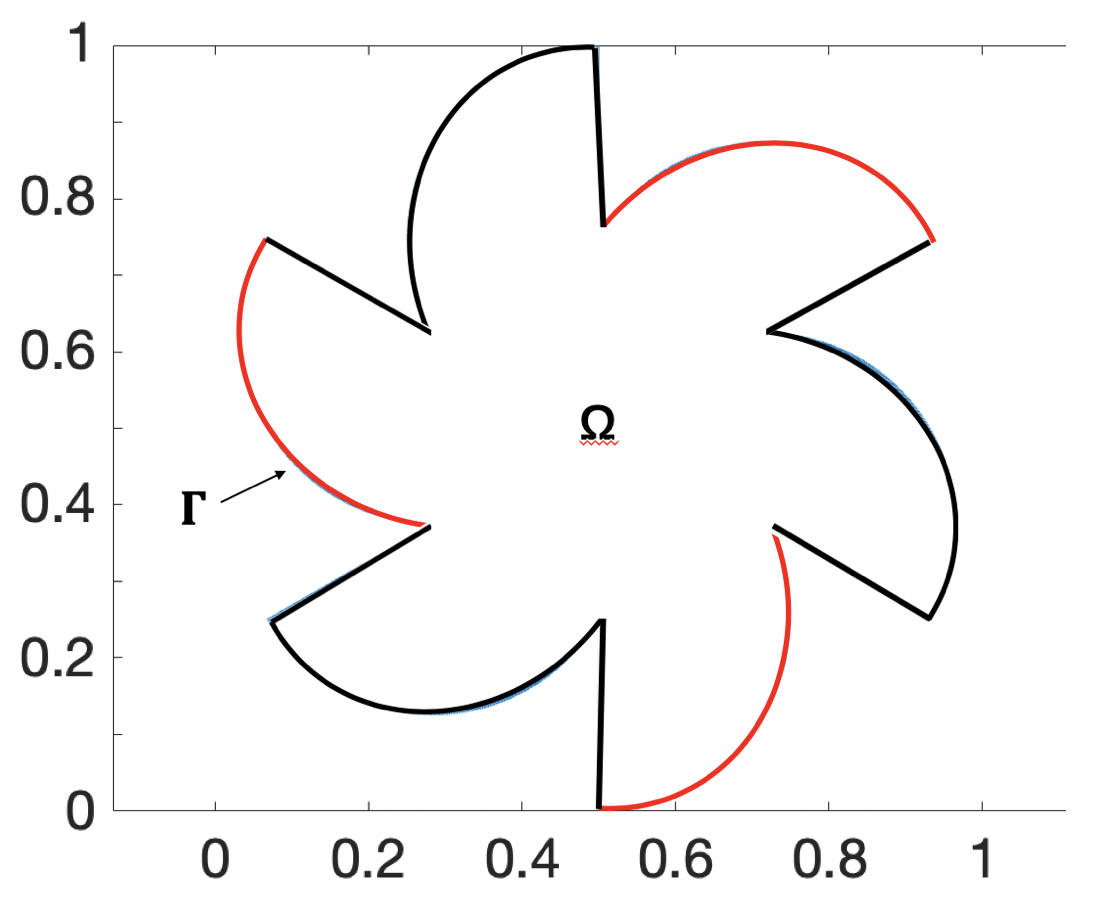}\qquad
\includegraphics[height=5.5cm]{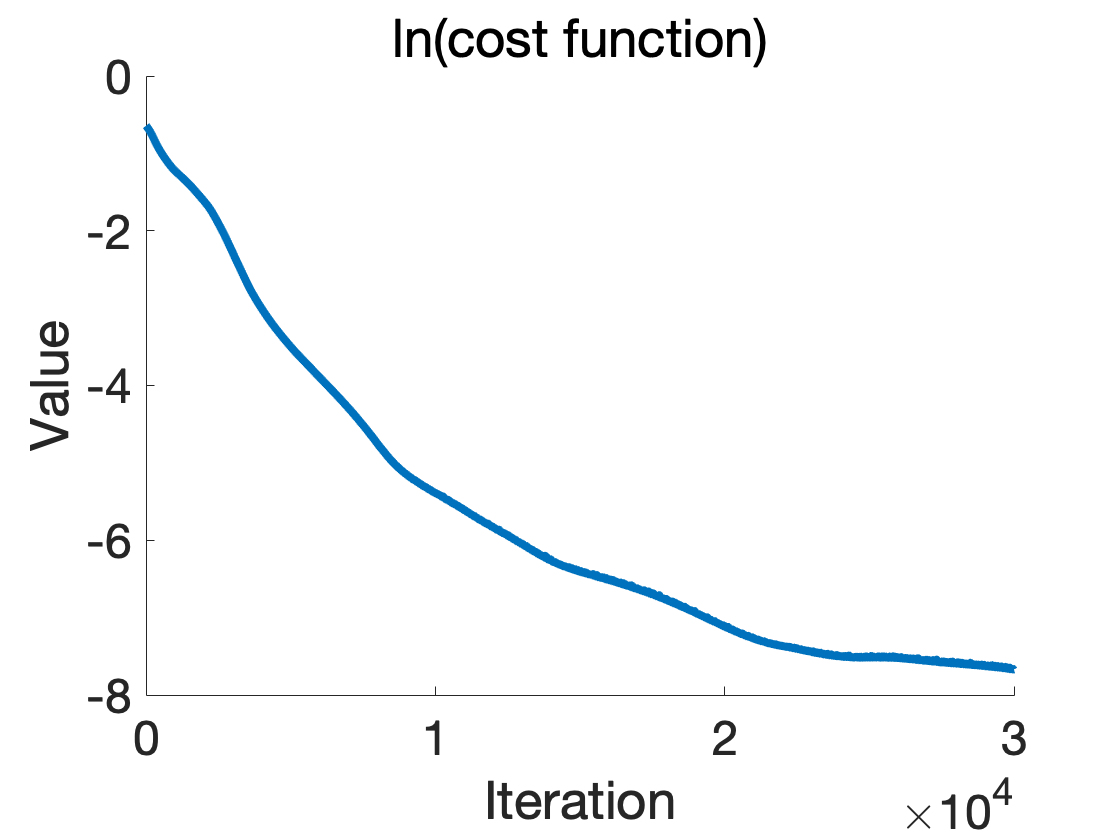}\qquad
\caption{(Left)The design area and boundary conditions of 2d case; (Right)The convergence history of cost function during iterations where time steps = $10^{-4}$.}
\label{fig:design_initial}
\end {figure}

Let operator $\mathcal{L}$ in problem \eqref{eq:elliptic} be $\Delta$(Laplace operator). The structure of ANN is chosen with layers $[d,120,20,14,12,10,1]$, where $d$ is the dimension of input data. By setting initial $\mathbf{w}$ and $\mathbf{b}$ randomly in $(-1,1)$ and $\beta_1 = 0.9,\beta_2 = 0.999$ for parameters in ADAM algorithm,  time-independent  problem in 2D case is illustrated in the following example at first.
\begin{example}
\label{em:2d_ori}
The equation of time-independent problem is given as
\begin{equation}
\begin{cases}
\Delta u(\mathbf{x}) = 0 & \mathbf{x} \ \ $in$  \ \ \Omega\\
u(\mathbf{x}) = e^{x_1}sin(x_2) & \mathbf{x} \ \ $on$ \ \ \Gamma\\
\frac{\partial u(\mathbf{x})}{\partial \mathbf{n}} = [e^{x_1}sin(x_2), e^{x_1}cos(x_2)]*\mathbf{n}& \mathbf{x} \ \ $on$  \ \ \Gamma
\end{cases}
\end{equation}
 domain $\Omega$ and boundary $\Gamma$ is shown in the Fig. \ref{fig:design_initial}(Left). 
\end{example}

\begin{figure}[H]
\centering
\includegraphics[height=5.7cm]{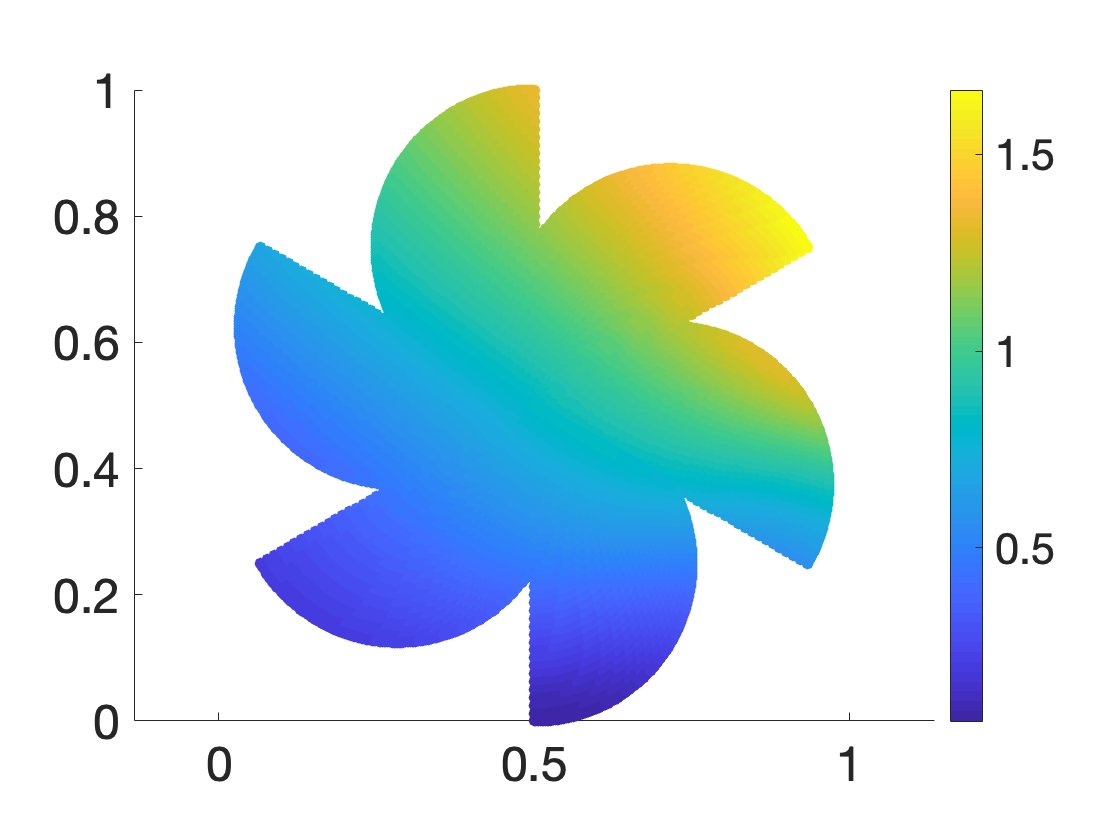}\qquad
\includegraphics[height=5.7cm]{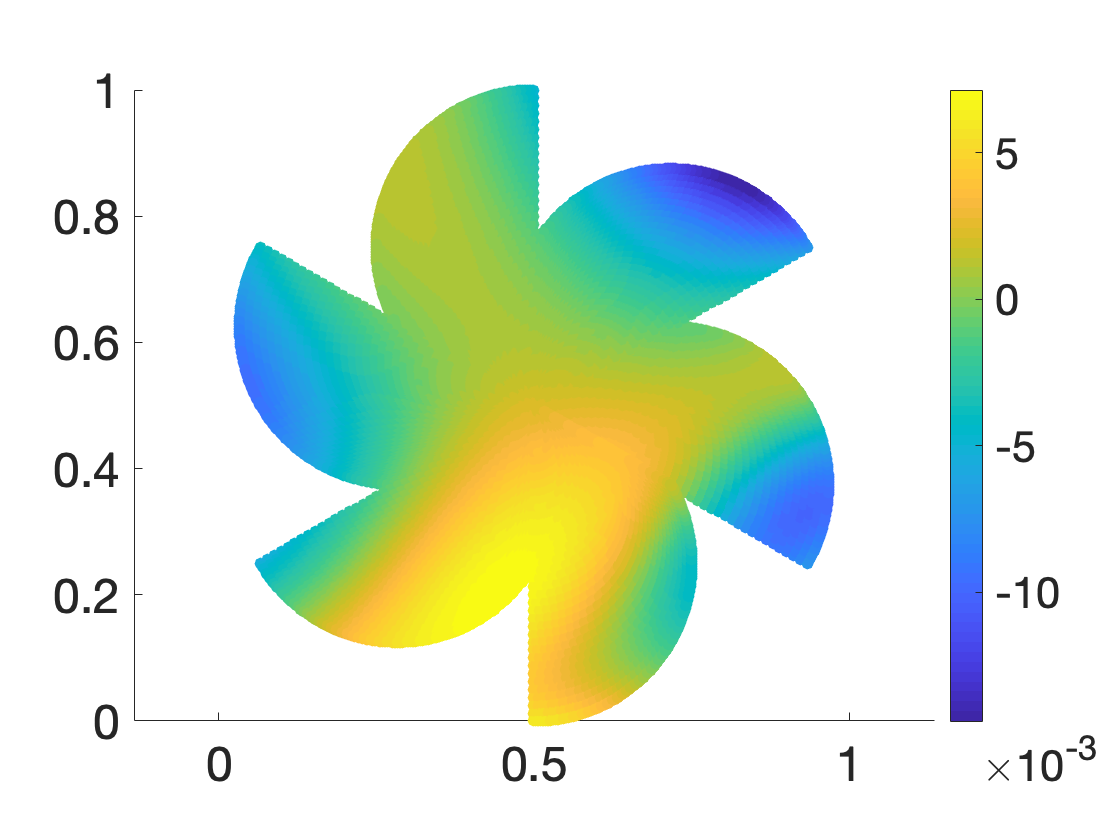}\qquad
\caption{(Left) The ANN solution to example \ref{em:2d_ori}; (Right) the errors  between exact and computed solution.}
\label{fig:approach}
\end {figure}

There are 2000 points and 500 points randomly sampling in $\Omega$ and $\Gamma$, respectively. The level of noise $\delta$ is set to be zero at first and the discussion of influence of noise is presented in following section \ref{subsec:prop}. Fig.~\ref{fig:design_initial} (right) shows the convergence history of the cost function and Fig.~\ref{fig:approach} presents  ANN solution and errors between exact solution and neural network approach. To display the depression curve more clearly, we add $\ln$ function to the cost function data. It can be observed that the optimization process is stable up to 30000 iteration steps, which verifies the accuracy of ANN for time-independent case in 2D.

\begin{figure}[H]
\centering
\includegraphics[height=5.5cm]{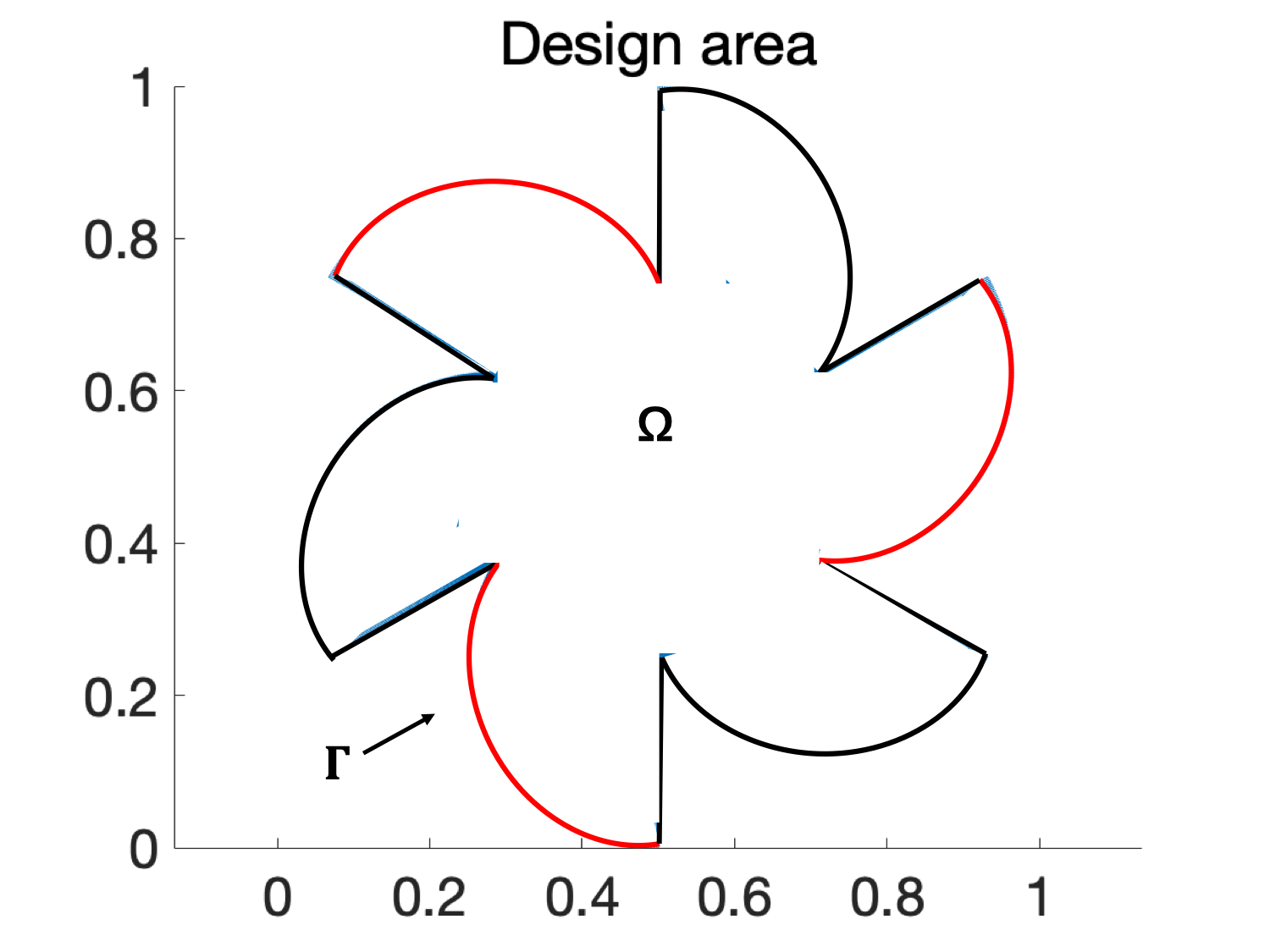}\qquad
\includegraphics[height=5.5cm]{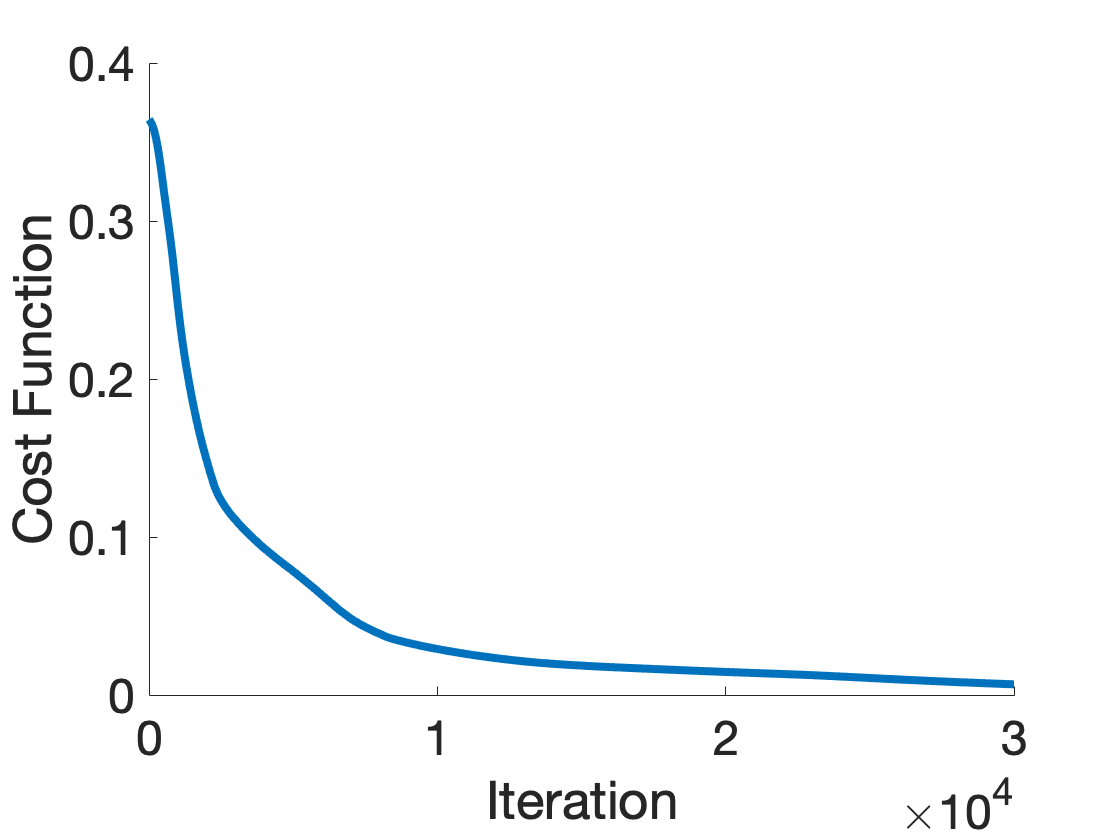}\qquad
\caption{(Left)The design area and boundary conditions of 2d case; (Right)The convergence history of cost function during iteration where steps = $10^{-4}$.}
\label{fig:design_sing}
\end {figure}

To illustrate the universality of this method, more examples for time-independent case in 2D are given. With same parameters as example \ref{em:2d_ori}, we present the following example:
\begin{example}
\label{em:2d_sing}
The equation of problem is given as
\begin{equation}
\begin{cases}
\Delta u(\mathbf{x}) = 0 & \mathbf{x} \ \ $in$  \ \ \Omega\\
u(\mathbf{x}) = \ln\left(\sqrt{(x_1-1)^2 + (x_2-1)^2}\right)& \mathbf{x} \ \ $on$ \ \ \Gamma\\
\frac{\partial u(\mathbf{x})}{\partial \mathbf{n}} = [\frac{2(x_1-1)}{(x_1-1)^2 + (x_2-1)^2}, \frac{2(x_2-1)}{(x_1-1)^2 + (x_2-1)^2}]*\mathbf{n}& \mathbf{x} \ \ $on$  \ \ \Gamma
\end{cases}
\end{equation}
where domain $\Omega$ and boundary $\Gamma$ is given in Fig. \ref{fig:design_sing}(Left). 
\end{example}

\begin{figure}[H]
\centering
\includegraphics[height=5.6cm]{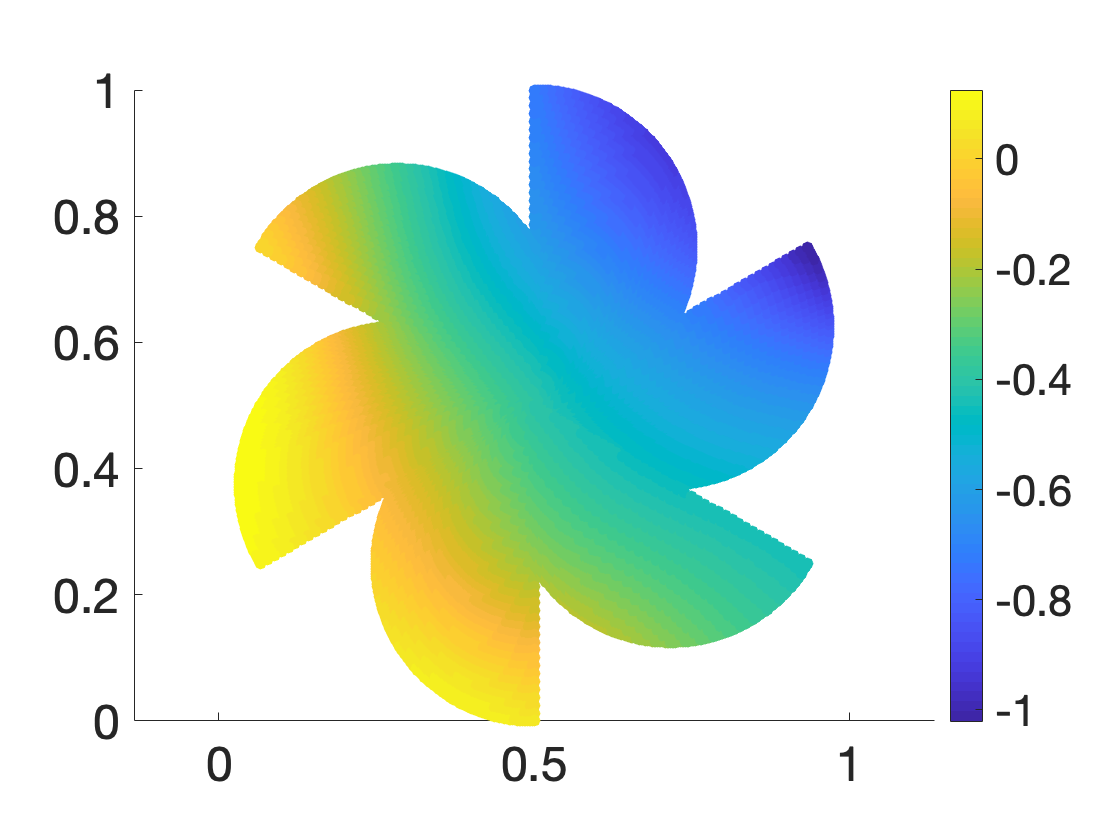}\qquad
\includegraphics[height=5.6cm]{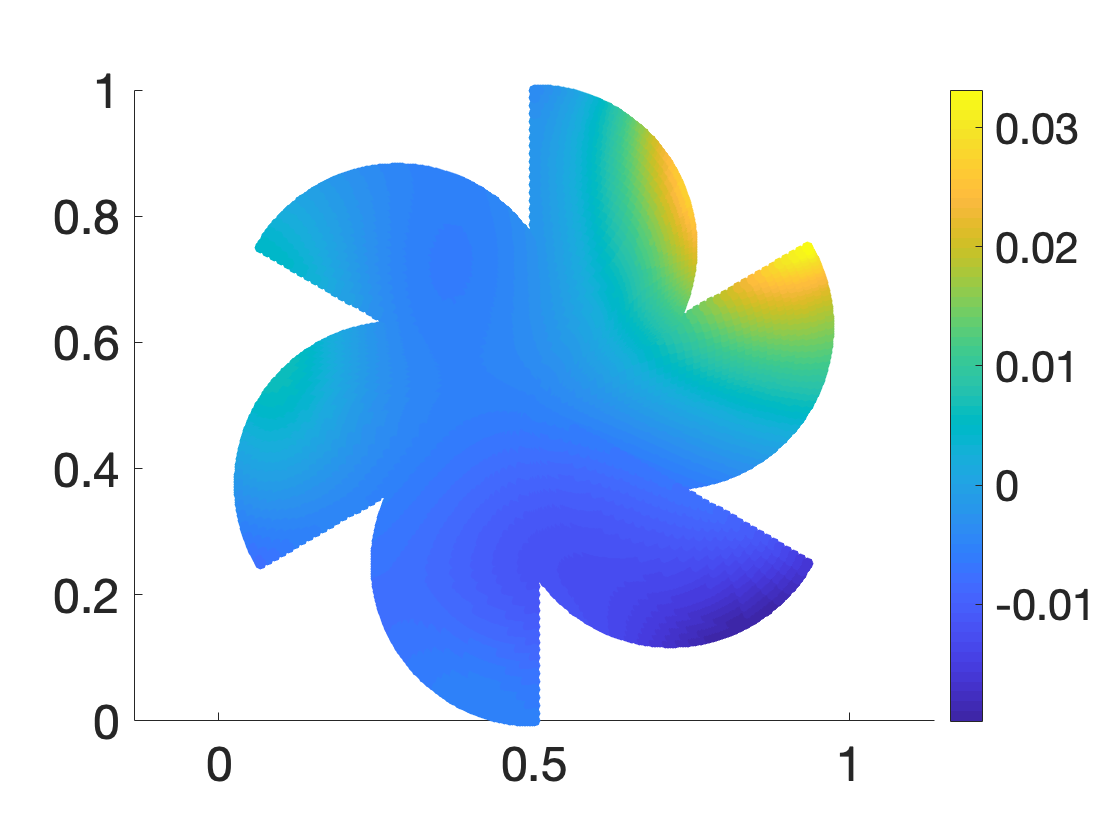}
\caption{(Left) The ANN solution to example \ref{em:2d_sing}; (Right) the errors  between exact and computed solution.}
\label{fig:approach_sing}
\end {figure}

The level of noise $\delta$ is set to be $\% 1$. Fig.~\ref{fig:design_sing} (Right) shows the convergence history of the cost function, and errors between exact and computed solution are presented in Fig.~\ref{fig:approach_sing}. As the above results show,  ANN method for Cauchy inverse problem of time-independent case is worked on 2D. Higher spatial dimensions cases are considered in the following examples.

\begin{example}[3D case]
\label{em:3d}
The equation of 3D case is given as
\begin{equation}
\begin{cases}
\Delta u(\mathbf{x}) = 0 & \mathbf{x} \ \ $in$  \ \ \Omega\\
u(\mathbf{x}) = \sinh(\sqrt{2} x_1)\sin( x_2)\sin( x_3) & \mathbf{x} \ \ $on$ \ \ \Gamma\\
\frac{\partial u(\mathbf{x})}{\partial \mathbf{n}} = \nabla\left(\sinh(\sqrt{2} x_1)\sin( x_2)\sin( x_3)\right)*\mathbf{n}& \mathbf{x} \ \ $on$  \ \ \Gamma
\end{cases}
\end{equation}

Let domain $\Omega$ be $(x_1-0.5)^2 + (x_2-0.5)^2 + (x_3-0.5)^2 < 0.25$ and boundary $\Gamma$ be $ (x_1-0.5)^2 + (x_2-0.5)^2 + (x_3-0.5)^2 = 0.25, x_3 \in [0,0.25]\cup[0.5,0.75]$.
\end{example}

\begin{figure}[H]
\centering
\subfigure[Cost Function]{ \label{fig:design_3d:a}
\includegraphics[width=2.8in]{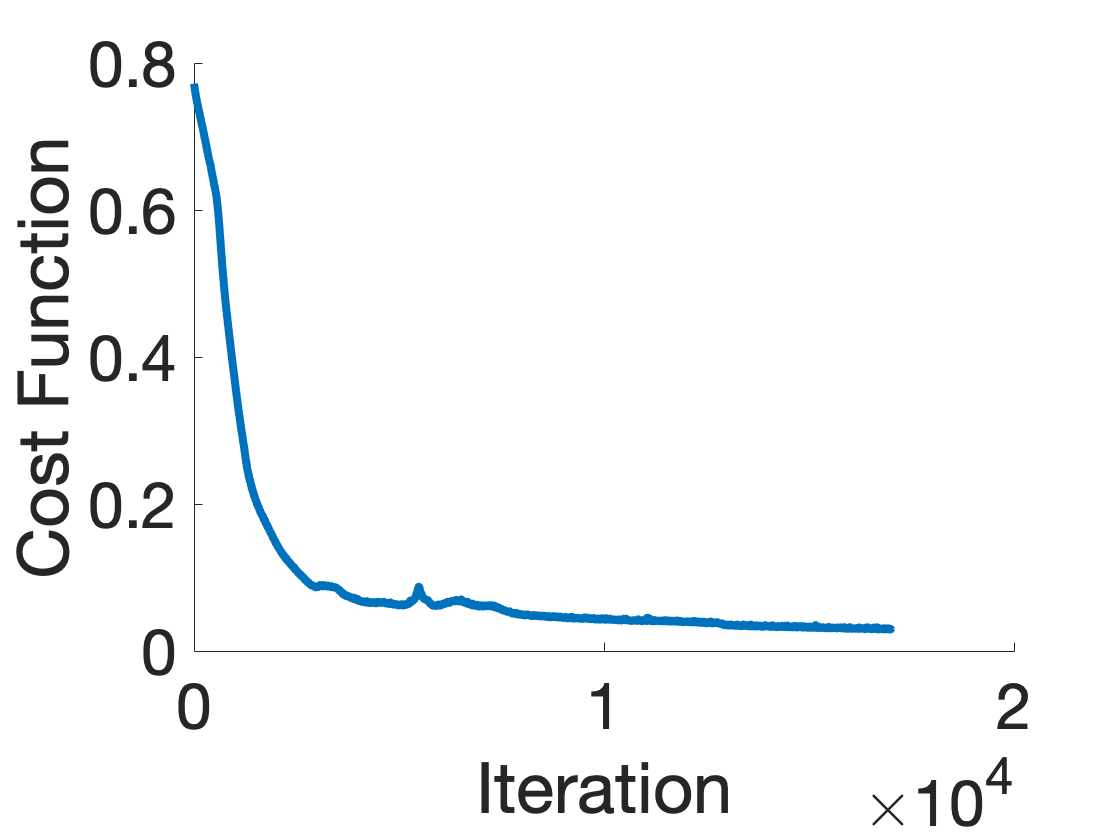}}
\subfigure[Edge Error]{ \label{fig:design_3d:b} 
\includegraphics[width=2.8in]{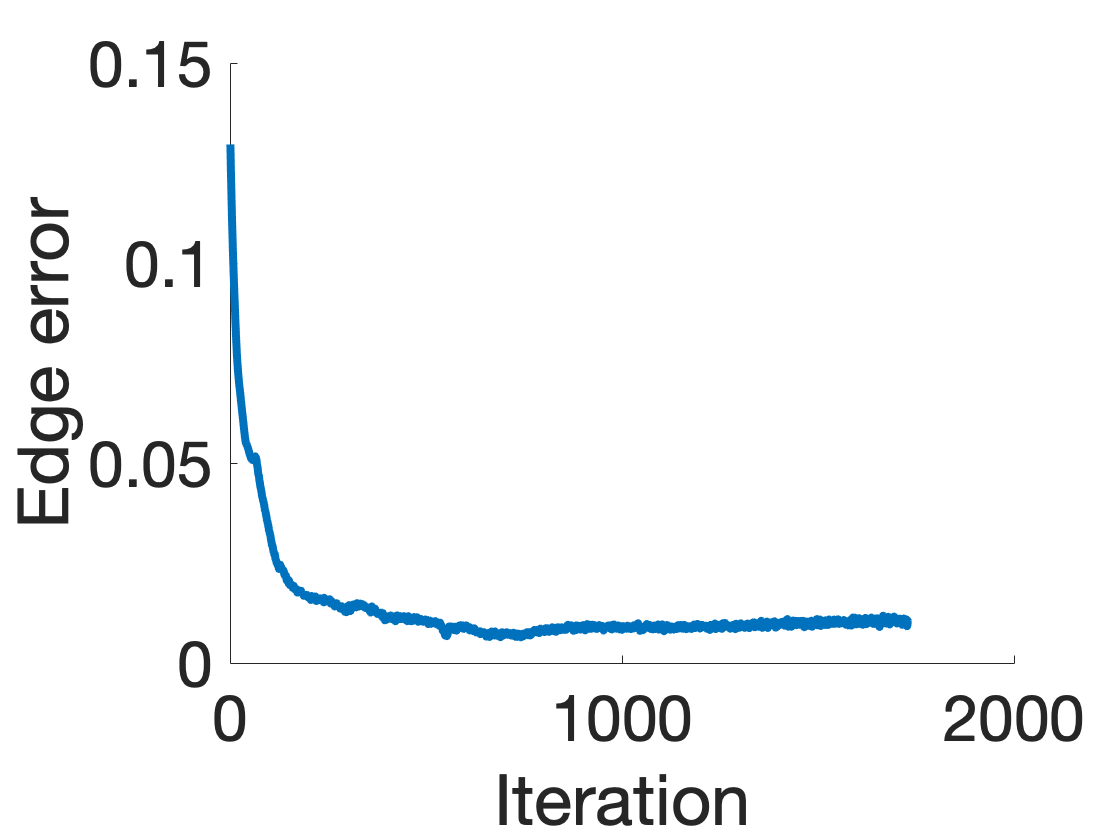}}

\caption{ (a) the convergence history of example \ref{em:3d} and (b) the errors in the edges where step = $5\times 10^{-4}$}
\label{fig:compare} 
\end{figure}

There are 10000 points and 2500 randomly sampling in $\Omega$  and $\Gamma$, respectively. The structure of ANN is chosen with layers $[3,120,20,14,12,10,1]$ and level of noise $\delta$ is set to be $\% 1$. Fig.~\ref{fig:design_3d:a} shows the convergence history of the cost function and errors between exact function and computed solution on edges are displayed in Fig.~\ref{fig:design_3d:b}. We choose five sections of $\Omega$ and the results are presented in Tab.~\ref{fig:approach_3d}.

\begin{table}[h!]
\begin{center}
\begin{tabular}{|c|c|c|c|c|}
\hline
\ 
&$x_3 = 0.8536  $    
&$x_3 = 0.75  $ 
&$x_3= 0.5 $
&$x_3 = 0.25$ \\
\hline
exact solution
&\raisebox{-.5\height}{\includegraphics[width=1.3in]{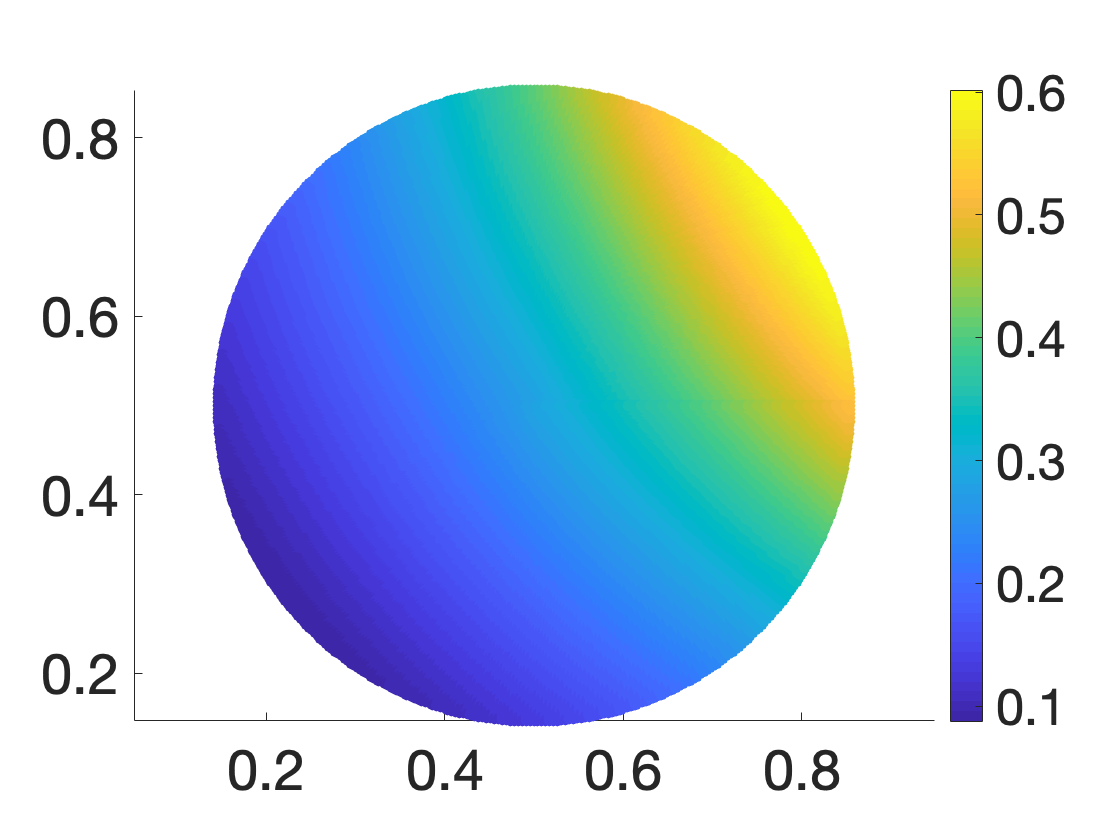}}                
&\raisebox{-.5\height}{\includegraphics[width=1.3in]{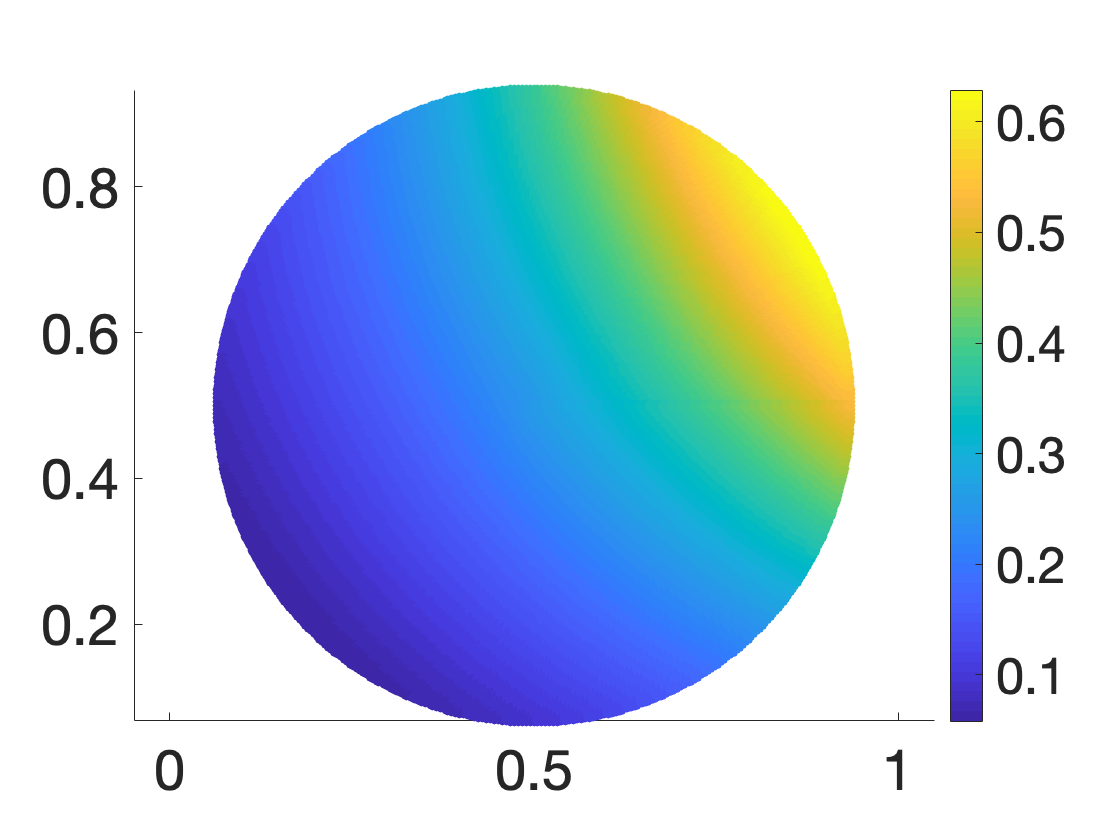}}  
&\raisebox{-.5\height}{\includegraphics[width=1.3in]{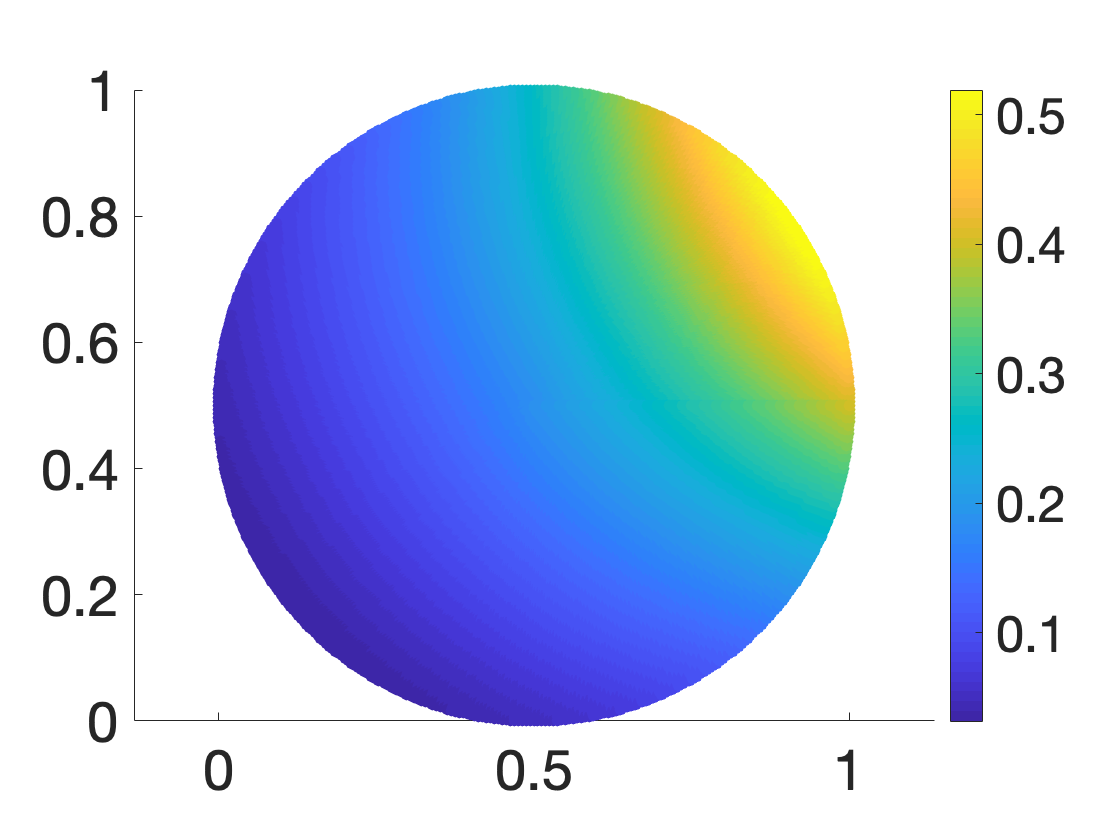}} 
&\raisebox{-.5\height}{\includegraphics[width=1.3in]{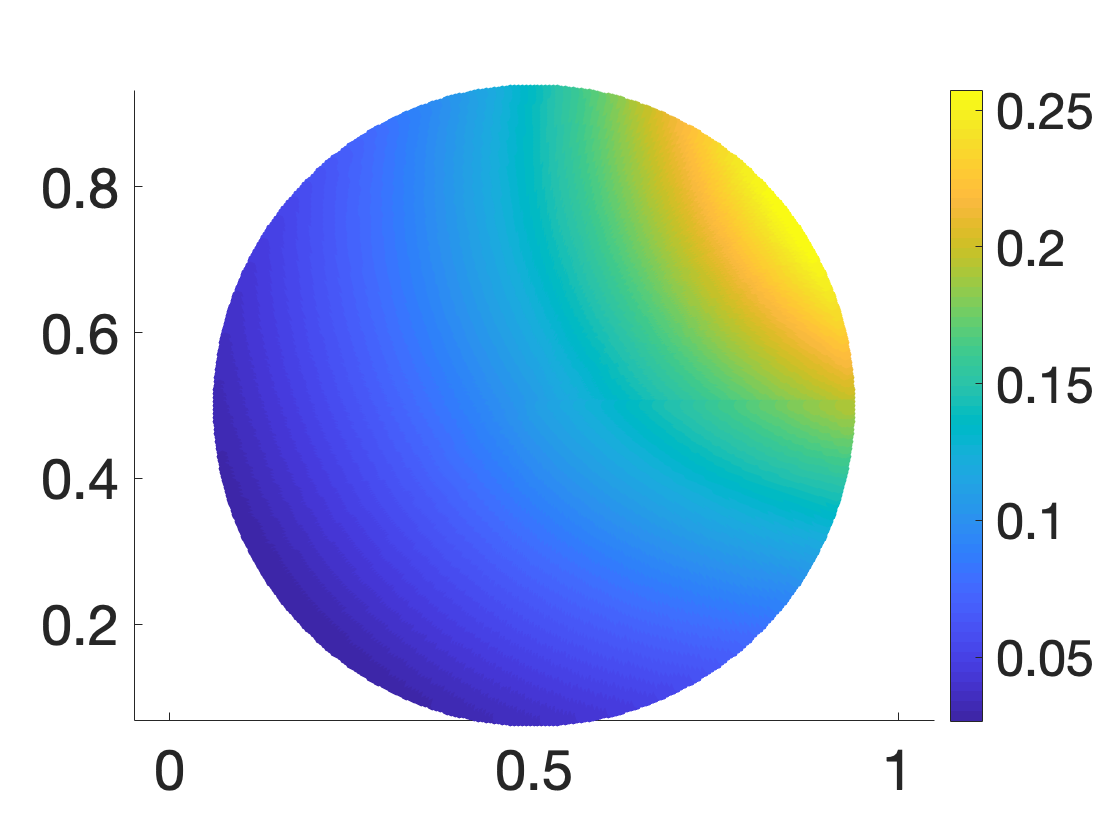}} \\
\hline
errors
&\raisebox{-.5\height}{\includegraphics[width=1.3in]{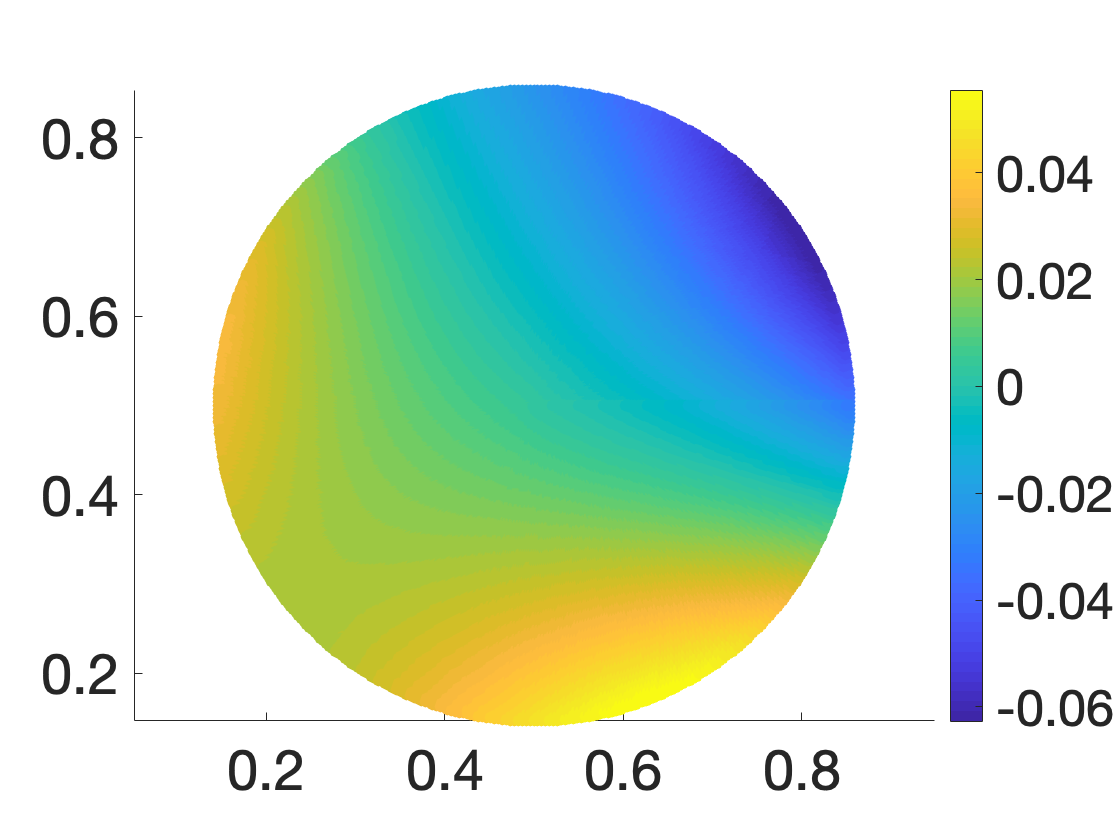}}                
&\raisebox{-.5\height}{\includegraphics[width=1.3in]{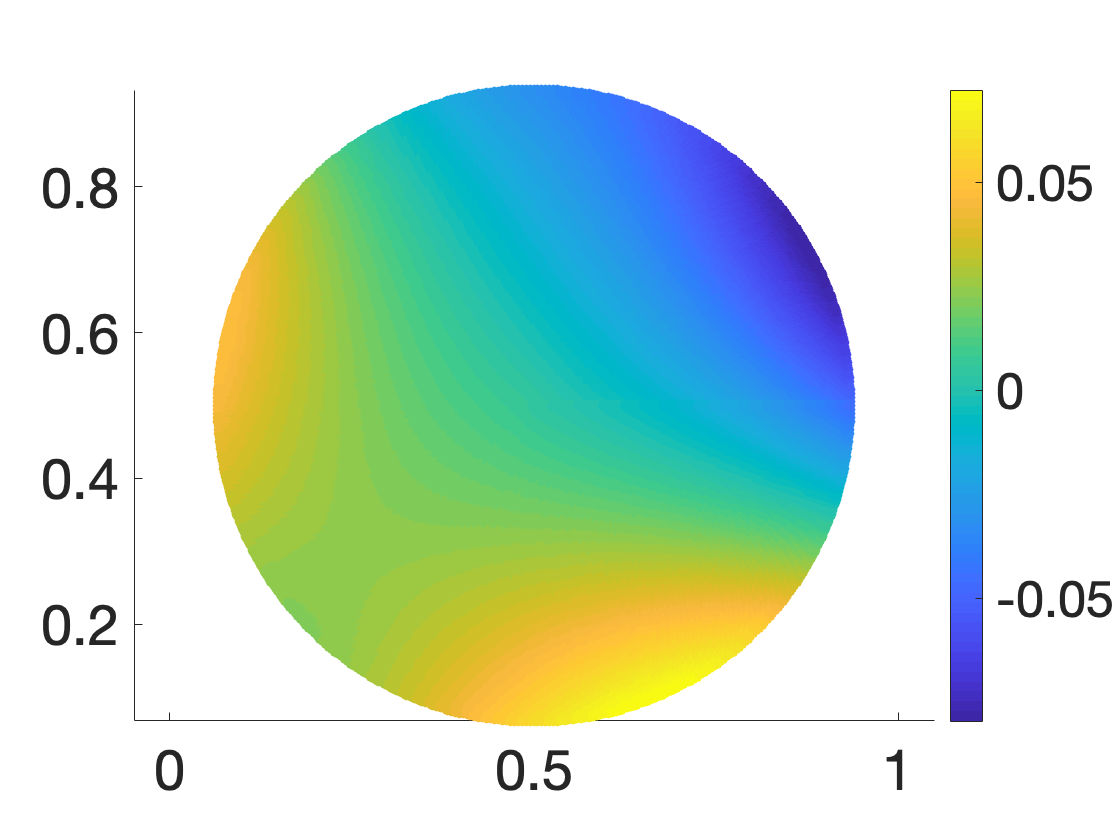}}  
&\raisebox{-.5\height}{\includegraphics[width=1.3in]{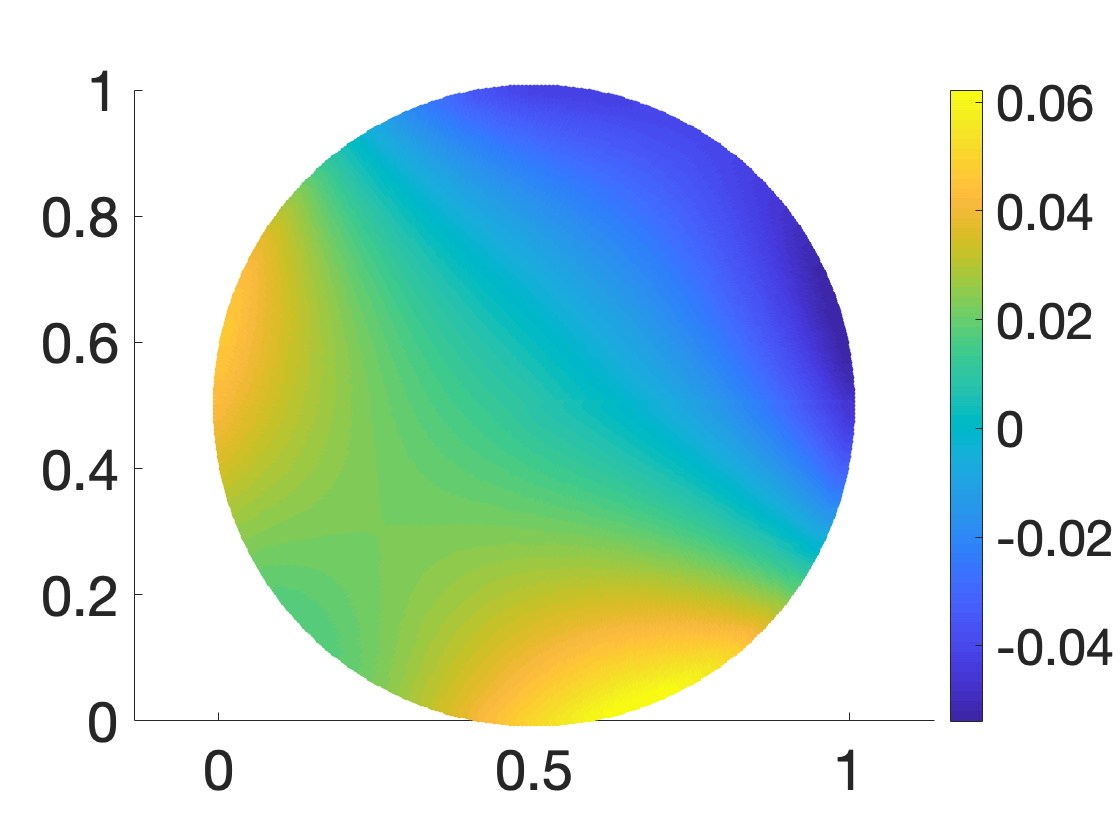}}  
&\raisebox{-.5\height}{\includegraphics[width=1.3in]{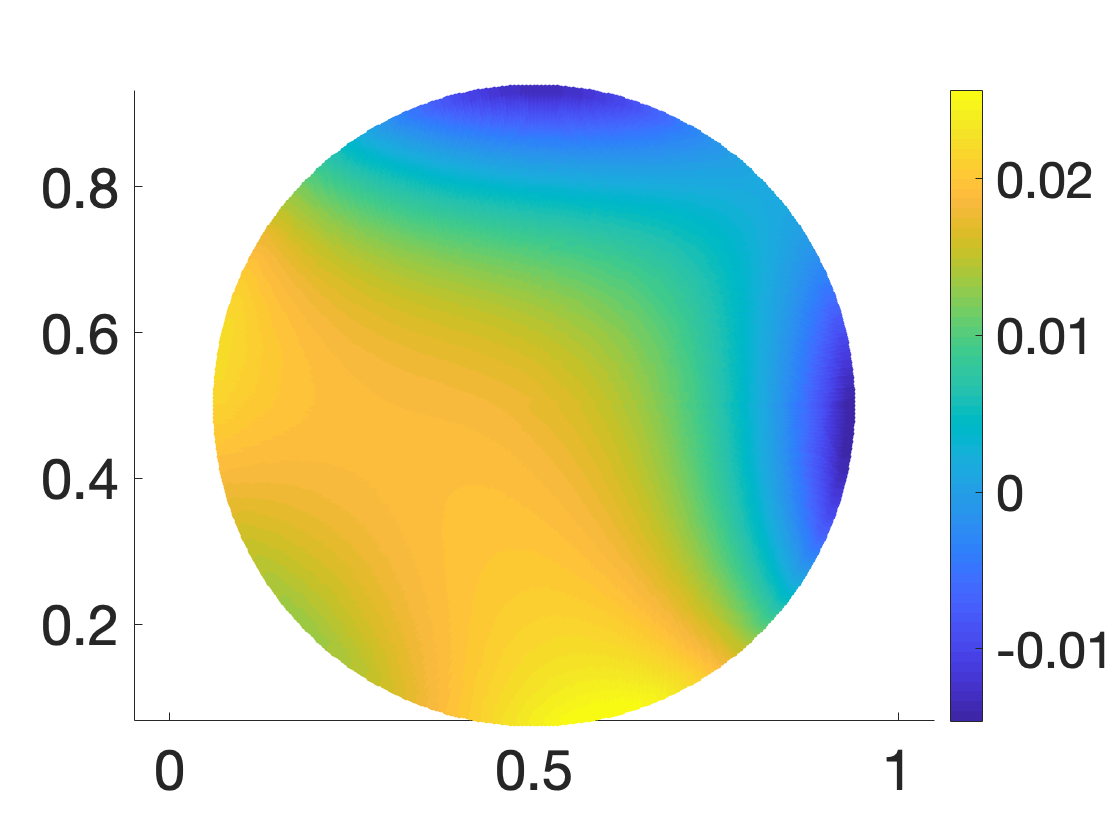}} \\
\hline

\end{tabular}
\end{center}
\caption{The exact solution in one section and errors  between exact and computed solution}
\label{fig:approach_3d}
\end{table}

%\begin{figure}
%\centering
%\subfigure[z = 0.8536]{    \label{fig:approach_3d:a}
%		\includegraphics[width=2.4in]{z1_3d_app.png}
%	}
%	\subfigure[z = 0.8536]{    \label{fig:approach_3:b}
%		\includegraphics[width=2.4in]{z1_3d_err.png}
%	}
%	\subfigure[z = 0.75]{  \label{fig:approach_3d:c}
%		\includegraphics[width=2.4in]{z2_3d_app.png}
%	}
%		\subfigure[z = 0.75]{  \label{fig:approach_3d:d}
%		\includegraphics[width=2.4in]{z2_3d_err.png}
%	}
%			\subfigure[z = 0.5]{  \label{fig:approach_3d:d}
%		\includegraphics[width=2.4in]{z3_3d_app.png}
%	}
%			\subfigure[z = 0.5]{  \label{fig:approach_3d:d}
%		\includegraphics[width=2.4in]{z3_3d_err.png}
%	}
%			\subfigure[z = 0.25]{  \label{fig:approach_3d:d}
%		\includegraphics[width=2.4in]{z4_3d_app.png}
%	}
%			\subfigure[z = 0.25]{  \label{fig:approach_3d:d}
%		\includegraphics[width=2.4in]{z4_3d_err.png}
%	}
%	\caption{(Left)The ANN solution in one section;(Right)errors  between exact and computed solution}
%	\label{fig:approach_3d} %% label for entire figure
%\end{figure}

\begin{figure}[H]
\centering
\subfigure[ln(Cost function)]{    \label{fig:ndcase:a}
		\includegraphics[width=2.7in]{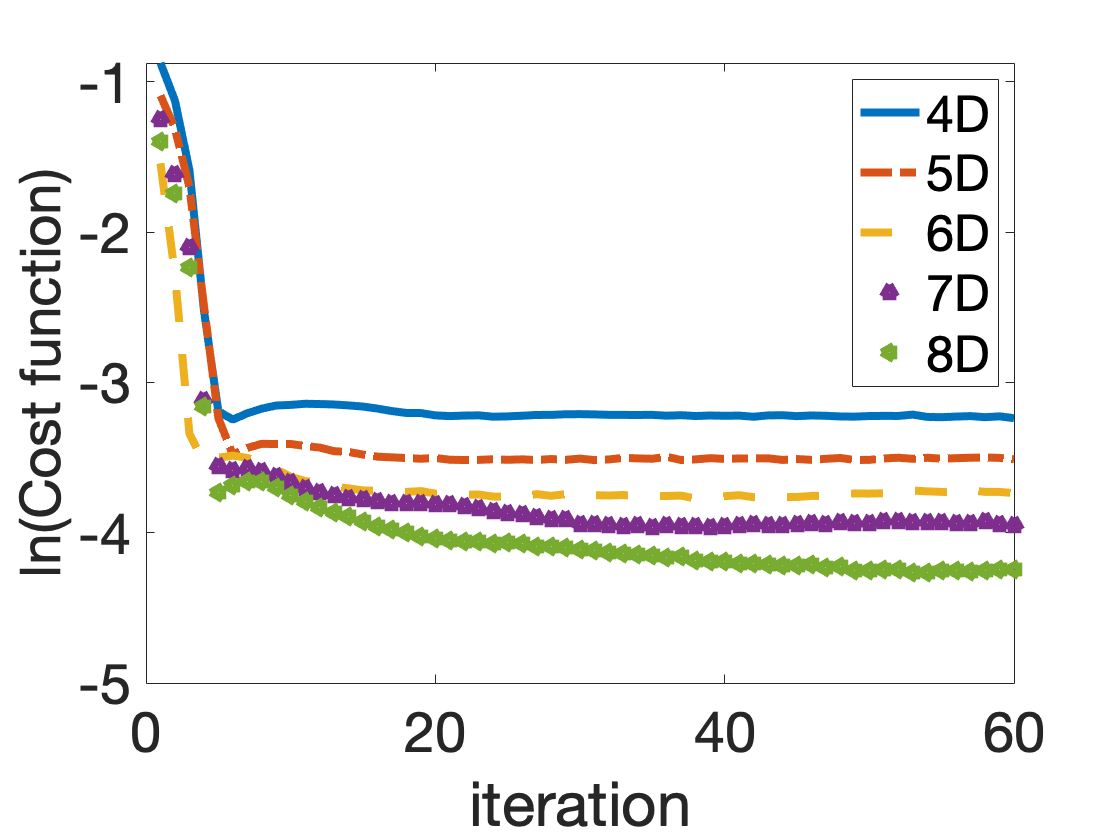}
	}
	\subfigure[ln(Edge error)]{    \label{fig:ndcase:b}
		\includegraphics[width=2.7in]{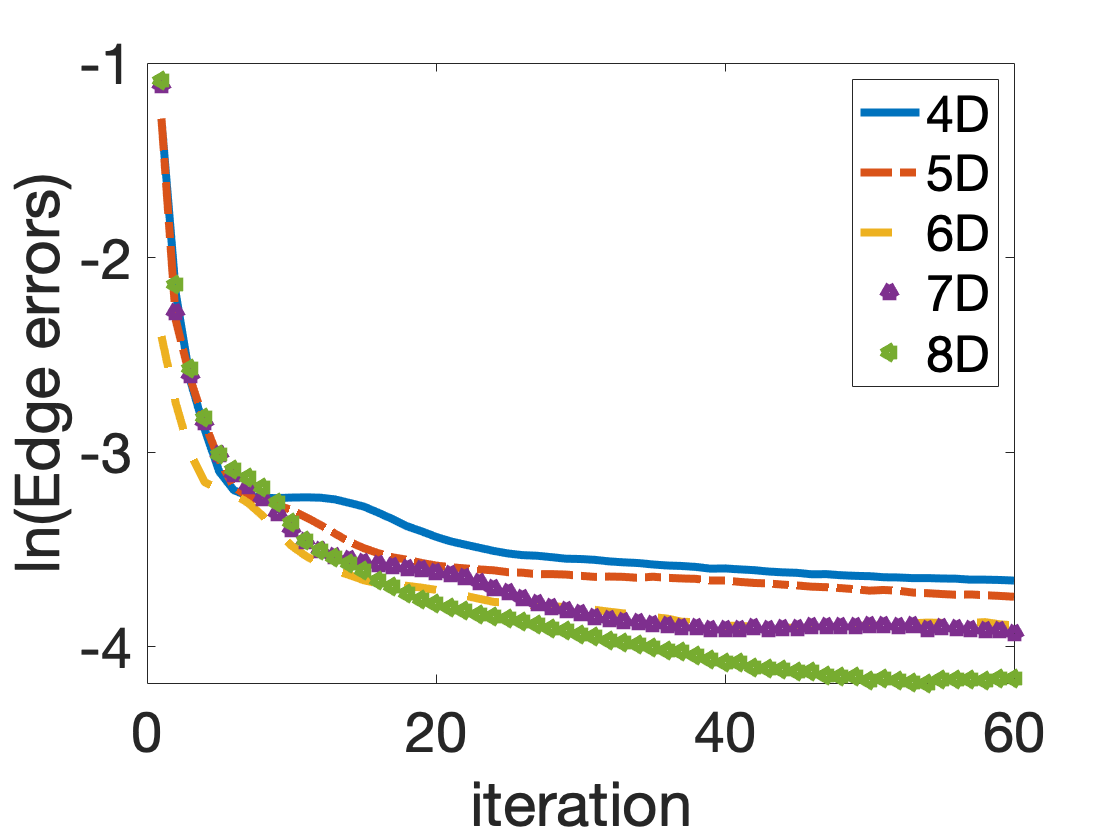}
	}
		\caption{The history of cost functions(Left) and edge errors(Right) in 4-8d problems with 49661 test points}
	\label{fig:ndcase} %% label for entire figure
\end{figure}

\begin{example}[4-8D cases]
\label{em:nd}
The equation of nD cases are given as
\begin{equation}
\begin{cases}
\Delta u(\mathbf{x}) = 0 & \mathbf{x} \ \ $in$  \ \ \Omega\\
u(\mathbf{x}) = x_1 + x_2 + \dots + x_n & \mathbf{x} \ \ $on$ \ \ \Gamma\\
\frac{\partial u(\mathbf{x})}{\partial \mathbf{n}} = \displaystyle\sum_i n_i& \mathbf{x} \ \ $on$  \ \ \Gamma
\end{cases}
\end{equation}

Let domain $\Omega$ be $x_1^2 + x_2^2 + \dots + x_n^2 < 0.25$ and boundary $\Gamma$ be $ x_1^2 + x_2^2 + \dots + x_n^2  = 0.25, x_n \in [-0.5,-0.25]\cup[0,0.25]$.
\end{example}
There are 10000 and 2500 points randomly sampling in $\Omega$ and $\Gamma$, respectively. Various dimensions(4-8) are studied to verify the accuracy of ANN method for Cauchy inverse problem. The training ended after 30000 steps, where steps = $10^{-4}$. Fig.~\ref{fig:ndcase} and Table \ref{fig:nderror} show the corresponding results.
%\begin{figure}[H]
%\centering
%\subfigure[4d case]{    \label{fig:ndcase:a}
%		\includegraphics[width=2.7in]{Perf4d}
%	}
%	\subfigure[5d case]{    \label{fig:ndcase:b}
%		\includegraphics[width=2.7in]{Perf5d}
%	}
%	\subfigure[6d case]{  \label{fig:ndcase:c}
%		\includegraphics[width=2.7in]{Perf6d}
%	}
%		\subfigure[7d case]{  \label{fig:ndcase:c}
%		\includegraphics[width=2.7in]{Perf7d}
%	}
%	\caption{The test errors in 4-7d problems with 49661 test points}
%	\label{fig:ndcase} %% label for entire figure
%\end{figure}

\begin{table}[h!]
\begin{center}
\begin{tabular}{|c|c|c|c|c|c|}
\hline
\ &4D case
&5D case            
&6D case
&7D case
&8D case \\
\hline
Training errors
&\raisebox{-.1\height}{0.0257}              
&\raisebox{-.1\height}{0.0237}  
&\raisebox{-.1\height}{0.0204}
&\raisebox{-.1\height}{0.0193} 
&\raisebox{-.1\height}{0.0142} \\
\hline
Test errors
&\raisebox{-.1\height}{0.0394}              
&\raisebox{-.1\height}{0.0302}  
&\raisebox{-.1\height}{ 0.0237}
&\raisebox{-.1\height}{0.0196} 
&\raisebox{-.1\height}{0.0153} \\
\hline
\end{tabular}
\end{center}
\caption{Errors in nd-case after 30000 iterations}
\label{fig:nderror}
\end{table}
As can be seen, there is a good agreement between the present result and the reference solution. The above results show that our method is worked for time-independent Cauchy inverse problem in high spatial dimension cases, and we discover more properties like stability about our method.
\subsection{Stability with noisy data and singular domain}
\label{subsec:prop}
The key point of numerical methods for Cauchy inverse problem is to treat the ill-posedness. Numerical experiments with 
noise on edge data are presented to show the stability of ANN method. The parameters and input data are same as example \ref{em:2d_ori} with two level of noise $\delta$($\ \% 1, \% 0.1$)., which is illustrated in Tab.~\ref{fig:noise}. We write down the history of cost function after 100000 iterations with noise $\% 1$. As can be seen, computed solution will always be convergent with some noise, illustrating that ANN method for Cauchy inverse problem is a good method to deal with ill-posedness.
\begin{table}[h!]
\begin{center}
\begin{tabular}{|c|c|c|c|}
\hline
$(\% 0.1)$ noise
&\raisebox{-.5\height}{\includegraphics[width=1.7in]{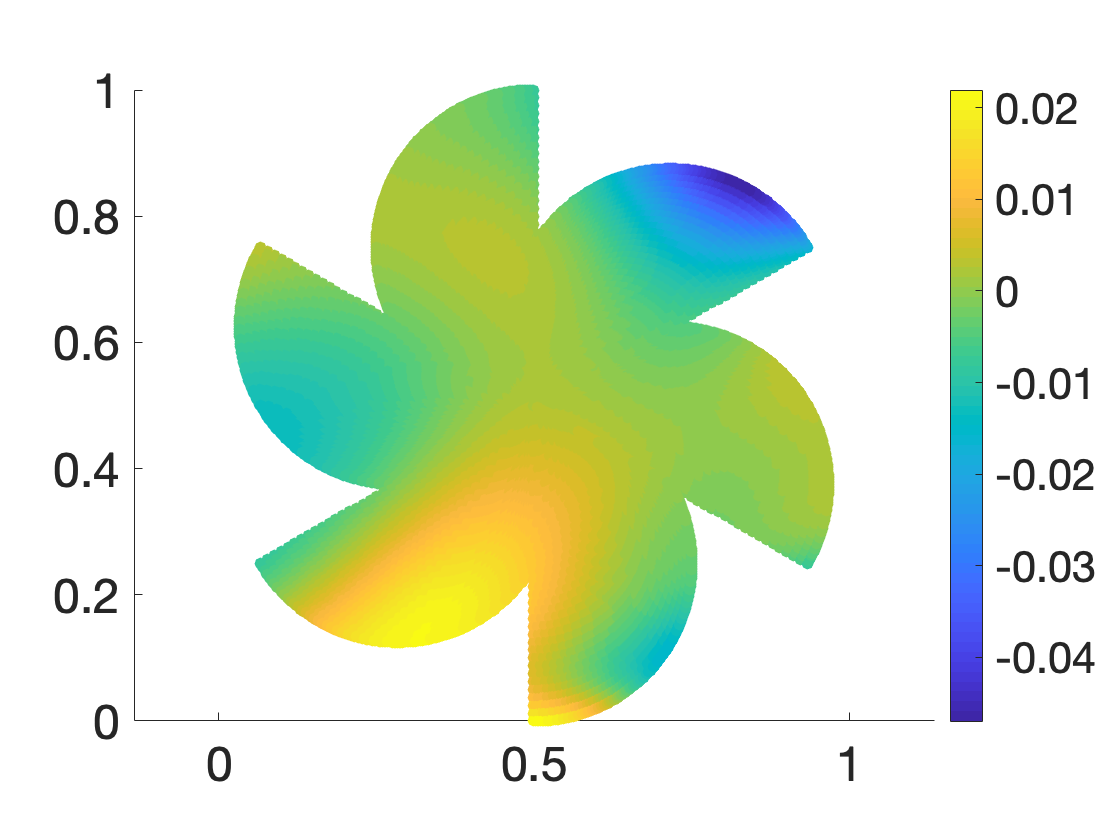} }              
&\raisebox{-.5\height}{\includegraphics[width=1.7in]{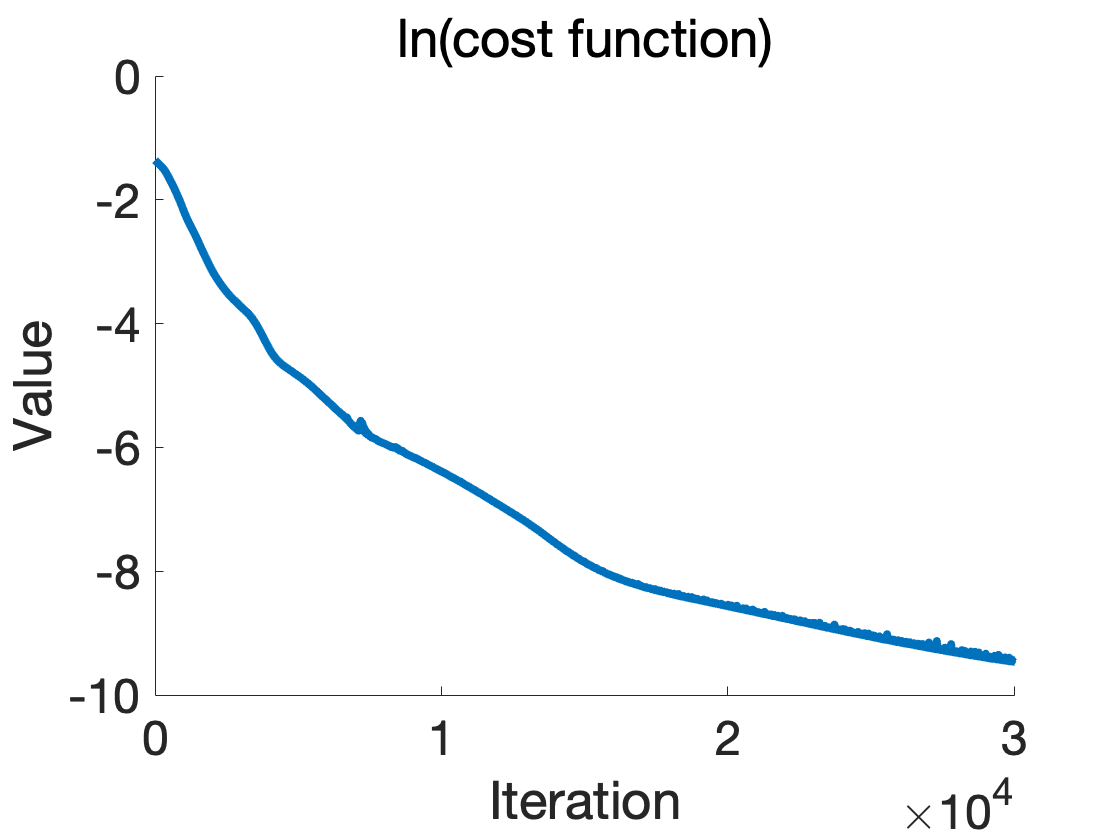}}  
&\raisebox{-.5\height}{\includegraphics[width=1.7in]{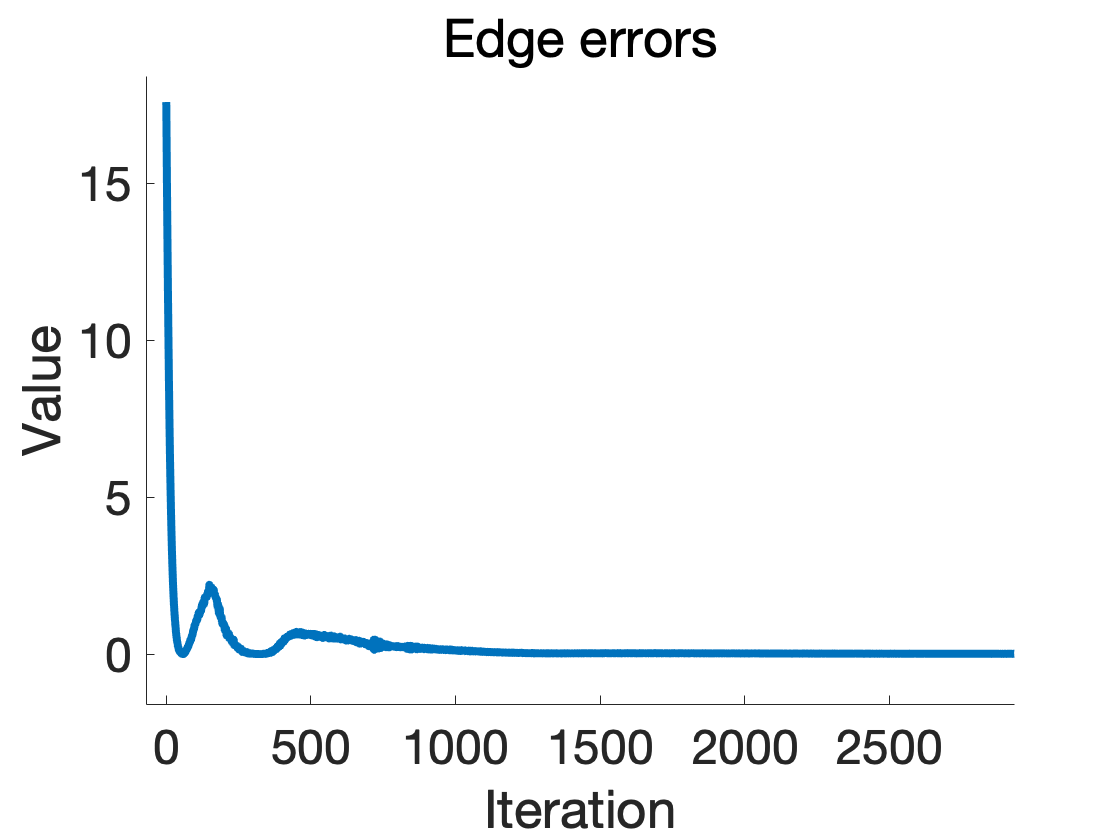}}  \\
\hline
$(\% 1)$ noise
&\raisebox{-.5\height}{\includegraphics[width=1.7in]{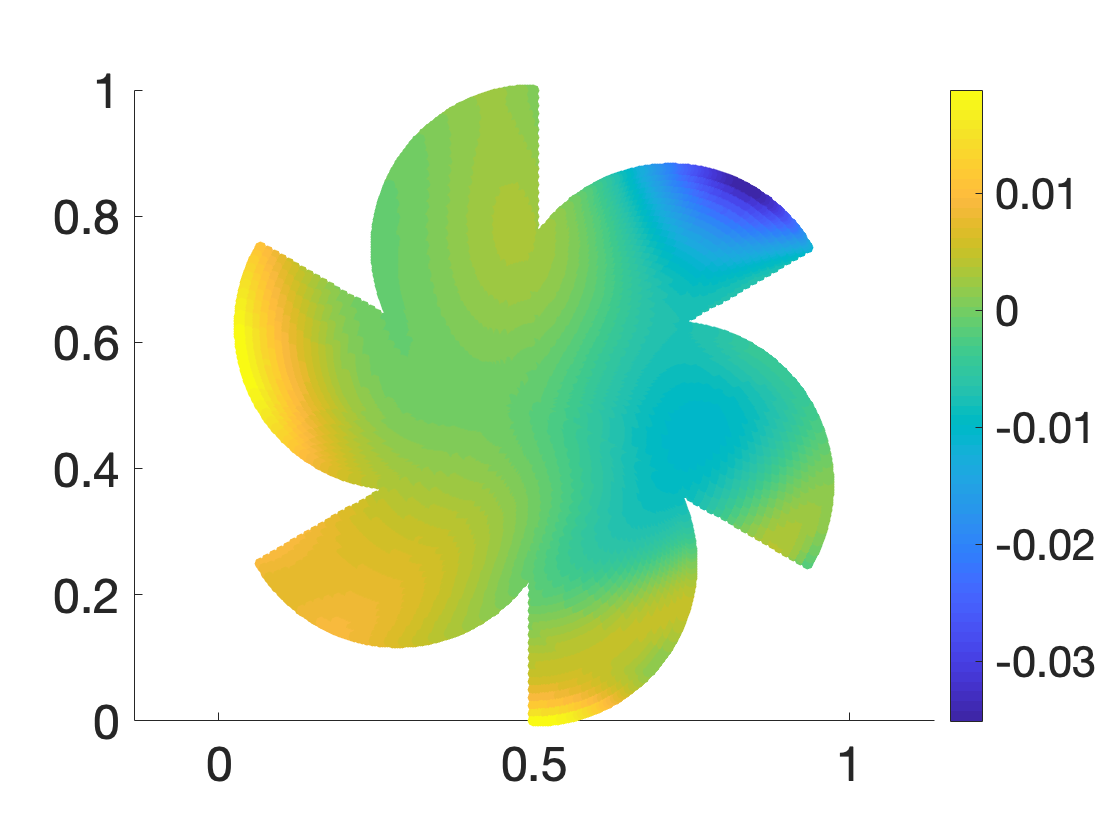}}                
&\raisebox{-.5\height}{\includegraphics[width=1.7in]{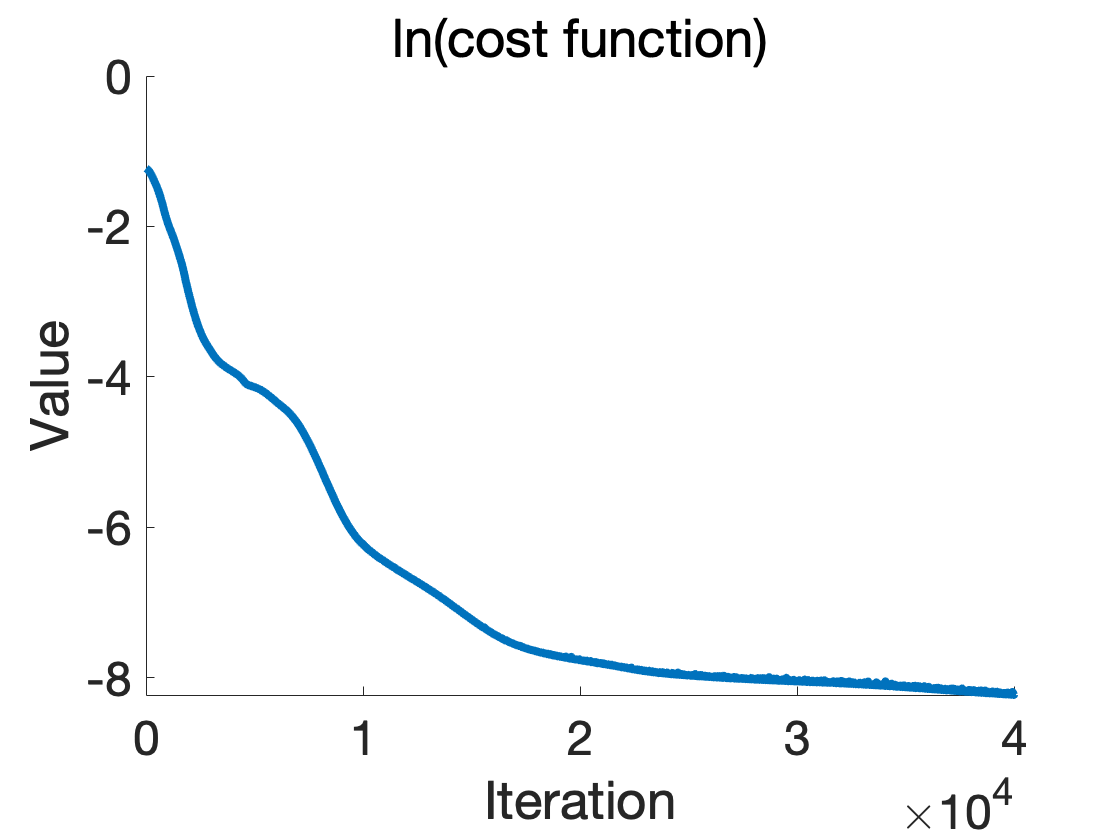}}  
&\raisebox{-.5\height}{\includegraphics[width=1.7in]{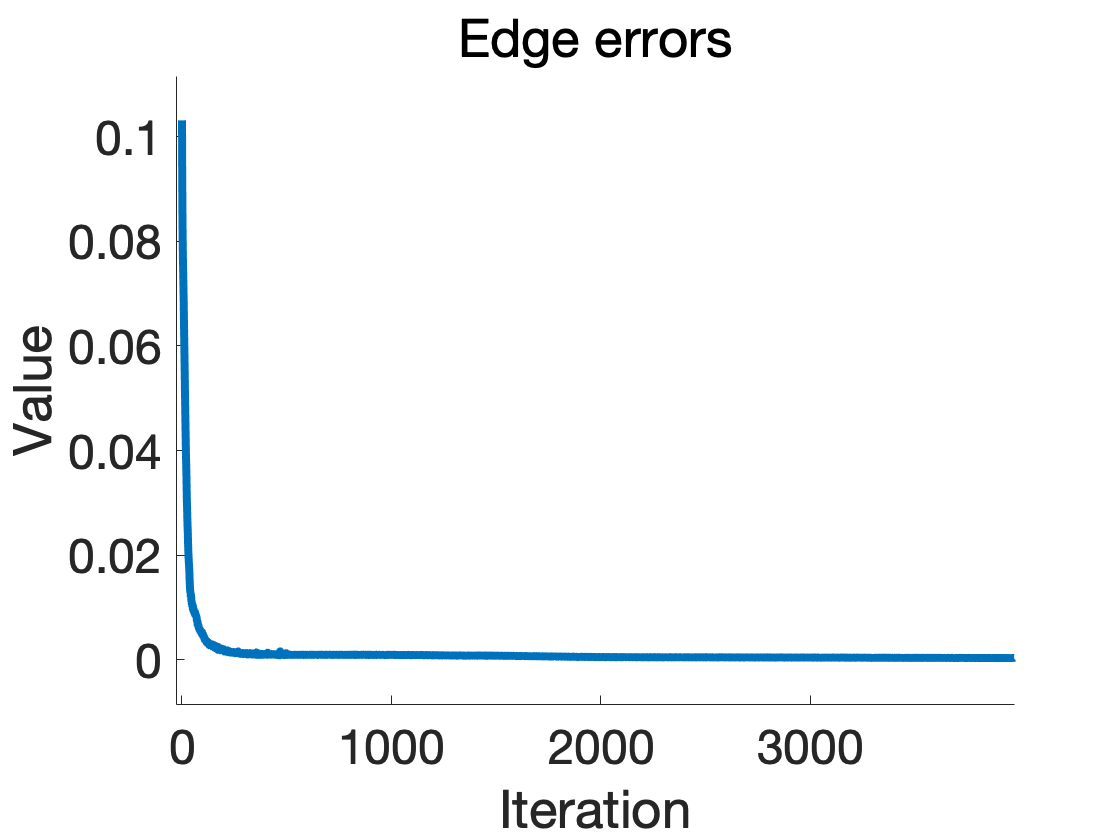}}  \\
\hline

\end{tabular}
\end{center}
\caption{(Left)errors  between exact and computed solution;(Middle)The convergence history of cost function;(Right)edge error between exact and computed solution during iterations}
\label{fig:noise}
\end{table}

\begin{figure}[H]
\centering
\includegraphics[height=5.5cm]{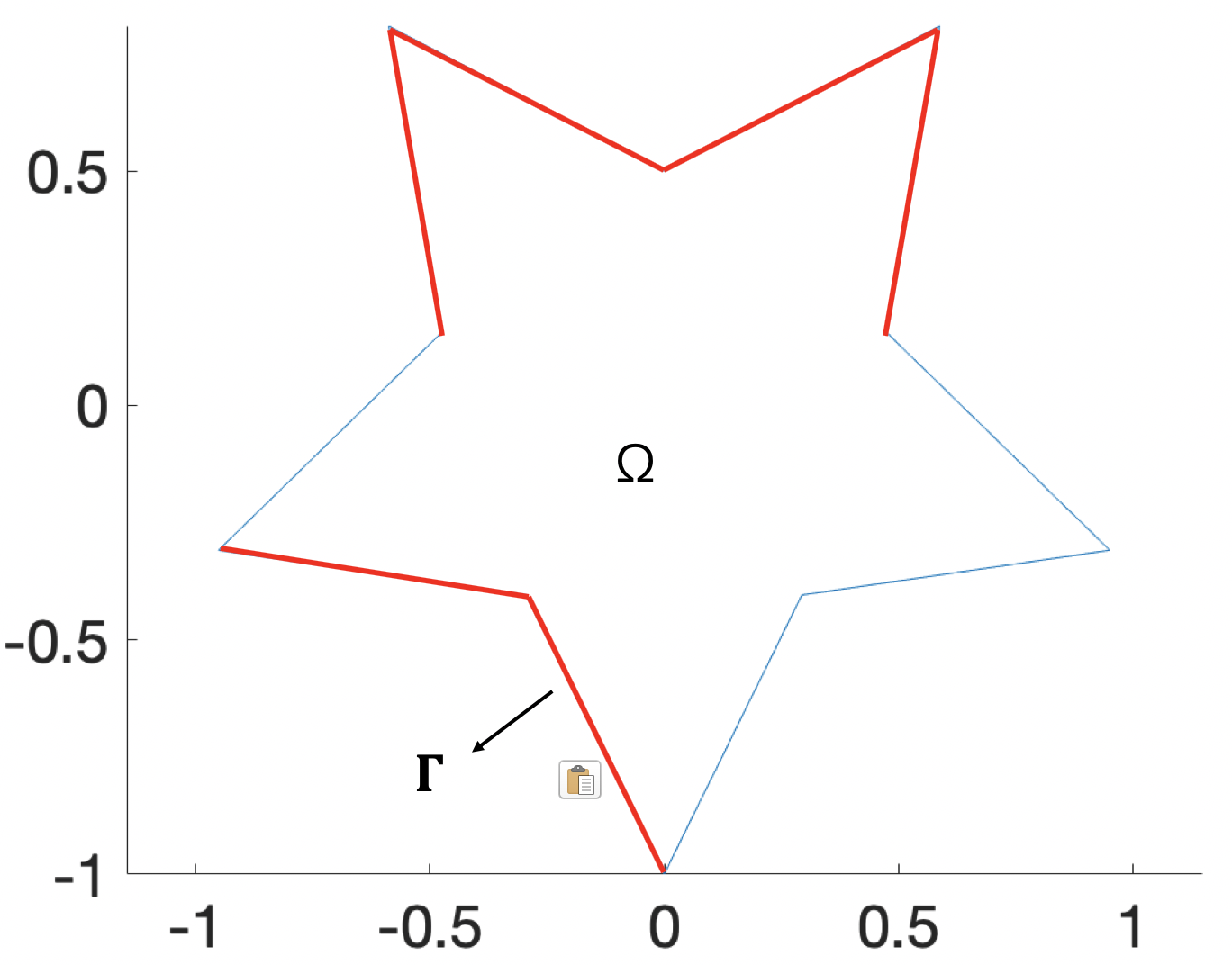}\qquad
\includegraphics[height=5.5cm]{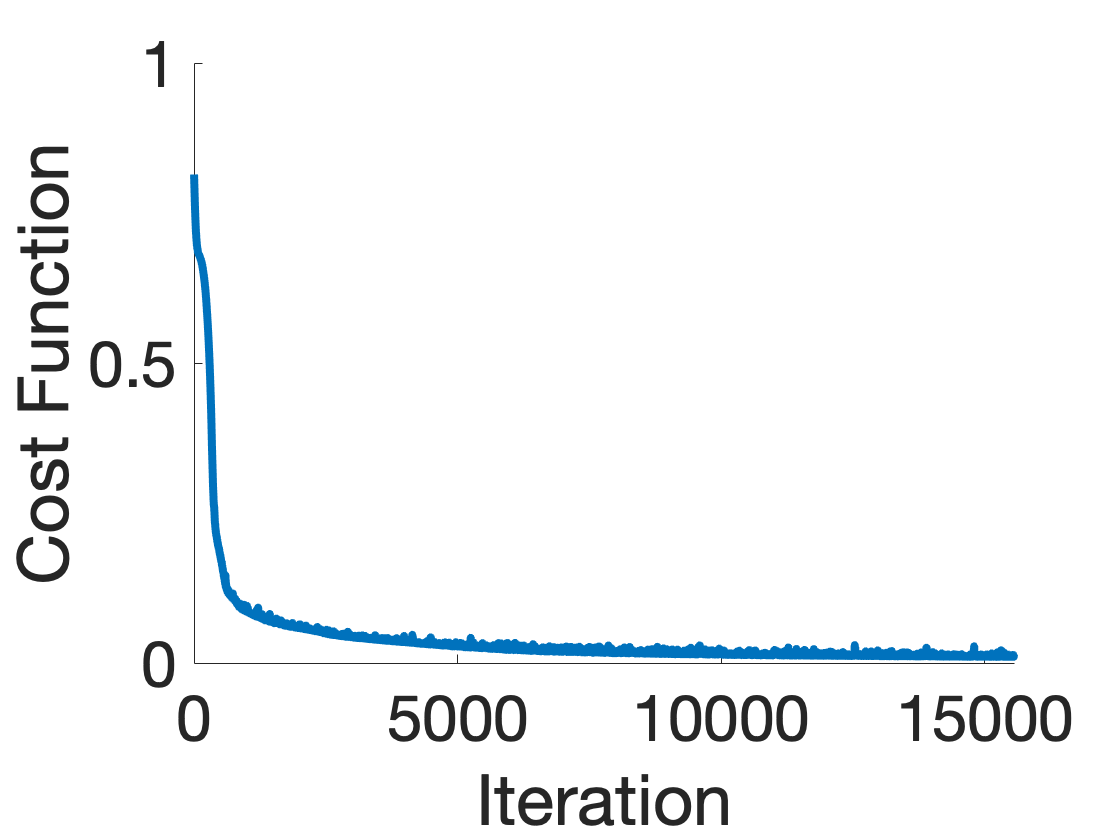}\qquad
\caption{(Left)The design area and boundary conditions of example \ref{em:2d_five}; (Right)The convergence history of cost function during iteration where steps = $5\times10^{-4}$.}
\label{fig:design_five}
\end {figure}

Moreover, it is well-known that some method for Cauchy inverse problem is sensitive with singular area, and the following example \ref{em:2d_five} will show the sensitivity about ANN method on a singular area. 
\begin{example}[domain with singularity]
\label{em:2d_five}
The equation of such problem is given as
\begin{equation}
\begin{cases}
\Delta u(\mathbf{x}) = 0 & \mathbf{x} \ \ $in$  \ \ \Omega\\
u(\mathbf{x}) = e^{x_1}sin(x_2) & \mathbf{x} \ \ $on$ \ \ \Gamma\\
\frac{\partial u(\mathbf{x})}{\partial \mathbf{n}} = [e^{x_1}sin(x_2), e^{x_1}cos(x_2)]*\mathbf{n}& \mathbf{x} \ \ $on$  \ \ \Gamma
\end{cases}
\end{equation}
We choose a singular domain $\Omega$ and boundary $\Gamma$ which is displayed in Fig. \ref{fig:design_five}(Left). 
\end{example}

\begin{figure}[H]
\centering
\includegraphics[height=5.6cm]{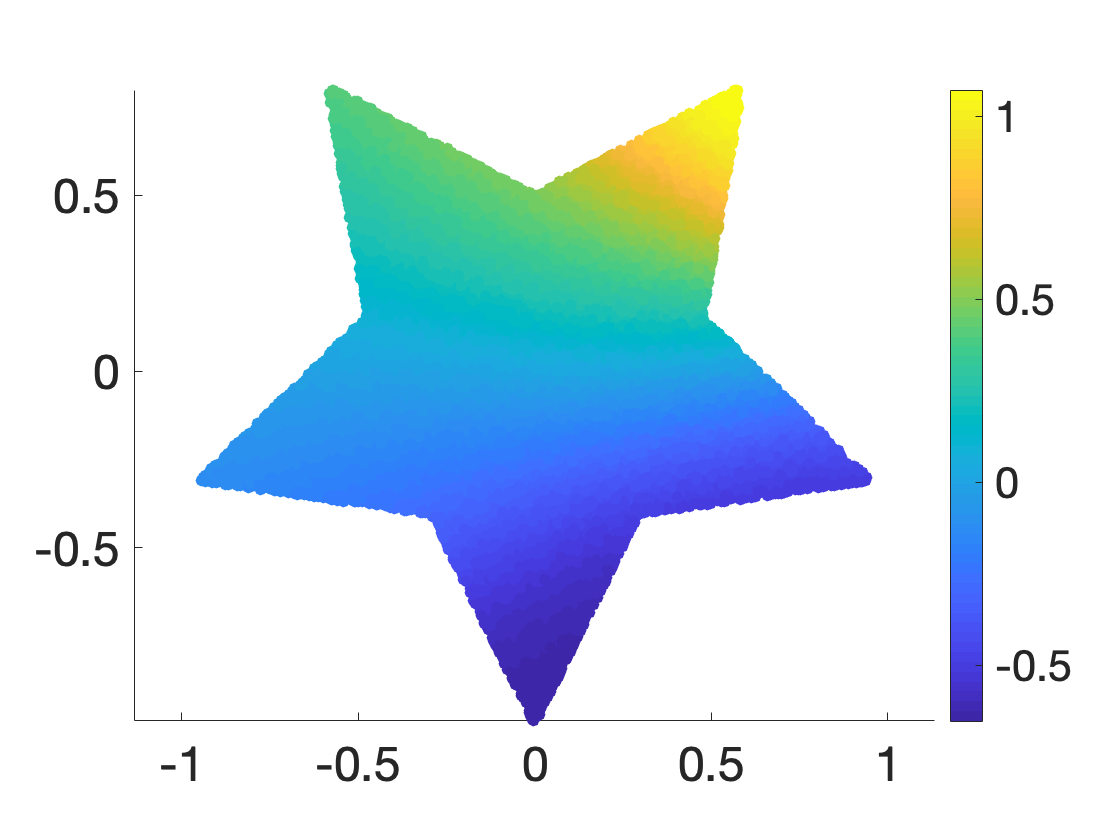}\qquad
\includegraphics[height=5.6cm]{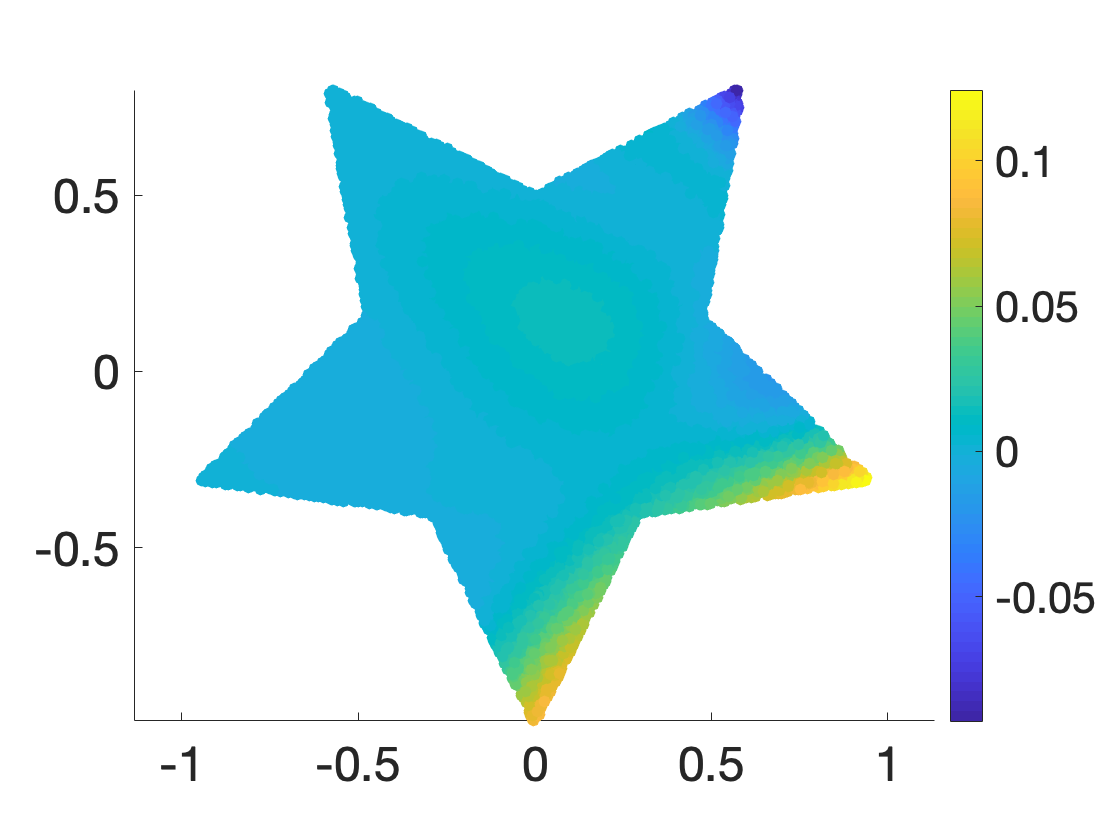}
\caption{(Left) The ANN solution to example \ref{em:2d_ori}; (Right) the errors  between exact and computed solution.}
\label{fig:approach_five}
\end {figure}

There are 3000 and 1000 points randomly sampling in $\Omega$ and $\Gamma$, respectivily. Let us choose the same structure of neural networks as example \ref{em:2d_ori}. The level of noise $\delta$ is set to be $\% 1$. Fig.~\ref{fig:design_five}(Right) shows the convergence history of the cost function  and errors between exact and computed solution are presented in Fig.~\ref{fig:approach_five}.

Graphically, there is a good agreement between the numerical and exact solutions, showing that ANN method for Cauchy inverse problem is insensitive with singular area. Moreover, the stability of it is verified. 

\subsection{The influence of depth and width of networks}
\begin{figure}[H]
\centering
\includegraphics[height=5.5cm]{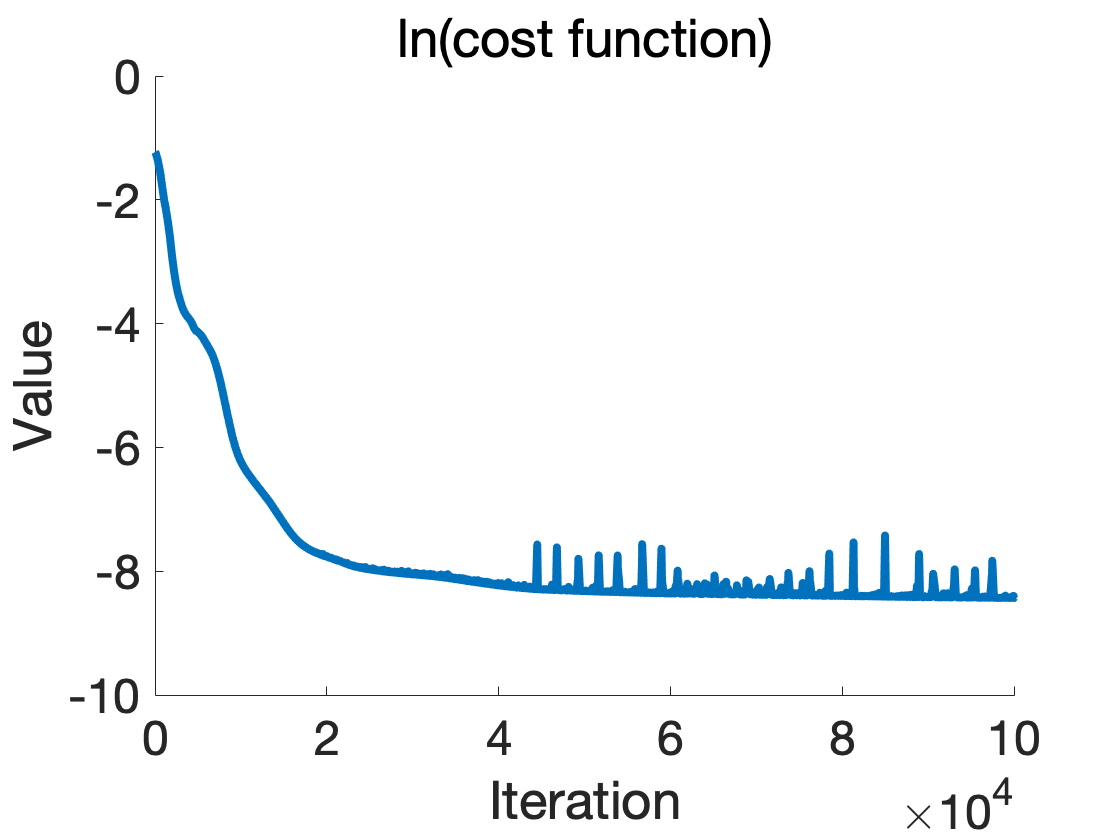}\qquad
\includegraphics[height=5.5cm]{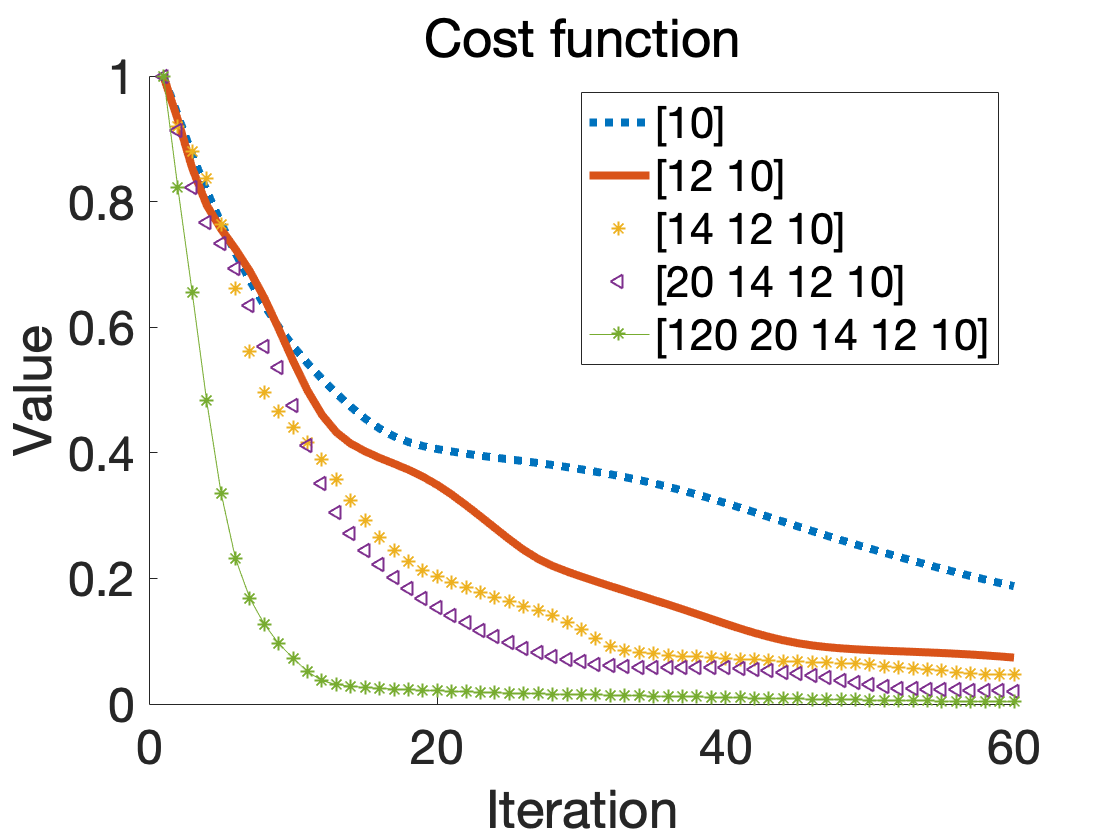}\qquad
\caption{(Left)The errors between exact and computed solution during iteration; (Right)The convergence history of cost function during iteration where steps = $10^{-4}$.}
\label{fig:width_layer}
\end {figure}

The influence of structures of ANN is an important issue which many researchers concern, at last we did some experiments on different hidden layers and neurons to discover it.  Firstly five neural networks with different hidden layers are tested with steps  $= 10^{-4}$ and 30000 iterations. Fig.~\ref{fig:width_layer}(Right) shows the corresponding cost function during iterations and Tab.~\ref{fig:hiddenlayers} shows errors(with training and testing) of five different hidden layers after 30000 iterations. As can be seen, the values of cost function decent faster with the parameters grows, and networks with more neurons and hidden layers have better approximation to the solution of PDEs, which satisfied theorems in section \ref{sec:analysis}. Secondly the influence of multi and single layers with same amounts of parameters is also tested with 30000 iterations. The result is displayed in Table.~\ref{fig:multilayer}, illustrating that approximation of multi hidden layers networks is always better than single hidden layers with same amount of parameters. It is worth to mention that when the neurons of one layer are too large(like 2968 parameters in Table \ref{fig:multilayer}), the method can not convergent with some initial $\mathbf{w}$ and $\mathbf{b}$ because of the value of sigmoid function will always approach to $1$. 
 
\begin{table}[h!]
\begin{center}
\begin{tabular}{|c|c|c|c|c|c|}
\hline
numbers of hidden layers 
&\raisebox{-.1\height}{1}              
&\raisebox{-.1\height}{2}  
&\raisebox{-.1\height}{3}
&\raisebox{-.1\height}{4} 
&\raisebox{-.1\height}{5} \\
\hline
Cost Function value
&\raisebox{-.1\height}{0.1118}              
&\raisebox{-.1\height}{0.04897}  
&\raisebox{-.1\height}{0.02745}
&\raisebox{-.1\height}{0.01484} 
&\raisebox{-.1\height}{0.002148} \\
\hline
Test errors ($10^5$ points)
&\raisebox{-.1\height}{0.08950}              
&\raisebox{-.1\height}{0.07075}  
&\raisebox{-.1\height}{0.06495}
&\raisebox{-.1\height}{0.06301} 
&\raisebox{-.1\height}{0.03483} \\
\hline
\end{tabular}
\end{center}
\caption{The cost function and test errors after 30000 iterations}
\label{fig:hiddenlayers}
\end{table}

\begin{table}[h!]
\begin{center}
\begin{tabular}{|c|c|c|c|c|c|}
\hline
Numbers of parameters 
&\raisebox{-.1\height}{10}              
&\raisebox{-.1\height}{120}  
&\raisebox{-.1\height}{288}
&\raisebox{-.1\height}{568} 
&\raisebox{-.1\height}{2968} \\
\hline
Single hidden layer case
&\raisebox{-.1\height}{0.1118}              
&\raisebox{-.1\height}{0.1053}  
&\raisebox{-.1\height}{0.0947}
&\raisebox{-.1\height}{0.0705} 
&\raisebox{-.1\height}{do not convergent} \\
%&\raisebox{-.1\height}{\textit{0.5831}} \\
\hline
Multi hidden layer case
&\raisebox{-.1\height}{0.1118}              
&\raisebox{-.1\height}{0.04897}  
&\raisebox{-.1\height}{0.02745}
&\raisebox{-.1\height}{0.01484} 
&\raisebox{-.1\height}{0.002148} \\
\hline
\end{tabular}
\end{center}
\caption{The cost function value of single and multi hidden layers after 30000 iterations}
\label{fig:multilayer}
\end{table}
With above numerical examples, we  can conclude that ANN method for Cauchy inverse problem with multi-hidden layers is a  stable method on various aspects, including well-posedness with noise on boundary and insensitivity on singular area.

%\begin{table}[h!]
%\begin{center}
%\begin{tabular}{|c|c|c|}
%\hline
%$z = 0.8536$\\
%&\raisebox{-.5\height}{\includegraphics[width=2in]{z1_3d_app.png} }              
%&\raisebox{-.5\height}{\includegraphics[width=2in]{z1_3d_err.png}}  \\
%\hline
%$z = 0.75$\\
%&\raisebox{-.5\height}{\includegraphics[width=2in]{z2_3d_app.png}}                
%&\raisebox{-.5\height}{\includegraphics[width=2in]{z2_3d_err.png}}  \\
%\hline
%$z = 0.5$\\
%&\raisebox{-.5\height}{\includegraphics[width=2in]{z3_3d_app.png}}                
%&\raisebox{-.5\height}{\includegraphics[width=2in]{z3_3d_err.png}}  \\
%\hline
%$z = 0.25$\\
%&\raisebox{-.5\height}{\includegraphics[width=2in]{z4_3d_app.png}}                
%&\raisebox{-.5\height}{\includegraphics[width=2in]{z4_3d_err.png}}  \\
%\hline
%$z = 0.1464$\\
%&\raisebox{-.5\height}{\includegraphics[width=2in]{z5_3d_app.png}}                
%&\raisebox{-.5\height}{\includegraphics[width=2in]{z5_3d_err.png}}  \\
%\hline
%\end{tabular}
%\end{center}
%\caption{(Left)The ANN solution in one section;(Right)errors  between exact and computed solution}
%\label{fig:approach_3d}
%\end{table}
\section{Conclusions}
\label{sec:conclusion}
We have investigated an artificial neural network approximation for solving Cauchy inverse problem. The inputs of network come from physical models, including the initial and boundary datas with noise, rather than any exact or experiment solutions. Denseness and m-denseness of networks are extended to multi-hidden layers in this paper, so that we can prove the equivalence between Cauchy inverse problem and its ANN approximation. Various numerical examples show that our method is  effective for both time-dependent and time-independent cases. It is worth noting that ANN methods for Cauchy inverse problem have some stable properties, such as approaching data with noise $\delta$, solving high-dimension cases with  lower computational cost and being insensitive on a singular domain. Moreover, the influence about the number of layers and neurons with neural network are also discovered by numerical examples that networks with more hidden layers have better performance. 

In addition, it would be of interest to extend neural network method into other inverse problem or parameter design problem. Otherwise, more novel structures and some convergence analysis of deep learning for solving PDE inverse problems is also an important question. We leave these questions for future work.

\section*{Appendix}
\begin{appendix}
In the appendix we introduce some computation of the derivation and back propagation based on the chain rule. Artificial neural network can be considered as an approximated function of input data.  As is known that operator $\mathcal{L}$ can be represented by $n_{th}$-order derivations of space $\Omega$. At first, let us compute some derivations of ANN.
\section{Back propagation of $n_{th}$-order derivations with ANN}
\label{App:derivation}
\subsection{Back propagation with ANN}
Back propagation of 0-order derivation is the original backpropagation with ANN. The back propagation with 0-order derivation($\frac{\partial\mathbf{y}^{L+1}}{\partial\mathbf{w}^1} , \frac{\partial\mathbf{y}^{L+1}}{\partial\mathbf{b}^l}$) can be computed by taking derivatives of $\mathbf{z}^l$ and $y^l$. Following the chain rule, we know that
\begin{equation}
\label{eq:backpropagation}
\begin{split}
\frac{\partial\mathbf{y}^{L+1}}{\partial\mathbf{w}^l}  & = \frac{\partial\mathbf{y}^{L+1}}{\partial\mathbf{z}^{l}} * \frac{\partial\mathbf{z}^{l}}{\partial\mathbf{w}^l} = \frac{\partial\mathbf{y}^{L+1}}{\partial\mathbf{z}^{l}} * \mathbf{y}^{l-1},\\
\frac{\partial\mathbf{y}^{L+1}}{\partial\mathbf{b}^l} & = \frac{\partial\mathbf{y}^{L+1}}{\partial\mathbf{z}^l} * \frac{\partial\mathbf{z}^{l}}{\partial\mathbf{b}^l} = \frac{\partial\mathbf{y}^{L+1}}{\partial\mathbf{z}^l}. 
\end{split}
\end{equation}
By defining $\delta^l := \frac{\partial\mathbf{y}^{L+1}}{\partial\mathbf{z}^l}, l = 1,2,\dots,L+1$, equation \eqref{eq:layerin} yields that 
\begin{equation*}
\label{eq:delta}
\begin{split}
&\delta^{L+1}  = \sigma_{L+1}'(\mathbf{z}^{L+1}),\\
&\delta^l  = \left(\mathbf{w}^{l+1}\right)^T * \delta^{l+1}\odot\sigma_l'({\mathbf{z}^l}), \ \ \ l = 1,2,\dots,L;
\end{split}
\end{equation*}
where $\odot$ represent multiplication by corresponding elements of two vectors. 

\subsection{Back propagation of 1-order derivation with ANN}
Assume that input $\mathbf{x} \in \Omega\times\mathcal{T} \subset\mathbb{R}^d\times\mathbb{R}^1$, the back propagation of 1-order derivation ($\frac{\partial^2\mathbf{y}^{L+1}}{\partial\mathbf{w}^1\partial x_i} , \frac{\partial^2\mathbf{y}^{L+1}}{\partial\mathbf{b}^l\partial x_i}$) can be obtained in the following equations:

\begin{eqnarray}
 \frac{\partial^2\mathbf{y}^{L+1}}{\partial\mathbf{x}_i\partial\mathbf{w}^l} & =  &\frac{\partial\delta^l}{\partial\mathbf{x}_i}*(\mathbf{y}^{l-1})^T + \delta^l * \frac{\partial(\mathbf{y}^{l-1})^T}{\partial\mathbf{x}_i}\label{eq:backpro_1_w}\\
  \frac{\partial^2\mathbf{y}^{L+1}}{\partial\mathbf{x}_i\partial\mathbf{b}^l} & = & \frac{\partial\delta^l}{\partial\mathbf{x}_i}\label{eq:backpro_1_b},
 \end{eqnarray}

where $\frac{\partial\delta^l}{\partial\mathbf{x}_i}$ can be computed as 
\begin{equation*}
\begin{split}
\frac{\partial\delta^{L+1}}{\partial\mathbf{x}_i} &= \sigma_{L+1}''(\mathbf{z}^{L+1}) \odot \frac{\partial\mathbf{z}^{L+1}}{\partial x_i}\\
\frac{\partial\delta^{l}}{\partial\mathbf{x}_i} &=  \left(\mathbf{w}^{l+1}\right)^T * \frac{\partial\delta^{l+1}}{\partial x_i}\odot\sigma_l'({\mathbf{z}^l}) +  \left(\mathbf{w}^{l+1}\right)^T * \delta^{l+1}\odot\sigma_l''({\mathbf{z}^l})\odot\frac{\partial\mathbf{z}^l} {\partial x_i}
\end{split}
\end{equation*}
The analytic formula of $\frac{\partial \mathbf{y}^{l}}{\partial x_i}$ and $\frac{\partial \mathbf{z}^{l}}{\partial x_i}$ in the above equations can be represented one by one with layers as:

\begin{equation*}
\begin{split}
&\frac{\partial \mathbf{z}^{1}}{\partial x_i} = \mathbf{w}^1\\
&\frac{\partial \mathbf{z}^{l+1}}{\partial x_i} = \mathbf{w}^{l+1} * \frac{\partial \mathbf{z}^{l}}{\partial x_i}\odot\sigma_{l}'(\mathbf{z}^l)\\
&\frac{\partial \mathbf{y}^{l}}{\partial x_i} = \frac{\partial \mathbf{z}^{l}}{\partial x_i} \odot \sigma_l'(\mathbf{z}^l)
\end{split}
\end{equation*}

Moreover, when $\Omega \subset \mathbb{R}^d, d > 3$, it's convenient to compute with matrix form. Before that we introduce some notations at first.

For a vector $\mathbf{f} = [f_1, f_2, \dots, f_m]$, the Jacobi matrix of $\mathbf{f}$ is defined as
\begin{equation*}
G(\mathbf{f}) = \begin{bmatrix}\frac{\partial f_1}{\partial x_1} &\frac{\partial f_1}{\partial x_2}  &\dots&\frac{\partial f_1}{\partial x_d}&\frac{\partial f_1}{\partial t}\\ \frac{\partial f_2}{\partial x_1}&\frac{\partial f_2}{\partial x_2}&\dots&\frac{\partial f_2}{\partial x_d}&\frac{\partial f_2}{\partial t}\\ \vdots&\vdots&\ddots&\vdots\\ \frac{\partial f_m}{\partial x_1}&\frac{\partial f_m}{\partial x_2}&\dots&\frac{\partial f_m}{\partial x_d}&\frac{\partial f_m}{\partial t}\end{bmatrix}
\end{equation*}
and for a vector $\mathbf{u} = [u_1, u_2, \dots, u_n]$, let the diagonal matrices be 
 \begin{equation*}
 diag(\mathbf{u}) = \begin{bmatrix} u_1 & 0 &\dots&0\\0&u_2&\dots&0\\ \vdots&\vdots&\ddots&\vdots\\0&0&\dots&u_n\end{bmatrix}
 \end{equation*}
 With the above definition there obviously establish in matrix form that 
  \begin{equation}
 \label{eq:Ddelta}
 \begin{split}
 & G(\delta^{L+1}) = diag(\sigma_{L+1}''(\mathbf{z}^{L+1}))*G(\mathbf{z}^{L+1}),\\
 &G(\delta^l)  = diag(\sigma_l'(\mathbf{z}^l))*(\mathbf{w}^{l+1})^T*G (\delta^{l+1}) + diag(\sigma_l''(\mathbf{z}^l))*(\mathbf{w}^{l+1})^T*\delta^{l+1}\odot G(\mathbf{z}^l),\\
 &\frac{\partial\delta^l}{\partial\mathbf{x}_i} = G(\delta^l)[i],
 \end{split}
 \end{equation}
where
\begin{equation*}
 \label{eq:Grad_y}
 \begin{split}
 &G(\mathbf{z}^1) = \mathbf{w}^1*I_{d+1},\\
& G(\mathbf{z}^l) =  \mathbf{w}^{l} * diag(\sigma_{l-1}(\mathbf{z}^{l-1}) )* G(\mathbf{z}^{l-1}),\ \ l = 2,3,\dots,L+1 \\
& G(\mathbf{y}^{l}) = diag(\sigma_{l}'(\mathbf{z}^{l}))*G(\mathbf{z}^{l}),\ \ l = 1,2,\dots,L+1 
 \end{split}
 \end{equation*}
 where $I_d$ is the $d\times d$ identity matrix. 
 
 \subsection{Back propagation with 2-order derivation of ANN}
Assume the input $\mathbf{x} \in \Omega\times\mathcal{T} \subset\mathbb{R}^d\times\mathbb{R}^1$, and back propagation of 2-order derivation ($\frac{\partial^3\mathbf{y}^{L+1}}{\partial\mathbf{w}^1\partial^2 x_i} , \frac{\partial^3\mathbf{y}^{L+1}}{\partial\mathbf{b}^l\partial^2 x_i}$) can be obtained in the following equations:

\begin{equation}
\label{eq:2_order_deri}
\begin{split}
\frac{\partial^3\mathbf{y}^{L+1}}{\partial\mathbf{w}^1\partial^2 x_i} &= \frac{\partial^2\delta^l}{\partial x_i^2} * (\mathbf{y}^{l-1})^T + 2\frac{\partial\delta^l}{\partial x_i} * \frac{\partial (\mathbf{y}^{l-1})^T }{\partial x_i} + \delta^l * \frac{\partial^2(\mathbf{y}^{l-1})^T}{\partial x_i^2}\\
\frac{\partial^3\mathbf{y}^{L+1}}{\partial\mathbf{b}^l\partial^2 x_i} &= \frac{\partial^2\delta^l}{\partial x_i^2},
\end{split}
\end{equation}
where $\frac{\partial^2\delta^l}{\partial x_i^2}$ can be computed as
\begin{equation*}
\begin{split}
\frac{\partial^2\delta^{L+1}}{\partial x_i^2} &= \sigma_{L+1}'''(\mathbf{z}^{L+1}) \odot \left(\frac{\partial\mathbf{z}^{L+1}}{\partial x_i}\right)^2 \\
\frac{\partial^2\delta^l}{\partial x_i^2} &= (\mathbf{w}^{l+1})^T * \left(\frac{\partial^2\delta^{l+1}}{\partial x_i^2} \odot \sigma_l'(\mathbf{z}^l) + 2 \frac{\partial\delta^{l+1}}{\partial x_i}\odot\sigma_l''(\mathbf{z}^l)\odot\frac{\partial\mathbf{z}^l}{\partial x_i}  \right) \\
& + (\mathbf{w}^{l+1})^T * \delta^l \odot\left(\sigma_l'''(\mathbf{z}^l) \odot \left(\frac{\partial\mathbf{z}^l}{\partial x_i}\right) + \sigma_l''(\mathbf{z}^l)\odot\frac{\partial^2\mathbf{z}^l}{\partial x_i^2} \right)
\end{split}
\end{equation*}

The analytic formula of $\frac{\partial^2 \mathbf{y}^{l}}{\partial x^2_i}$ and $\frac{\partial^2 \mathbf{z}^{l}}{\partial x^2_i}$ in the above euqations can be represented one by one with layers:
\begin{eqnarray*}
\frac{\partial^2 \mathbf{z}^{1}}{\partial x^2_i} &= &0 \\
\frac{\partial^2 \mathbf{z}^{l+1}}{\partial x^2_i} & = & \mathbf{w}^{l+1} * \left( \frac{\partial^2\mathbf{z}^l}{\partial x_i^2}\odot\sigma_l'(\mathbf{z}^l) + \left(\frac{\partial\mathbf{z}^l}{\partial x_i}\right)^2 \odot \sigma_l''(\mathbf{z}^l)\right)\\
\frac{\partial^2 \mathbf{y}^{l}}{\partial x^2_i} & = &  \frac{\partial^2\mathbf{z}^l}{\partial x_i^2}\odot\sigma_l'(\mathbf{z}^l) + \left(\frac{\partial\mathbf{z}^l}{\partial x_i}\right)^2 \odot \sigma_l''(\mathbf{z}^l)
\end{eqnarray*}

Moreover, when $\Omega \subset \mathbb{R}^d, d > 3$, it's convenient to compute with matrix form. We introduce some notations at first.

For a vector $\mathbf{u}$, $G^2(\mathbf{u})$ is defined as: 

\begin{equation*}
G^2(\mathbf{u}) = \begin{bmatrix} \frac{\partial^2u_1}{\partial x_1^2}&  \frac{\partial^2u_1}{\partial x_2^2}& \dots,& \frac{\partial^2u_1}{\partial x_d^2}\\   \frac{\partial^2u_2}{\partial x_1^2}&  \frac{\partial^2u_2}{\partial x_2^2}& \dots,& \frac{\partial^2u_2}{\partial x_d^2}\\ \vdots&\vdots&\ddots&\vdots \\ \frac{\partial^2u_n}{\partial x_1^2}&  \frac{\partial^2u_n}{\partial x_2^2}& \dots,& \frac{\partial^2u_n}{\partial x_d^2}\end{bmatrix}
\end{equation*}

The matrix form of 2-order derivations are shown as 
\begin{eqnarray}
& &G^2(\delta^{L+1}) = diag(\sigma_{L+1}^{(3)}) * (\nabla\mathbf{z}^{L+1}\odot  \nabla\mathbf{z}^{L+1})+diag(\sigma_{L+1}''(\mathbf{z}^{L+1}))*G^2(\mathbf{z}^{L+1})\label{eq:Lap_y}\\
& & \begin{split}
G^2(\delta^l) = & diag(\sigma_l'(\mathbf{z}^l))*(\mathbf{w}^{l+1})^T*G^2(\delta^{l+1}) + \\
&diag(\sigma_l''(\mathbf{z}^l))*(\mathbf{w}^{l+1})^T*(2\nabla\delta^{l+1})\odot\nabla\mathbf{z}^l + \\
&diag(\sigma_l^{(3)}(\mathbf{z}^l))*(\mathbf{w}^{l+1})^T*\delta^{l+1}\odot(\nabla\mathbf{z}^l\odot\nabla\mathbf{z}^l) + \\
&diag(\sigma_l''(\mathbf{z}^l))*(\mathbf{w}^{l+1})^T*\delta^{l+1}\odot J(\mathbf{z}^l), \ \ l = 1,2,\dots,L
\end{split}\\
& & \frac{\partial^2\delta^l}{\partial x_i^2} = G^2(\delta^l)[i],
\end{eqnarray}

There obviously establishes following equation for $G^2(\mathbf{z}^l)$ and $ G^2(\mathbf{z}^l)$: 
\begin{equation*}
\begin{split}
& G^2(\mathbf{z}^1) = \mathbf{w}^1*O_{d+1}, \\
& G^2(\mathbf{z}^l) = \mathbf{w}^l*(diag(\sigma_{l-1}''(\mathbf{z}^{l-1}))*(\nabla\mathbf{z}^{l-1}\odot\nabla\mathbf{z}^{l-1}) + diag(\sigma_{l-1}'(\mathbf{z}^{l-1})) * G^2(\mathbf{z}^{l-1})),\\
& G^2(\mathbf{y}^{l}) = diag(\sigma_{l}''(\mathbf{z}^{l}))*(\nabla\mathbf{z}^{l}\odot\nabla\mathbf{z}^{l}) + diag(\sigma_{l}'(\mathbf{z}^{l})) * G^2(\mathbf{z}^{L+1})),
\end{split}
\end{equation*}
where $O_d$ is the zero matrix of dimension $d\times d$.

\section{The analytic formula of back propagation for Cauchy inverse problem}
\label{App:backpro}
We use the result and notation in the last appendix \ref{App:derivation}. Firstly, we show the analytic formula of back propagation for problem \eqref{eq:parabolic}, and problem \eqref{eq:elliptic} is similar.

The back propagation for Dirichlet boundary condition and initial condition can be easily obtained in equation \eqref{eq:backpropagation}. Neumann boundary condition and the state equation is considered with problem \eqref{eq:parabolic}.
\subsection{back propagation for Neumann boundary condition}
We know that for time dependent problem, there establish that
\begin{equation}
 G(\mathbf{y}^{L+1})  =  \begin{bmatrix}\frac{\partial \mathbf{y}^{L+1}}{\partial x_1} &\frac{\partial \mathbf{y}^{L+1}}{\partial x_2}  &\dots&\frac{\partial \mathbf{y}^{L+1}}{\partial x_{d}}&\frac{\partial\mathbf{y}^{L+1}}{\partial t}\end{bmatrix} 
  \end{equation}
Neumann boundary condition can be represented by 
\begin{equation}
\label{eq:grap}
 \frac{\partial \mathbf{y}^{L+1}}{\partial\mathbf{n}} = \nabla \mathbf{y}^{L+1} * D^T(\mathbf{x}) = G( \mathbf{y}^{L+1})[1:d]* D^T(\mathbf{x}),
 \end{equation}
where $D(\mathbf{x}) = (\beta_{1},\beta_{2},\dots,\beta_{d})$ represents the direction at points $\mathbf{x}$.  Similarly the back propagation of Neumann boundary condition is shown as following
\begin{eqnarray}
 \frac{\partial^2 \mathbf{y}^{L+1}}{\partial\mathbf{n}\partial\mathbf{w}^l}& =& \frac{\partial\nabla \mathbf{y}^{L+1}}{\partial\mathbf{w}^l} * D^T(\mathbf{x}) = \nabla\frac{\partial \mathbf{y}^{L+1}}{\partial\mathbf{w}^l} * D^T(\mathbf{x}), \\
 \frac{\partial^2 \mathbf{y}^{L+1}}{\partial\mathbf{n}\partial\mathbf{b}^l} &=& \frac{\partial\nabla \mathbf{y}^{L+1}}{\partial\mathbf{b}^l} * D^T(\mathbf{x}) = \nabla\frac{\partial \mathbf{y}^{L+1}}{\partial\mathbf{b}^l} * D^T(\mathbf{x}), 
 \end{eqnarray}
Let $sum(A)$ be the vector whose corresponding element is the sum of $A$'s columns, there establishes that
 \begin{equation}
 \begin{split}
  \nabla\frac{\partial \mathbf{y}^{L+1}}{\partial\mathbf{w}^l} *D^T(\mathbf{x}) &= sum(diag(D(\mathbf{x}))*G(\delta^l)[1:d] )* (\mathbf{y}^{l-1})^T\\
  & + \delta^l*sum(diag(D(\mathbf{x}))*G(\mathbf{y}^{l-1})^T[1:d]),
  \end{split}
 \end{equation}
 
 \begin{equation}
 \label{eq:Neumannlast}
 \nabla\frac{\partial \mathbf{y}^{L+1}}{\partial\mathbf{b}^l} * D^T(\mathbf{x}) = sum(diag(D(\mathbf{x}))*G(\delta^l )[1:d]),
 \end{equation}
So that we can compute back propagation for Neumann boundary condition with equation \eqref{eq:grap} to \eqref{eq:Neumannlast}

\subsection{back propagation for state equation}

It is obvious that the term of time($\frac{\partial\mathbf{y}^L+1}{\partial t}$) in problem \eqref{eq:parabolic} can be computed directly by equation \ref{eq:backpro_1_w}  and \ref{eq:backpro_1_b}. The more important thing to compute is the back propagation for operator $\mathcal{L}$. There are lots of choice for $\mathcal{L}$ and in this paper we just show the gradient operator $\nabla$ and Laplace operator $\Delta$. More operator of high order can be computed similar as them.

We know that $\nabla\mathbf{y}^{L+1} = G(\mathbf{y}^{L+1})[1:d]$, so that the back propagation can be computed directly by
equation \eqref{eq:Grad_y}. Then let's consider the Laplace operator. We know that
\begin{equation}
\label{eq:Deltay}
\Delta\mathbf{y}^{L+1} = \displaystyle\sum_{i=1}^d \frac{\partial^2\mathbf{y}^{L+1}}{\partial x_i^2}
\end{equation}

Following equation \eqref{eq:Lap_y}, there establish $\Delta\mathbf{u} = sum(G^2(\mathbf{u})[1:d])$. Define $G_d(\mathbf{u}):=G^2(\mathbf{u})[1:d]$, then we have the back propagation as
\begin{eqnarray}
& &
\begin{split}
\frac{\partial \Delta\mathbf{y}^{L+1}}{\partial\mathbf{w}^l} &= \displaystyle\sum_{i = 1}^d \frac{\partial^3\mathbf{y}^{L+1}}{\partial x_i^2\partial\mathbf{w}^l}\\
& = \sum_{i=1}^d \frac{\partial^2}{\partial x_i^2}\left(\delta^l*(\mathbf{y}^{l-1})^T\right)\\
& = \left(\sum_{i=1}^d\frac{\partial^2\delta^l}{\partial x_i^2}\right)*(\mathbf{y}^{l-1})^T + 2\sum_{i=1}^d\frac{\partial\delta^l}{\partial x_i}*\frac{\partial (\mathbf{y}^{l-1})^T}{\partial x_i}+ \delta^l *  \left(\sum_{i=1}^d\frac{\partial^2(\mathbf{y}^{l-1})^T}{\partial x_i^2}\right)\\
& = sum(G_d(\delta^l))*(\mathbf{y}^{l-1})^T + 2\nabla\delta^l *\nabla(\mathbf{y}^{l-1})^T + \delta^l * sum(G_d^T(\mathbf{y}^{l-1})))
\end{split}\\
& &
\begin{split}
\frac{\partial \Delta \mathbf{y}^{L+1}}{\partial\mathbf{b}^l} &= \displaystyle\sum_{i = 1}^d \frac{\partial^3\mathbf{y}^{L+1}}{\partial x_i^2\partial\mathbf{b}^l} = sum(G_d(\delta^l)),
\end{split}\label{eq:backpro_lap}
\end{eqnarray}

With all the above equation  \eqref{eq:Deltay} to \eqref{eq:backpro_lap}, we can compute the back propagation for state equation for various operator $\mathcal{L}$.

\end{appendix}

%%%%%%%%%%%%%%%%%%%%%%%%%%%%%%%%%%%%%%%%

\section*{References}
\bibliographystyle{abbrv}
\bibliography{CauchyNet.bib}

\begin{thebibliography}{10}

\bibitem{Aarts2001}
L.~P. Aarts and P.~van~der Veer.
\newblock Neural network method for solving partial differential equations.
\newblock {\em Neural Processing Letters}, 14(3):261--271, 2001.

\bibitem{Besala1966}
D.~A. amd P.~Besala.
\newblock Uniqueness of solutions of the cauchy problem for parabolic
  equations.
\newblock {\em Journal of Mathematical Analysis and Applications}, 13:516--526,
  1966.

\bibitem{Berntsson2014}
F.~Berntsson, V.~Kozlov, L.~Mpinganzima, and B.~Turesson.
\newblock An accelerated alternating procedure for the cauchy problem for the
  helmholtz equation.
\newblock {\em Computers and Mathematics with Applications}, 68:44--60, 2014.

\bibitem{Berntsson2017}
F.~Berntsson, V.~Kozlov, L.~Mpinganzima, and B.~Turesson.
\newblock Iterative tikhonov regularization for the cauchy problem for the
  helmholtz equation.
\newblock {\em Computers and Mathematics with Applications}, 73:163--172, 2017.

\bibitem{Bourgeois_2005}
L.~Bourgeois.
\newblock A mixed formulation of quasi-reversibility to solve the cauchy
  problem for laplace's equation.
\newblock {\em Inverse Problems}, 21(3):1087--1104, 2005.

\bibitem{Cannon1967}
J.~R. Cannon and J.~Douglas, Jr.
\newblock The cauchy problem for the heat hquation.
\newblock {\em SIAM Journal on Numerical Analysis}, 4(3):317--336, 1967.

\bibitem{Giuseppe2017}
G.~Carleo and M.~Troyer.
\newblock Solving the quantum many-body problem with artificial neural
  networks.
\newblock {\em Science}, 355(6325):602--606, 2017.

\bibitem{Chakib2001}
A.~Chakib and A.~Nachaoui.
\newblock Convergence analysis for finite element approximation to an inverse
  cauchy problem.
\newblock {\em Inverse Problems}, 22(4), 2006.

\bibitem{Cheng_2014}
X.~Cheng, R.~Gong, W.~Han, and X.~Zheng.
\newblock A novel coupled complex boundary method for solving inverse source
  problems.
\newblock {\em Inverse Problems}, 30(5):055002, 2014.

\bibitem{Cimeti_re_2001}
A.~Cimeti{\`{e}}re, F.~Delvare, M.~Jaoua, and F.~Pons.
\newblock Solution of the cauchy problem using iterated tikhonov
  regularization.
\newblock {\em Inverse Problems}, 17(3):553--570, 2001.

\bibitem{Cybenko1989}
G.~Cybenko.
\newblock Approximation by superpositions of a sigmoidal function.
\newblock {\em Mathematics of Control, Signals, and Systems}, 2:303--314, 1989.

\bibitem{Elden1987}
L.~Elden.
\newblock Approximations for a cauchy problem for the heat equation.
\newblock {\em Inverse Problems}, 3(2):263--273, 1987.

\bibitem{Feng_2013}
X.~Feng and L.~Eld{\'{e}}n.
\newblock Solving a cauchy problem for a 3d elliptic {PDE} with variable
  coefficients by a quasi-boundary-value method.
\newblock {\em Inverse Problems}, 30(1):015005, 2013.

\bibitem{Hon_2001}
Y.~C. Hon and T.~Wei.
\newblock Backus-gilbert algorithm for the cauchy problem of the laplace
  equation.
\newblock {\em Inverse Problems}, 17(2):261--271, 2001.

\bibitem{Kurt1989}
K.~Hornik.
\newblock Multilayer feedforward networks are universal approximators.
\newblock {\em Neural Networks}, 2:359--366, 1989.

\bibitem{Kurt1991}
K.~Hornik.
\newblock Approximation capabilities of multilayer feedforward networks.
\newblock {\em Neural Networks}, 4:251--257, 1991.

\bibitem{Lesnic2000}
D.~N. Hào and D.~Lesnic.
\newblock The cauchy problem for laplace’s equation via the conjugate
  gradient method.
\newblock {\em IMA Journal of Applied Mathematics}, 65(2):199--217, 2000.

\bibitem{Isakov1998}
V.~Isakov.
\newblock {\em Inverse problem for partial differential equations}.
\newblock Springer-Verlag, New York, 1998.

\bibitem{Jin_2008}
B.~Jin and J.~Zou.
\newblock Augmented tikhonov regularization.
\newblock {\em Inverse Problems}, 25(2):025001, 2008.

\bibitem{Jin2010}
B.~Jin and J.~Zou.
\newblock Hierarchical bayesian inference for ill-posed problems via
  variational method.
\newblock {\em Journal of Computational Physics}, 229:7317--7343, 2010.

\bibitem{Nachaoui2002}
M.~Jourhmane and A.~Nachaoui.
\newblock Convergence of an alternating method to solve the cauchy problem for
  poisson's equation.
\newblock {\em Applicable Analysis}, 81(5):1065--1083, 2002.

\bibitem{Ito2014}
I.~Kazufumi and B.~Jin.
\newblock {\em Inverse problem. Tikhonov theory and algorithm}.
\newblock Series on Applied Mathematics, vol. 22, World Scientific, Singapore,
  2014.

\bibitem{Santosa1991}
M.~Klibanov and F.~Santosa.
\newblock A computational quasi-reversibility method for cauchy problems for
  laplace’s equation.
\newblock {\em SIAM Journal on Applied Mathematics}, 51(6):1653--1675, 1991.

\bibitem{Lagaris1998}
I.~E. Lagaris, A.~Likas, and D.~I. Fotiadis.
\newblock Artificial neural networks for solving ordinary and partial
  differential equations.
\newblock {\em IEEE Transactions on Neural Networks}, 9(5):987--1000, 1998.

\bibitem{Lagaris2000}
I.~E. {Lagaris}, A.~C. {Likas}, and D.~G. {Papageorgiou}.
\newblock Neural-network methods for boundary value problems with irregular
  boundaries.
\newblock {\em IEEE Transactions on Neural Networks}, 11(5):1041--1049, 2000.

\bibitem{Li2018}
Q.~Li, L.~Chen, C.~Tai, and W.~E.
\newblock Maximum principle based algorithms for deep learning.
\newblock {\em Journal of Machine Learning Research}, 18(165):1--29, 2018.

\bibitem{Dong2019}
Z.~Long, Y.~Lu, and B.~Dong.
\newblock Pde-net 2.0: Learning pdes from data with a numeric-symbolic hybrid
  deep network.
\newblock {\em Journal of Computational Physics}, 399:108925, 2019.

\bibitem{Malek2006}
A.~Malek and R.~S. Beidokhti.
\newblock Numerical solution for high order differential equations using a
  hybrid neural network—optimization method.
\newblock {\em Applied Mathematics and Computation}, 183(1):260--271, 2006.

\bibitem{Marin2003}
L.~Marin, L.~Elliott, P.~Heggs, D.~Ingham, D.~Lesnic, and X.~Wen.
\newblock An alternating iterative algorithm for the cauchy problem associated
  to the helmholtz equation.
\newblock {\em Computer Methods in Applied Mechanics and Engineering},
  192(5-6):709--722, 2003.

\bibitem{Mishra2018}
S.~Mishra.
\newblock A machine learning framework for data driven acceleration of
  computations of differential equations.
\newblock {\em Mathematics in Engineering}, 1(1):118--146, 2018.

\bibitem{Nachaoui2004}
A.~Nachaoui.
\newblock Numerical linear algebra for reconstruction inverse problems.
\newblock {\em Journal of Computational and Applied Mathematics},
  162(1):147--164, 2004.

\bibitem{Pang2019}
G.~Pang, L.~Lu, and G.~E. Karniadakis.
\newblock fpinns: Fractional physics-informed neural networks.
\newblock {\em SIAM Journal on Scientific Computing}, 41(4):A2603--A2626, 2019.

\bibitem{Qin2010}
H.~Qin and T.~Wei.
\newblock Two regularization methods for the cauchy problems of the helmholtz
  equation.
\newblock {\em Applied Mathematical Modelling}, 34(4):947--967, 2010.

\bibitem{Xiu2019}
T.~Qin, K.~Wu, and D.~Xiu.
\newblock Data driven governing equations approximation using deep neural
  networks.
\newblock {\em Journal of Computational Physics}, 395:620--635, 2019.

\bibitem{Raissi2019}
M.~Raissia, P.~Perdikarisb, and G.~Karniadakis.
\newblock Physics-informed neural networks: A deep learning framework for
  solving forward and inverse problems involving nonlinear partial differential
  equations.
\newblock {\em Journal of Computational Physics}, 378:686--707, 2019.

\bibitem{Reinhardt1999}
H.~Reinhardt, H.~Han, and D.~Hào.
\newblock Stability and regularization of a discrete approximation to the
  cauchy problem for laplace's equation.
\newblock {\em SIAM Journal on Numerical Analysis}, 36(3):890--905, 1999.

\bibitem{Sakamoto2011}
K.~Sakamoto and M.~Yamamoto.
\newblock Initial value/boundary value problems for fractional diffusion-wave
  equations and applications to some inverse problems.
\newblock {\em Journal of Mathematical Analysis and Applications},
  382:426--447, 2011.

\bibitem{Justin2018}
J.~Sirignano and K.~Spiliopoulos.
\newblock Dgm: A deep learning algorithm for solving partial differential
  equations.
\newblock {\em Journal of Computational Physics}, 375:1339--1364, 2018.

\bibitem{Tomoya2008}
T.~Takeuchi and M.~Yamamoto.
\newblock Tikhonov regularization by a reproducing kernel hilbert space for the
  cauchy problem for an elliptic equation.
\newblock {\em SIAM Journal on Scientific Computing}, 31(1):112--142, 2008.

\bibitem{Wei2013}
T.~Wei, Y.~Chen, and J.~Liu.
\newblock A variational-type method of fundamental solutions for a cauchy
  problem of laplace’s equation.
\newblock {\em Applied Mathematical Modelling}, 37:1039--1053, 2013.

\bibitem{White2019}
D.~A. White, W.~J. Arrighi, J.~Kudo, and S.~E. Watts.
\newblock Multiscale topology optimization using neural network surrogate
  models.
\newblock {\em Computer Methods in Applied Mechanics and Engineering},
  345(1):1118--1135, 2019.

\bibitem{Yan2019}
L.~Yan and T.~Zhou.
\newblock Adaptive multi-fidelity polynomial chaos approach to bayesian
  inference in inverse problems.
\newblock {\em Journal of Computational Physics}, 381:110--128, 2019.

\bibitem{SLi2018}
Y.S.Li and T.~Wei.
\newblock An inverse time-dependent source problem for a time–space
  fractional diffusion equation.
\newblock {\em Applied Mathematics and Computation}, 336:257--271, 2018.

\bibitem{Zhang2019}
D.~Zhang, L.~Lu, and G.~E. Karniadakis.
\newblock Quantifying total uncertainty in physics-informed neural networks for
  solving forward and inverse stochastic problems.
\newblock {\em Journal of Computational Physics}, 399:108925, 2019.

\bibitem{Zhang_2018}
Y.~Zhang, R.~Gong, X.~Cheng, and M.~Gulliksson.
\newblock A dynamical regularization algorithm for solving inverse source
  problems of elliptic partial differential equations.
\newblock {\em Inverse Problems}, 34(6):065001, 2018.

\end{thebibliography}

\end{document}